\documentclass[ps,preprint,dvips,twoside]{imsart}

\usepackage{amsthm,amsmath}
\usepackage{amsfonts}
\usepackage{amssymb}
\usepackage{epsfig}
\usepackage{verbatim}
\usepackage{graphics}
\usepackage{hyperref,verbatim}

\IfFileExists{epsf.def}{\input epsf.def}{\usepackage{epsf}}

\usepackage{mathrsfs}
\RequirePackage[numbers]{natbib}
\RequirePackage[colorlinks,citecolor=blue,urlcolor=blue]{hyperref}

\doi{10.1214/11-PS190}
\pubyear{2011}
\volume{8}
\issue{0}
\firstpage{294}
\lastpage{373}

\startlocaldefs

\numberwithin{equation}{section}
\numberwithin{figure}{section}

\newtheorem{theorem}{Theorem}[section]
\newtheorem{lemma}[theorem]{Lemma}
\newtheorem{proposition}[theorem]{Proposition}

\newtheorem{corollary}[theorem]{Corollary}

\newtheorem{assumption}[theorem]{Assumptions}
\newtheorem{problem}[theorem]{Problem}
\newtheorem{hwproblem}[theorem]{Exercise}
\theoremstyle{definition}
\newtheorem{definition}[theorem]{Definition}

\newenvironment{proofsect}[1] {\vskip0cm\noindent{\rmfamily\itshape#1.}}{\hfill\qed\vspace{0.15cm}}

\DeclareFontFamily{OT1}{eusb}{} \DeclareFontShape{OT1}{eusb}{m}{n} {<5> <6> <7> <8> <9> <10> <11> <12> <14.4> eusb10}{}
\DeclareMathAlphabet{\eusb}{OT1}{eusb}{m}{n}

\DeclareFontFamily{OT1}{eusm}{} \DeclareFontShape{OT1}{eusm}{m}{n} {<5> <6> <7> <8> <9> <10> <11> <12> <14.4> eusm10}{}
\DeclareMathAlphabet{\eusm}{OT1}{eusm}{m}{n}

\DeclareFontFamily{OT1}{eufm}{} \DeclareFontShape{OT1}{eufm}{m}{n} {<5> <6> <7> <8> <9> <10> <11> <12> <14.4> eufm10}{}
\DeclareMathAlphabet{\mathfrak}{OT1}{eufm}{m}{n}

\DeclareFontFamily{OT1}{fraktura}{}
\DeclareFontShape{OT1}{fraktura}{m}{n} {<5> <6> <7> <8> <9> <10> <11> <12> <13> <14.4> [1.1] eufm10}{}
\DeclareMathAlphabet{\fraktura}{OT1}{fraktura}{m}{n}

\DeclareFontFamily{OT1}{cmfi}{} \DeclareFontShape{OT1}{cmfi}{m}{n} {<5> <6> <7> <8> <9> <10> <11> <12> <13> <14.4> [0.9] cmfi10}{}
\DeclareMathAlphabet{\cmfi}{OT1}{cmfi}{b}{n}

\DeclareFontFamily{OT1}{cmss}{} \DeclareFontShape{OT1}{cmss}{m}{n} {<5> <6> <7> <8> <9> <10> <11> <12> <13> <14.4> cmss10}{}
\DeclareMathAlphabet{\cmss}{OT1}{cmss}{m}{n}

\newcommand{\textd}{\text{\rm d}\mkern0.5mu}
\newcommand{\texti}{\text{\rm  i}\mkern0.7mu}
\newcommand{\texte}{\text{\rm  e}\mkern0.7mu}

\newcommand{\Var}{\text{\rm Var}}
\newcommand{\Cov}{\text{\rm \Cov}}
\newcommand{\1}{{1\mkern-4.5mu\textrm{l}}}
\renewcommand{\1}{\text{\sf 1}}

\newcommand{\CC}{\mathscr C}

\newcommand{\EE}{\mathcal E}

\newcommand{\HH}{\mathcal H}

\newcommand{\NN}{\mathcal N}

\newcommand{\B}{\mathbb B}

\newcommand{\E}{\mathbb E}

\newcommand{\G}{\mathbb G}

\newcommand{\N}{\mathbb N}

\newcommand{\BbbP}{\mathbb P}
\newcommand{\Q}{\mathbb Q}
\newcommand{\R}{\mathbb R}

\newcommand{\Z}{\mathbb Z}

\newcommand{\scrA}{\mathscr{A}}

\newcommand{\scrE}{\mathscr{E}}
\newcommand{\scrF}{\mathscr{F}}
\newcommand{\scrG}{\mathscr{G}}

\newcommand{\scrL}{\mathscr{L}}

\newcommand{\scrO}{\mathscr{O}}
\newcommand{\scrP}{\mathscr{P}}

\newcommand{\scrV}{\mathscr{V}}

\newcommand{\scrX}{\mathscr{X}}

\newcommand{\twoeqref}[2]{(\ref{#1}--\ref{#2})}

\newcommand{\cc}{{\text{\rm c}}}

\def\myffrac#1#2 in #3{\raise 2.6pt\hbox{$#3 #1$}\mkern-1.5mu\raise 0.8pt\hbox{$#3/$}\mkern-1.1mu\lower 1.5pt\hbox{$#3 #2$}}
\newcommand{\ffrac}[2]{\mathchoice%
    {\myffrac{#1}{#2} in \scriptstyle}
    {\myffrac{#1}{#2} in \scriptstyle}
    {\myffrac{#1}{#2} in \scriptscriptstyle}
    {\myffrac{#1}{#2} in \scriptscriptstyle}
}

\newcommand{\hate}{\hat{\texte}}

\endlocaldefs

\begin{document}

\begin{frontmatter}

\title{Recent progress on the Random\\Conductance Model\thanksref{T1}}
\thankstext{T1}{\copyright\,\textrm{2011} \textrm{M.~Biskup}.
Reproduction, by any means, of the entire article for non-commercial purposes is permitted without charge.}
\runtitle{Random Conductance Model}

\begin{aug}
\author{\fnms{Marek} \snm{Biskup}\ead[label=e1]{biskup@math.ucla.edu}}
\address{Department of Mathematics, UCLA, Los Angeles, California, USA\\School of Economics, University of South Bohemia, \v Cesk\'e Bud\v ejovice, Czech Republic\\\printead{e1}}

\runauthor{M.~Biskup}
\end{aug}

\begin{abstract}
Recent progress on the understanding of the Random Conductance Model is
reviewed and commented. A particular emphasis is on the results on the
scaling limit of the random walk among random conductances for almost
every realization of the environment, observations on the behavior of
the effective resistance as well as the scaling limit of certain models
of gradient fields with non-convex interactions. The text is an
expanded version of the lecture notes for a course delivered at the
2011 Cornell Summer School on Probability.
\end{abstract}

\begin{keyword}[class=AMS]
\kwd[Primary ]{60K37} 
\kwd{60F17}  
\kwd[; secondary ]{60J45}  
\kwd{82B43}  
\kwd{80M40} 
\end{keyword}

\begin{keyword}
\kwd{Random conductance model} 
\kwd{elliptic environment}
\kwd{quenched invariance principle}
\kwd{corrector}
\kwd{heat kernel bounds} 
\kwd{effective resistivity} 
\kwd{homogenization}
\end{keyword}

\received{\smonth{12} \syear{2011}}

\tableofcontents

\end{frontmatter}



\section*{Prologue}
\noindent
Random walks in random environments have been at the center of the probabilists' interest for several decades. A specific class of such random walks goes under the banner of the \emph{Random Conductance Model}. What makes this class special is the fact that the corresponding Markov chains are reversible. This somewhat restrictive feature has the benefit of fruitful connections to other, seemingly unrelated fields: the random resistor networks and gradient fields. At the technical level, many of the problems are thus naturally embedded into the larger area of harmonic analysis and homogenization theory.

This survey article is an expanded version of the set of lecture notes written for a course on the Random Conductance Model that the author delivered at the 2011 Cornell Summer School on Probability. A personal point of view promoted here is that the Random Conductance Model belongs to the collection of ``paradigm'' problems such as percolation, Ising model, exclusion process, etc, that are characterized by a simple definition and yet feature interesting and non-trivial phenomena (and, of course, pose interesting questions in mathematics). The text below attempts to summarize the important developments in the understanding of the Random Conductance Model. While paying most attention to recent results, much of what is discussed draws on by-now classical work.

The text retains the layout of lecture notes that have been spiced
up with comments and references to related subjects. The general
structure is as follows: The first section introduces the three rather
different areas where the Random Conductance Model naturally appears.
Sections~\ref{sec2}--\ref{sec5} then deal predominantly with the first
such area --- namely, the various aspects of the limit behavior of
random walks in reversible random environments. Section~\ref{sec6} then
applies the introduced machinery to the remaining problems. A number of
Problems are mentioned throughout the text; these refer to questions
that are either solved directly in the text or remain a subject of
research interest until present day. Easier questions are phrased as
Exercises; these are of varied difficulty but should all be generally
accessible to graduate students.

\section*{Acknowledgments}
\noindent
This text would not exist without the generous invitation from Rick Durrett to speak at the~2011 Cornell Summer School on Probability. The author is equally grateful to Geoffrey Grimmett, who suggested rather persuasively that the preliminary and incomplete notes be made into a proper survey article --- rather than stay preliminary and incomplete forever. Much credit goes also to the coauthors N.~Berger, O.~Boukhardra, C.~Hoffman, G.~Kozma, O.~Louidor, T.~Prescott, A.~Rozinov, H.~Spohn and A.~Vanden\-berg-Rodes of various joint projects whose results are reviewed in these notes, and to numerous other colleagues for discussions that helped improve the author's understanding of the subject. T.~Kumagai was very kind to provide valuable comments on the section dealing with heat-kernel estimates, J.~Dyre suggested interesting pointers to the physics literature concerning the random resistance problem and M.~Salvi  offered a lot of feedback and suggestions on the material in Section~\ref{sec3}.
Many thanks go also to an anonymous referee for a quick and efficient report.
The research reported on in these notes has partially been supported by the NSF grants DMS-0949250 and DMS-1106850, the NSA grant~NSA-AMS 091113 and the GA\v CR project P201-11-1558.


\section{Overview and main questions}
\label{sec1}
\vspace{-6mm}
\noindent
\subsection{Random conductance model}
\noindent
We begin with the definition of the problem in the context of random walks in random environments. Consider a countable set $\scrV$ and suppose that we are given a collection of numbers $(\omega_{xy})_{x,y\in\scrV}$ with the following properties: $\omega_{xy}\ge0$~with
\begin{equation}
\label{E:pi>0}
\pi_\omega(x):=\sum_{y\in\scrV}\omega_{xy}\in(0,\infty),\qquad x\in\scrV,
\end{equation}
and the symmetry condition
\begin{equation}
\label{E:omega-sym}
\omega_{xy}=\omega_{yx},\qquad x,y\in\scrV.
\end{equation}
We will predominantly take~$\scrV$ to be the hypercubic lattice $\Z^d$ naturally embedded in~$\R^d$. The quantity $\omega_{xy}$ is called the \emph{conductance} of the pair $(x,y)$ --- the use of the term will be clarified in the subsection dealing with resistor networks.

When~$\scrV$ has an unoriented-graph structure with edge set $\scrE$, we often enforce $\omega_{xy}=0$ whenever $(x,y)\not\in\scrE$; in that case we speak of the \emph{nearest-neighbor model}. Such a model is then called \emph{uniformly elliptic} if there is $\alpha\in(0,1)$ for which
\begin{equation}
\label{E:elliptic}
\alpha<\omega_{xy}<\frac1\alpha,\qquad (x,y)\in\scrE.
\end{equation}
When $\scrV:=\Z^d$, we use the phrase ``nearest-neighbor model'' for the situation when $\scrE$ is the set of pairs of vertices that are at the Euclidean distance one from each other.

\begin{figure}[t]
\centerline{\includegraphics[width=0.86\textwidth]{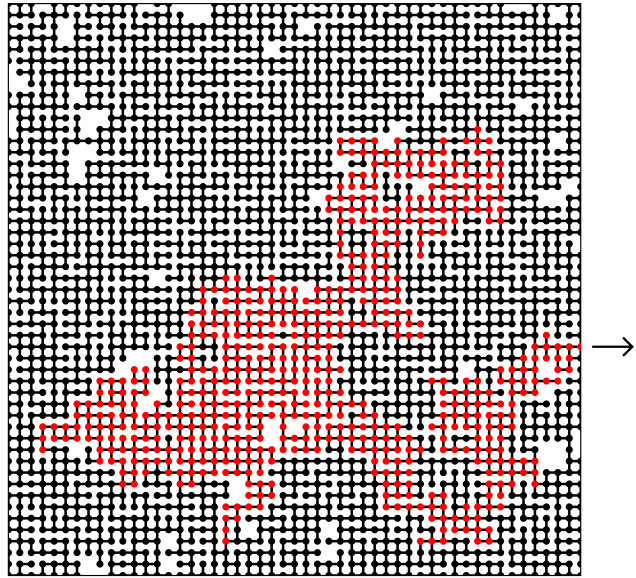}}
  \caption{Naturally included in the family of  Random Conductance Models is the simple random walk on a supercritical percolation cluster (for details on the definition and properties of percolation, see the monograph of Grimmett~\cite{Grimmett}). Here the conductances are nearest-neighbor only and take values either zero or one independently at random with the same probabilities everywhere. The density of conductance-one edges exceeds the percolation threshold so there is an infinite connected component of vertices joined by conductance-one edges; only the vertices in this component are retained in the figure. At each time the random walk chooses a neighbor of its current position at random and passes to it. The marked vertices depict those visited by a sample path of such random walk started at the center before it exits the box at the point on the right indicated by the arrow.}
  \label{fig1}
\end{figure}

The aforementioned ``random walk'' in environment~$\omega$ is technically a discrete-time Markov chain with state-space~$\scrV$ and transition kernel
\begin{equation}
\cmss P_\omega(x,y):=\frac{\omega_{xy}}{\pi_\omega(x)},\qquad x,y\in\scrV.
\end{equation}
In plain words, the ``walk'' at site~$x$ chooses its next position~$y$ proportionally to the value of the conductance~$\omega_{xy}$. The non-degeneracy condition \eqref{E:pi>0} guarantees that this chain is well defined everywhere; when positivity of $\pi_\omega$ fails at some vertices --- as, e.g., for the simple random walk on the supercritical percolation cluster, cf Fig.~\ref{fig1} --- one simply restricts the chain to the subset of~$\scrV$ where $\pi_\omega(x)>0$.

\smallskip
A key consequence of the symmetry condition \eqref{E:omega-sym} is:

\begin{lemma}
$\pi_\omega$ is a stationary and reversible measure for the Markov chain.
\end{lemma}

\begin{proofsect}{Proof}
Invoking the above definitions we get
\begin{equation}
\pi_\omega(x)\cmss P_\omega(x,y)=\omega_{xy}=\omega_{yx}=\pi_\omega(y)\cmss P_\omega(y,x),
\end{equation}
which is the condition of reversibility (a.k.a.~the detailed balance condition). The fact that $\pi_\omega$ is stationary follows by summing the extreme ends of this equality on~$x$.
\end{proofsect}

\newcounter{obrazek}

Note that for the nearest-neighbor model on~$\Z^d$ with conductances $\omega_{xy}=1$ if $|x-y|=1$ and $\omega_{xy}=0$ otherwise, the above Markov chain reduces to the ordinary simple (symmetric) random walk. In this case the increments of the walk are i.i.d.\ which permits derivation of many deep conclusions --- e.g., Donsker's Invariance Principle, Law of Iterated Logarithm, etc. However, when~$\omega$ is non-constant, the increments of the chain are no longer independent; worse yet, they are not even stationary. As we will see, this can be overcome but only at the cost of taking~$\omega$ to be a sample from a shift-invariant distribution. This reasoning underpins the large area of random walks in random environment of which the above chain is only a rather specific example.

Let $\Omega$ be the space of all configurations $(\omega_{xy})$ of the conductances. This space is naturally endowed with a product $\sigma$-algebra~$\scrF$. A \emph{shift by~$x$} is the map $\tau_x\colon\Omega\to\Omega$ acting so that
\begin{equation}
(\tau_x\omega)_{yz}:=\omega_{y+x,z+x},\qquad x,y,z\in\Z^d.
\end{equation}
We will henceforth assume that~$\BbbP$ is a probability measure on $(\Omega,\scrF)$ which is \emph{translation invariant} in the sense that
\begin{equation}
\BbbP\circ\tau_x^{-1}=\BbbP,\qquad x\in\Z^d.
\end{equation}
We recall that this measure is said to be \emph{ergodic} if $\BbbP(A)\in\{0,1\}$ for any event~$A$ with the property $\tau_x^{-1}(A)=A$ for all~$x\in\Z^d$. A canonical example of an ergodic~$\BbbP$ would be the nearest-neighbor model where the values of conductances are chosen independently at random from the same distribution. We will use $\E$ to denote expectation with respect to~$\BbbP$.

\smallskip
Let us now turn to the main questions one may wish to ask concerning the above setup. For this let~$X=(X_n)$ denote a sample path of the above Markov chain and let $P^x_\omega$ denote the law of~$X$ subject to the initial condition
\begin{equation}
P_\omega^x(X_0=x):=1.
\end{equation}
Let $\cmss P_\omega^n$ denote the $n$-th power of the transition kernel~$\cmss P_\omega$, i.e.,
\begin{equation}
\cmss P_\omega^n(x,y)=P_\omega^x(X_n=y).
\end{equation}
The aforementioned connection with the special case of simple symmetric random walk leads to the following questions:

\begin{problem}
\label{P:SLLN}
Does the limit
\begin{equation}
\lim_{n\to\infty}\frac{X_n}n
\end{equation}
exist almost surely? Under what conditions is it zero (as it is for the simple random~walk)?
\end{problem}

\begin{problem}
\label{P:CLT}
Under what conditions does the path obey an invariance principle --- i.e., does its law tend to Brownian motion under diffusive scaling of space and time? And if so, what is the rate of convergence?
\end{problem}

\begin{problem}
Does one have a local CLT as well in the sense that
\begin{equation}
\cmss P_\omega^n(x,y)\approx\frac{c_1}{n^{d/2}}\texte^{-c_2|x-y|^2/n}
\end{equation}
whenever $y$ can be ``comfortably'' reached by the random walk from~$x$ in~$n$ steps?
\end{problem}

As it turns out, there are subtle but important  differences in the precise technical {sense} in which these asymptotic statements might be true, or at least provably true. Indeed, there are two natural laws on the path space that are considered in the literature: the aforementioned \emph{quenched} law $P_\omega^x(-)$ and the \emph{annealed} or, more accurately, \emph{averaged} law $E_\Q P_\omega^x(-)$ where~$\Q$ is a specific (natural) measure on environments (similar to~$\BbbP$). An advantage of the annealed law is that, thanks to averaging, it allows for an easier control of the irregularities of the environment; a drawback is that the path law under it is no longer Markovian.
As we will see, one of the main challenges for the Random Conductance Model that prevail to the present day is the
resolution~of:

\begin{problem}
Does the annealed invariance principle imply the quenched invariance principle? (Here and henceforth the words annealed and quenched designate the path distribution that is considered for the scaling limit.)
\end{problem}

\begin{figure}[t]
\centerline{\includegraphics[width=0.6\textwidth]{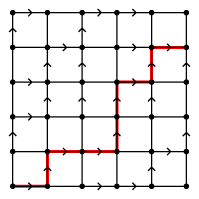}}
  \caption{An example of a random walk in a random environment where, at each vertex, one of the North or East arrows is chosen independently at random (with equal probabilities). The random walk is then forced to follow the arrows. The extremity of this example is seen from the fact that while the quenched law of the path is deterministic --- and no invariance principle can hold for fluctuations --- the  averaged law looks like an ordinary North\,\&\,East random walk whose fluctuations are described by the Central Limit Theorem.}
  \label{fig1b}
\end{figure}

We remark that, for general random walks in random environments, the annealed and quenched law can be dramatically different. See~Fig.~\ref{fig1b}.

As soon as the above ``fundamental'' questions have been resolved, one can try to imitate various derivations that have over years been accomplished in the context of the simple random walk. This leads to further rather interesting questions, for~instance:

\begin{problem}
What are the intersection exponents --- i.e., the decay exponents for the probability of non-interection up to the first-exit time from a ball of a large radius --- of several independent copies of such random walks?
\end{problem}

\begin{problem}
Does the (chronological) loop-erasure of the walk have the same scaling limit as the simple random walk?
And how many steps of the walk are needed to generate~$n$ steps of the loop-erased walk?\end{problem}

\begin{problem}
Is there a scaling limit for the trace of the walk in Fig.~\ref{fig1} as the size of the box tends to infinity?
\end{problem}

The last question naturally puts us into a bounded domain where, as it turns out, many additional technical difficulties arise compared to the full lattice. However, even the following questions are quite relevant:

\begin{problem}
\label{pb1.8}
Is there a scaling limit of the random walk among random conductances restricted to the half-space, quarter space or a wedge (i.e., for the problem with conductances ``leading'' outside these regions set to zero)?
\end{problem}

\subsection{Digression on continuous time}
\noindent
Although the discrete-time Markov chain is very natural, one is often interested in a \emph{continuous-time} version thereof. We will therefore introduce these objects right away and discuss some of the technical issues that come up in this context.

There are two natural ways how to make the time flow continuously. First, we may simply Poissonize the discrete time and consider the transition kernel
\begin{equation}
\cmss Q_\omega^t(x,y):=\sum_{n\ge0}\frac{t^n}{n!}\texte^{-t}\,\cmss P_\omega^n(x,y).
\end{equation}
The corresponding (continuous-time) Markov process is then referred to as \emph{constant-speed} random walk  among random conductances (CSRW), where the adjective highlights the fact that the jumps happen at the same rate regardless of the current
position.\looseness=1

Another natural way how to make time flow continuously is by attaching a clock to each pair $(x,y)$ that rings after exponential waiting times with expectation $1/\omega_{xy}$. This can just as well be done by prescribing the generator
\begin{equation}
\label{L-omega}
(\cmss L_\omega f)(x):=\sum_y\omega_{xy}\bigl[f(y)-f(x)\bigr],
\end{equation}
and demanding that the corresponding transition kernel $\cmss R_\omega^t$ is the (unique) stochastic solution of the backward Kolmogorov equations,
\begin{equation}
\label{E:backward-Kolmogorov}
\frac{\textd}{\textd t}\cmss R_\omega^t(x,y)=\sum_z\cmss L_\omega(x,z)\cmss R_\omega^t(z,y)
\end{equation}
with initial condition
\begin{equation}
\label{E:initial-cond}
\cmss R_\omega^t(x,y)=\delta_x(y).
\end{equation}
Here $\delta_x(z)$ equals one when~$x=z$ and zero otherwise. This leads to the \emph{variable speed} random walk among random conductances (VSRW), because the resulting Markov chain at~$x$ makes a new jump at rate~$\pi_\omega(x)$.

A specific problem with the VSRW is that the walk may escape to infinity in finite time --- a \emph{blow-up} occurs. (This will not happen for the discrete-time walk and thus also the CSRW.) A simple criterion to check is:

\begin{hwproblem}
\label{ex1.11}
Consider a configuration $\omega$ of conductances such that $\pi_\omega(x)\in(0,\infty)$ for each~$x$. Let~$(X_k)$ be the path of the discrete-time random walk among conductances~$\omega$ and let $T_0,T_1,\dots$ be the times between the successive jumps of the corresponding VSRW. Show that
\begin{equation}
P_\omega^0\biggl(\,\sum_{k=0}^\infty T_k<\infty\biggr)=P_\omega^0\biggl(\,\sum_{k=0}^\infty \frac1{\pi_\omega(X_k)}<\infty\biggr)
\end{equation}
\end{hwproblem}

The upshot of this Exercise is that the question of blow-ups in VSRW can be resolved purely in the context of the discrete-time walk. We refer to, e.g., Liggett~\cite[Chapter~2]{Liggett} for a thorough discussion of such situations. See also Exercise~\ref{ex2.8} in Sect.~\ref{sec2.2}.

The above transition kernels are distinguished by their invariant measures and natural function spaces they act on. Indeed, we can write~$\cmss R_\omega^t$ as
\begin{equation}
\label{VS}
\cmss R_\omega^t(x,y):=\langle \delta_y,\texte^{t\,\cmss L_\omega}\delta_x\rangle_{\ell^2(\Z^d)},
\end{equation}
where we think of~$\ell^2(\Z^d)$ as endowed by the counting measure. On the other hand, the constant speed Markov chain admits the representation
\begin{equation}
\cmss Q^t_\omega(x,y):=\frac1{\pi_\omega(x)}\langle\delta_y,\texte^{t(\cmss P_\omega-1)}\delta_x\rangle_{\ell^2(\pi_\omega)}
\end{equation}
where $\ell^2(\pi_\omega)$ is the space of functions $f\colon\Z^d\to\R$ that are square integrable with respect to the measure $\pi_\omega$ on~$\Z^d$. In this case the generator of the Markov chain is simply $\cmss P_\omega-1$. The reason why one uses different underlying measure on~$\Z^d$ in the two cases is seen via:

\begin{hwproblem}
Show that $\cmss L_\omega$ is symmetric on $\ell^2(\Z^d)$ while $\cmss P_\omega-1$ is symmetric on~$\ell^2(\pi_\omega)$. In particular, the VSRW is reversible with respect to the counting measure on~$\Z^d$ while the CSRW is reversible with respect to~$\pi_\omega$.
\end{hwproblem}

It is clear that the constant-speed chain will follow the discrete-time chain very closely, but the variable-speed chain may deviate considerably because its time parametrization depends on the entire path. This discrepancy will be particularly obvious in the places where, in comparison with the neighbors, $\pi_\omega(x)$ is either very small (VSRW gets stuck but CSRW departs easily) or very large (VSRW departs easily but CSRW gets stuck). This may or may not be a disadvantage depending on the context.

\subsection{Harmonic analysis and resistor networks}
\label{sec1.2}\noindent
The above (discrete-time) Markov chain is in a class of models for which we can apply a well-known connection between reversible Markov processes and harmonic analysis/electrostatic theory. This connection goes back to the work of Kirchhoff in mid 1800s (Kirchhoff~\cite{Kirchhoff}) and it underlies many modern treatments of Markov processes. For our purposes the best general introductory text seems to be the monograph by Doyle and Snell~\cite{Doyle-Snell}.

We begin by introducing some relevant notions for the full lattice; the finite-volume counterparts will be dealt with later. For a configuration of the conductances $(\omega_{xy})$ and a function $f\colon\Z^d\to\R$, let us define
\begin{equation}
\label{E:Dir-E}
\EE(f):=\frac12\sum_{x,y}\omega_{xy}\bigl[f(y)-f(x)\bigr]^2.
\end{equation}
In physics vernacular, this is the electrostatic or Dirichlet energy corresponding to the electrostatic potential~$f$. We then define the \emph{effective (point-to-point) resistance} $R(x,y)$ between~$x$ and~$y$ by the formula
\begin{equation}
\label{E:1.20}
R(x,y)^{-1}:=\inf\bigl\{\EE(f)\colon f\in\ell^2(\pi_\omega),\, f(x)=1,\, f(y)=0\bigr\}.
\end{equation}
More generally, we define an \emph{effective point-to-set resistance} $R(x,A)$ by requiring $f(y)=0$ for all~$y\in A$ in the formula above. Of course, both $\EE(f)$ and $R(x,y)$ depend on~$\omega$, but we leave that notationally implicit.

A key problem now is a computation, an analysis of various scaling properties, of the effective resistance. As a warm-up, consider now the homogeneous problem when the conductances are equal to one for nearest neighbors and zero otherwise. Leaving aside some technical issues, any minimizer of the Dirichlet energy in \eqref{E:1.20} will then~obey
\begin{equation}
\label{E:harmonic}
\sum_{u\colon|u-v|=1} \bigl[f(u)-f(v)\bigr]=0,\qquad v\in\Z^d\setminus\{x,y\},
\end{equation}
with
\begin{equation}
f(y)-f(x)=1.
\end{equation}
In other words, $f$ is \emph{discrete harmonic} everywhere away from~$x$ and~$y$. It is an interesting exercise in upper-division analysis to solve:

\begin{hwproblem}
Fix~$I\in\R$. For the homogeneous nearest-neighbor problem, use Fourier transform to solve the equation
\begin{equation}
\sum_{u\colon|u-v|=1} \bigl[f(u)-f(v)\bigr]=I\bigl[\delta_x(v)-\delta_y(v)\bigr],\qquad v\in\Z^d,
\end{equation}
and then adjust~$I$ so that $f(x)-f(y)=1$. Use this to derive an integral formula for~$R(x,y)$.
\end{hwproblem}

We can thus check that while the following problem may appear hard, it is at least not ill posed:

\begin{hwproblem}
For the homogeneous nearest-neighbor problem on~$\Z^2$, show without relying on Fourier transform that $R(x,y)=\ffrac12$ whenever $x$ and $y$ are nearest neighbors.
\end{hwproblem}

\begin{figure}[t]
\centerline{\includegraphics[width=3.5in]{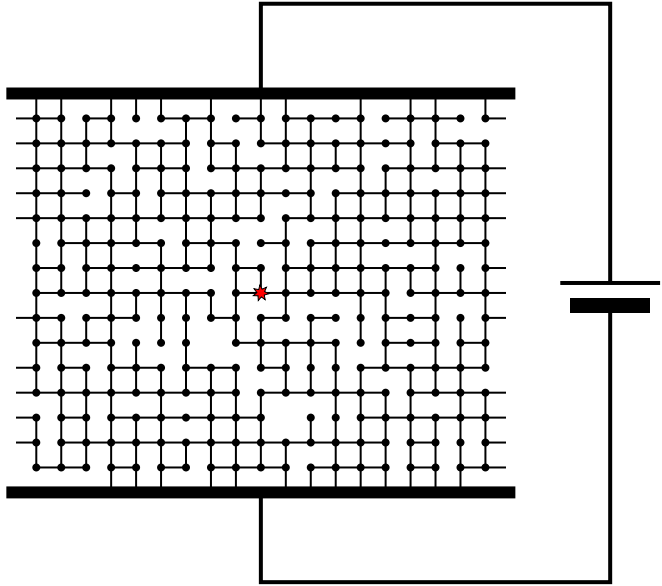}}
\caption{An example of an electrostatic problem connected to the Random Conductance Model. Here part of the percolation cluster in a slab with vertical coordinates in the interval $[-N,N]$ is attached to metal plates with a given voltage difference. The edges present in the cluster have resistivity one, the edges that are absent are total insulators. A key question is to find the total current density --- per unit area of the plates --- running through the system. Another question is the value of the electrostatic potential at the origin.}
\label{fig2}
\end{figure}

Returning to the full-fledged Random Conductance Model, let us now discuss the (somewhat degenerate) example of the supercritical percolation cluster depicted in Fig.~\ref{fig2}. Assuming the potential is fixed to $\varphi\equiv-1$ at a conducting plate at ``height'' $-N$ and to~$\varphi\equiv+1$ at the corresponding  plate at ``height'' $+N$, the question is what is the electrostatic potential right at the center. As before, this potential is a minimizer of the Dirichlet energy $\EE(f)$ in \eqref{E:Dir-E} subject to the conditions that $f(x):=1$ when $\hate_2\cdot x\ge N$ and $f(x):=-1$ when $\hate_2\cdot x\le -N$. Here $\hate_i$ is the coordinate unit vector in the $i$-th lattice direction.

What makes this problem relevant for probabilists is the existence of a direct probabilistic ``solution:'' Let $\tau_\pm^{(N)}$ be the first hitting time of the upper, resp., lower metal plate,
\begin{equation}
\tau_{\pm}^{(N)}:=\inf\{n\ge0\colon X_n\cdot\hate_2=\pm N\}.
\end{equation}
Then the electric potential at vertex~$x$ turns out to be given by the formula
\begin{equation}
\label{E:1.13}
\varphi(x):=P^x_\omega\bigl(\tau_+^{(N)}<\tau_-^{(N)}\bigr)-P^x_\omega\bigl(\tau_-^{(N)}<\tau_+^{(N)}\bigr),
\end{equation}
where $P_\omega^x$ is our notation for the law on paths $(X_n)$ of the random walk on environment~$\omega$ such that
$P_\omega^x(X_0=x)=1$. The key point is that the function $\varphi$ defined by \eqref{E:1.13} is \emph{harmonic} with respect to the generator of the continuous time Markov chain \eqref{L-omega} with the boundary values given as above. Here a function is said to be harmonic at~$x$ when $\cmss L_\omega f(x)=0$.

\begin{hwproblem}
Prove the formula \eqref{E:1.13} by showing that such a harmonic function is uniquely determined by its boundary data.
\end{hwproblem}

Notice that, as soon as the conductances are non-constant, there is no reason why the potential~$\varphi$ at the symmetry point should be equal to zero --- as it would be, thanks to symmetry considerations, for the case of homogeneous networks. Obviously, this is quite related to Problem~\ref{pb1.8}.

\smallskip
The concept of effective resistance is closely related to the question of \emph{recurrence} and \emph{transience} of the corresponding Markov chain. Let
\begin{equation}
\label{E:1.16a}
\tilde\tau_x:=\inf\{n\ge1: X_n=x\}
\end{equation}
and, for~$A\subset\Z^d$,
\begin{equation}
\label{E:1.17a}
\tau_A:=\inf\{n\ge0\colon X_n\in A\}.
\end{equation}
Set $\Lambda_N:=[-N,N]^d\cap\Z^d$. The chain will then be recurrent if $P_\omega^0(\tilde\tau_0<\tau_{\Lambda_N^\cc})\to1$ as~$N\to\infty$ and transient otherwise. The connection with effective resistance shows that the tendency to recurrence decreases with increasing conductances. Explicitly, we have:

\begin{hwproblem}
Show that the function
\begin{equation}
\varphi(x):=\begin{cases}
P_\omega^x(\tilde\tau_0<\tau_{\Lambda_N^\cc}),\qquad&\text{if }x\ne0,
\\
1,\qquad&\text{if }x=0.
\end{cases}
\end{equation}
is the unique minimizer of the Dirichlet energy for the boundary conditions corresponding to point-to-set resistance $R(x,\Lambda_N^\cc)$ and use this to derive
\begin{equation}
R(x,\Lambda_N^\cc)^{-1}=\pi_\omega(0)P_\omega^0(\tilde\tau_0\ge\tau_{\Lambda_N^\cc}).
\end{equation}
Conclude that $\pi_\omega(0)P_\omega^0(\tilde\tau_0\ge\tau_{\Lambda_N^\cc})$ is monotone increasing in each~$\omega_{xy}$.
\end{hwproblem}

The upshot of this observation is that if $\omega_{xy}\le\omega'_{xy}$ for all pairs~$x,y$, then
\begin{equation}
\pi_\omega(0)P^0_\omega(\tilde\tau_0\ge\tau_{\Lambda_N^\cc})\le \pi_{\omega'}(0)P^0_{\omega'}(\tilde\tau_0\ge\tau_{\Lambda_N^\cc}).
\end{equation}
In particular, if the random walk is recurrent in the environment~$\omega'$ then so it is in $\omega$, and \emph{vice versa} for the question of transience. For (say) nearest-neighbor Random Conductance Models subject to the ellipticity condition \eqref{E:elliptic}, recurrence is thus equivalent to the recurrence of the simple random walk. However, as soon as ellipticity is violated, interesting problems arise.

Consider for illustration the random walk on the supercritical percolation cluster. There the conductances are bounded above but not below. This still permits us to conclude that the random walk is is recurrent in spatial dimension $d=2$, and if it is transient in dimension $d=3$, then it is transient in all dimensions $d\ge3$. A key question to resolve is thus:

\begin{problem}
\label{pb1.15}
Is the random walk on almost every realization of the three-dimensional supercritical percolation cluster transient?
\end{problem}

The following question should ideally be solved before tackling Problem~\ref{pb1.15}:

\begin{problem}
Let $\omega_b\in\{0,1\}$ and let $\CC_\infty(\omega)$ denote the set of vertices in~$\Z^d$ that lie in an infinite self-avoiding path using only edges with~$\omega_b=1$. Let $\omega'$ differ from~$\omega$ in a finite number of coordinates so that $\omega'_b\ge\omega_b$ for all~$b$. Assuming that~$0\in\CC_\infty(\omega)$, show that
\begin{equation}
P_\omega^0(X \text{ \rm is transient})=P_{\omega'}^0(X \text{ \rm is transient})
\end{equation}
and conclude that $\{X \text{ \rm is transient on }\CC_\infty\}$ is a tail event. (In particular, for Bernoulli $\omega_b$'s, it is also a zero-one event.)
\end{problem}

There are a good number of variations on the  problem depicted in Fig.~\ref{fig2}, but here is one that has been particularly perplexing for a number of years --- in spite of an existing solution claimed in the book of Jikov, Kozlov and Oleinik~\cite{JKO}. The formulation goes back to Kesten's monograph on percolation (Kesten~\cite{Kesten}). Consider the square box $\Lambda_N:=[-N,N]^2\cap\Z^d$ and let $\scrG_N$ be the set of those edges whose both endpoints lie in the infinite bond-percolation cluster and also in~$\Lambda_N$. Define the effective resistance
\begin{equation}
R_N^{-1}:=\inf\biggl\{\sum_{(x,y)\in\scrG_N}\bigl[f(x)-f(y)\bigr]^2\colon f(x)=-\hate_1\cdot x\,\,\,\text{when}\,\,\,\hate_1\cdot x=\pm N\biggr\}
\end{equation}
corresponding to the boundary conditions $-N$ on the ``left'' side of the box and $+N$ on the ``right'' side of the boundary; no boundary condition is prescribed at the remaining portions of the boundary. It is not hard to convince oneself that $R_N^{-1}$ is at most of order~$N^d$, but identifying a precise rate is far more challenging:

\begin{problem}
\label{pb1.18a}
Prove that for almost every realization of the supercritical percolation cluster, the limit
\begin{equation}
\lim_{N\to\infty}\frac{R_N^{-1}}{N^d}
\end{equation}
exists and is independent of the realization. Characterize its value.
\end{problem}

Of course, once this has been settled, one may want to go beyond a LLN-type of information and study the fluctuations. Interestingly, as observed already a while ago by Wehr~\cite{Wehr}, the variance of~$R_N^{-1}$ is order at most~$N^d$ --- at least in the elliptic setting --- which suggests the following question:

\begin{problem}
Show that the law of $N^{-d/2}[\,R_N^{-1}-\E R_N^{-1}]$ tends to Gaussian as~$N\to\infty$.
\end{problem}

Recently, thanks to the work of Gloria and Otto~\cite{Gloria-Otto}, we even know that the variance is actually of order~$N^d$ (at least in $d\ge3$) so the time seems ripe for resolving this problem as well.

\subsection{Gradient models}
\noindent
The third and somewhat unexpected context in which one naturally encounters the Random Conductance Model is that of gradient fields. In our formulation, a \emph{gradient field} is a collection of $\R$-valued random variables~$\phi_x$ indexed by the vertices $x\in\Z^d$. We impose the following law:
\begin{equation}
\label{E:GGM}
\mu_\Lambda^{\bar\phi}(\textd\phi):=\frac1{Z^{\bar\phi}_\Lambda}\exp\biggl\{\,\,-\!\!\!\sum_{\langle x,y\rangle\in\B(\Lambda)}V(\phi_x-\phi_y)\biggr\}\prod_{x\in\Lambda}\textd\phi_x\,\prod_{x\not\in\Lambda}\delta_{\bar\phi_x}(\textd\phi_x).
\end{equation}
Here~$\Lambda\subset\Z^d$ is a finite set and $\B(\Lambda)$ is the set of all edges with at least one endpoint in~$\Lambda$. The function $V\colon\R\to\R$ is the \emph{potential} which we take to be a continuous, even function with sufficient (e.g., quadratic) growth at infinity. The measure depends on the values immediately outside~$\Lambda$ which are set to the \emph{boundary condition} $\bar\phi$ by the product of delta-masses.

Gradient models are ubiquitous in physical sciences where they arise as effective-interface models, with $\phi_x$ giving the height of a surface above a reference plane, or in descriptions of the fluctuation fields in critical statistical mechanical (spin) models. A higher-dimensional variant, particularly, $\phi_x\in\R^d$, has the interpretation of a deformation field representing the displacements of atoms in a crystal from their ideal positions. Further applications can be found in field theory and material physics. The reviews by Giacomin~\cite{Giacomin}, Velenik~\cite{Velenik}, Funaki~\cite{Funaki} and Sheffield~\cite{Sheffield-RS} give more information and further connections.

\smallskip
We will actually consider the measure \eqref{E:GGM} to be a law on the sigma-field  of \emph{gradient events}
\begin{equation}
\label{E:grad-sigma-alg}
\scrF:=\sigma\bigl(\phi_x-\phi_y\colon x,y\in\Z^d\bigr),
\end{equation}
which is legitimate since the corresponding restriction of $\mu_\Lambda^{\bar\phi}$ does not depend on the values of~$\bar\phi$ but only on their differences. This restriction is dictated by practical reasons --- the actual ``height'' of an interface is usually of lesser importance than the ``shape'' of its configuration --- but also due to technical restrictions in low spatial dimensions. We say that a measure $\mu$ on $(\R^{\Z^d},\scrF)$ is a \emph{gradient Gibbs measure} (GGM) if for every~$A\in\scrF$ and any finite $\Lambda\subset\Z^d$,
\begin{equation}
\mu(A)=E_\mu\bigl(\mu_\Lambda^{\bar\phi}(A)\bigr),
\end{equation}
where the expectation is over the boundary condition $\bar\phi$. Put another way, this says that the conditional probability of~$A$ given the configuration $\bar\phi$ outside~$\Lambda$ is exactly the measure \eqref{E:GGM}.

\smallskip
Before we start discussing the relevant problems arising in this subject area, it is interesting to note two special instances of the above formalism. The first one is the $d=1$ case. Let us assume that~$\Lambda$ is connected and, in fact, $\Lambda:=(-N,N)\cap\Z$. Then the law of of the gradients, $(\phi_{x+1}-\phi_x)_{x=-N}^{N-1}$ is i.i.d.\ --- with marginal law proportional to $\texte^{-V(\eta)}\textd\eta$ --- conditional on
\begin{equation}
\sum_{x=-N}^{N-1}(\phi_{x+1}-\phi_x)=\bar\phi_{N}-\bar\phi_{-N}.
\end{equation}
This situation can be analyzed with the help of standard methods of large-deviation theory (cf, e.g., Dembo and Zeitouni~\cite{Dembo-Zeitouni}, den Hollander~\cite{denHollander}) --- in fact, Cram\'er's theorem more or less suffices --- and so one can prove:

\begin{hwproblem}
\label{HW3}
Suppose $d=1$ and a linear boundary condition, i.e., $\bar\phi_x:=t x$ for some $t\in\R$. Show that, for any continuous, even potential~$V$ growing superlinearly at infinity, the law of
\begin{equation}
t\mapsto N^{-1/2}[\phi_{\lfloor tN\rfloor}-tN],\qquad -1\le t\le1,
\end{equation}
linearly interpolated into a continuous function, scales to a Brownian bridge as~$N\to\infty$. Characterize the  variance at $t=0$.
\end{hwproblem}

Another instance of special interest is that when~$V$ is quadratic,
\begin{equation}
\label{V-Gauss}
V(\eta):=\frac\kappa2\eta^2,
\end{equation}
for some \emph{stiffness} $\kappa>0$. In this case the above measure is Gaussian and so it is amenable to explicit calculations. In fact, for (say) zero boundary condition $\bar\phi_x=0$, one can even pass to the limit $\Lambda\uparrow\Z^d$, provided one restricts to the sigma-algebra of gradient events \eqref{E:grad-sigma-alg}.
This restriction is necessary because in dimensions~$d=1,2$, the law of~$\phi_0$ is not tight in this limit. To see this in more explicit terms, note that
\begin{equation}
\text{Cov}_{\mu_\Lambda^0}(\phi_x,\phi_y)=\cmss G_\Lambda(x,y),
\end{equation}
where $\cmss G_\Lambda(x,y)$ is the Green's function associated with the discrete Laplacian with Dirichlet boundary condition on~$\partial\Lambda$. In probabilist's terms, $\cmss G_\Lambda(x,y)$ is the expected number of visits to~$y$ by the simple random walk started at~$x$ before it exits from~$\Lambda$. The classical formula
\begin{equation}
\label{E:1.38}
\cmss G_\Lambda(x,y)=\frac{P_\omega^x(\tilde\tau_y<\tau_{\Lambda^\cc})}{1-P_\omega^x(\tilde\tau_x<\tau_{\Lambda^\cc})},
\end{equation}
see, e.g., Spitzer~\cite{Spitzer} or Lawler~\cite{Lawler-book}, using the notation \twoeqref{E:1.16a}{E:1.17a}, provides an explicit connection to the issues discussed in the previous subsection.

\smallskip
An analogue of Exercise~\ref{HW3} in $d\ge2$ will then be:

\begin{hwproblem}
\label{HW4}
Consider the Gaussian gradient model with \eqref{V-Gauss} with~$\kappa:=1$. For a sample of the field $(\phi_x)$ from the infinite-volume limit $\mu:=\lim_{\Lambda\uparrow\Z^d}\mu_\Lambda^0$, and a smooth $f\colon\R^d\to\R$ with compact support and $\int f(x)\textd x=0$ define
\begin{equation}
\label{E:phi_epsilon}
\phi_\epsilon(f):=\epsilon^{1+d/2}\int \phi_{\lfloor x\rfloor}\,f(\epsilon x)\,\textd x.
\end{equation}
Show that in the limit $\epsilon\downarrow0$, the law of~$\phi_\epsilon(f)$ is a Gaussian $\NN(0,\sigma_f^2)$, where
\begin{equation}
\sigma_f^2:=(f,-\Delta^{-1}f)_{L^2(\R^d)}=\int\textd x\,\textd y\, f(x) G(x-y)f(y),
\end{equation}
where~$G:=\lim_{\Lambda\uparrow\Z^d}G_\Lambda$ is the infinite-volume Green's function. (The expression on the right is well-define because $f\in\text{\rm Dom}(\Delta^{-1})$.)
\end{hwproblem}

The problem is meaningful in all $d\ge1$ but only in $d=1$ we have a hope to describe the limit as a (real-valued) process. This is because the limiting continuum object, the \emph{Gaussian Free Field} (GFF), is very rough in $d\ge2$ and, in fact, can only be interpreted in the sense of distribution theory --- hence our formulation using a linear functional $\phi_\epsilon$ in \eqref{E:phi_epsilon}. We refer to, e.g., Sheffield~\cite{Sheffield-GFF} for more information on the tightness issues and other aspects of the GFF.

\smallskip
Having dealt with these instructive examples, let us move on to general potentials~$V$. A remarkable feature of gradient models is that much of what has already been said about the quadratic case applies to any gradient model for which~$V$ is \emph{uniformly strictly convex} --- i.e., when~$V$ is~$C^2$ with~$V''$ positive and uniformly bounded away from zero and infinity. (We will expound on the specifics in the discussion of dynamical environments in Section~\ref{sec4.4}). Unfortunately, convex potentials are not what one typically finds in models coming from realistic systems and/or applications and so  the last decade has witnessed a major push to obtain a similar level of control also for non-convex interactions. This has so far succeeded only partially because most of the existing techniques fail as soon as~$V$ is non-convex anywhere, regardless how unlikely (or energetically unfavorable) a configuration for which this happens may be.

\begin{figure}[t]
\centerline{\includegraphics[width=0.6\textwidth]{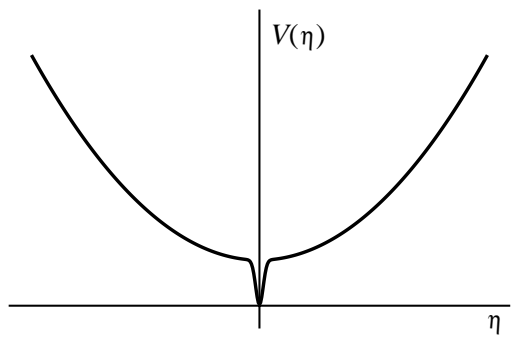}}
\caption{A plot of the potential of the form in \eqref{V-def} with $\rho:=p\delta_{\kappa_1}+(1-p)\delta_{\kappa_2}$. Once $0<\kappa_2\ll\kappa_1$ the potential is not convex.
\normalsize}
\label{fig3}
\end{figure}

Notwithstanding, there is a family of models with non-convex~$V$ that can be studied by way of a connection to the Random Conductance Model. These models are defined generally by requiring that~$V$ be given by
\begin{equation}
\label{V-def}
\texte^{-V(\eta)}:=\int_{(0,\infty)}\rho(\textd\kappa)\,\texte^{-\frac12\kappa\eta^2},
\end{equation}
where~$\rho$ is a positive measure on positive reals. Notice that when $\rho$ is supported at a single point, then~$V$ is quadratic, but as soon as $\rho$ has at least two points in its support, $V$ can be non-convex, see~Fig.~\ref{fig3}. (Nontheless~$\eta\mapsto V(\eta)$ will always be increasing on positive~$\eta$'s.)

An essential feature of the assumption \eqref{V-def} is that it permits us to consider $\mu_\Lambda^{\bar\phi}(\textd\phi)$ in \eqref{E:GGM} as the $\phi$-marginal of the measure $\mu_\Lambda^{\bar\phi}(\textd\phi\textd\kappa)$ on $\R^{\Z^d}\times(0,\infty)^{\B(\Lambda)}$ which is given by
\begin{equation}
\label{E:1.43}
\mu_\Lambda^{\bar\phi}(\textd\phi\textd\kappa):=\frac1{Z^{\bar\phi}_\Lambda}\,\texte^{-H_\Lambda(\phi,\kappa)}\,\prod_{x\in\Lambda}\textd\phi_x\,\prod_{x\not\in\Lambda}\delta_{\bar\phi_x}(\textd\phi_x)\!\!\!\prod_{\langle x,y\rangle\in\B(\Lambda)}\!\!\!\rho(\textd\kappa_{xy}),
\end{equation}
where
\begin{equation}
H_\Lambda(\phi,\kappa):=\frac12\sum_{\langle x,y\rangle\in\B(\Lambda)}\!\!\!\kappa_{xy}(\phi_x-\phi_y)^2.
\end{equation}
To see why this holds, introduce a ``private'' variable $\kappa_{xy}=\kappa_{yx}$ for each edge $\langle x,y\rangle\in\B(\Lambda)$ and use the additive structure of the interaction to write the exponential weight
in \eqref{E:GGM} as the exponential weight in \eqref{E:1.43} integrated over the product of the $\rho$'s. A key point is that, by regarding the $\kappa_{xy}$'s as genuine random variables and conditioning on their values, the law of the~$\phi$'s is \emph{again Gaussian}, albeit now with a spatially inhomogeneous covariance structure.

The above constructions can be performed in infinite volume; see Biskup and Spohn~\cite{Biskup-Spohn} for details. We will only communicate the salient conclusions: First, one can represent every gradient measure $\mu$ for the potential~$V$ in \eqref{V-def} as the $\phi$-marginal of an \emph{extended} measure~$\nu$ on pairs of configurations $(\phi,\kappa)$ such that the following holds:
\begin{enumerate}
\item[(1)]
Conditional on the~$\phi$'s, the individual $\kappa$'s are independent with $\kappa_{xy}$ having the marginal law proportional to $\texte^{-\frac12\kappa_{xy}(\phi_x-\phi_y)^2}\rho(\textd\kappa_{xy})$.
\item[(2)]
Conditional on the $\kappa$'s, the $\phi$'s are then Gaussian with covariance given by the inverse of (the negative of) the generator
\begin{equation}
\cmss L_\kappa f(x):=\sum_{y\colon |y-x|=1}\kappa_{xy}\bigl[f(y)-f(x)\bigr]
\end{equation}
of the Random Conductance Model with nearest-neighbor conductances $(\kappa_{xy})$. (The mean can be characterized too, but we will discuss this in the proof of Theorem~\ref{thm-main}.)
\item[(3)]
The $\kappa$-marginal is generally strongly correlated, but if the initial gradient measure is ergodic with respect to translations, then the extended is ergodic as well.
\end{enumerate}
For those familiar with the Random Cluster Model (see, e.g., the monograph by Grimmett~\cite{Grimmett-RCM}) and the Fortuin-Kasteleyn represenation of the Potts model (Fortuin and Kasteleyn~\cite{FK}), the above should be quite reminiscent of the so called Edwards-Sokal coupling of these two processes (Edwards and Sokal~\cite{ES}).

The structure described above offers the possibility to study the gradient model with non-convex interaction of the type \eqref{V-def} by conditioning on the $\kappa$'s. The proof of scaling of the gradient field to the Gaussian Free Field at large scales then boils down to solving:

\begin{problem}
\label{pb19}
Let~$(\phi_x)$ be a collection of Gaussian fields with mean zero and covariance given for any $g\colon\Z^d\to\R$ with finite support and $\sum_xg(x)=0$ by
\begin{equation}
\text{\rm Var}\biggl(\,\sum_x g(x)\phi_x\biggr)=\sum_{x,y\in\Z^d}g(x)g(y)(-\cmss L_\kappa)^{-1}(x,y),
\end{equation}
where $(-\cmss L_\kappa)^{-1}(x,y)$ --- the inverse of the operator $-\cmss L_\kappa$ --- can equivalently be described as the full-lattice Green's function of the random walk among nearest-neighbor random conductances~$\kappa$. Show that, for any ergodic law $\BbbP$ on the~$\kappa$'s, the random functional $\phi_\epsilon(f)$ in \eqref{E:phi_epsilon} tends to a Gaussian random variable $\BbbP$-a.s. Characterize its variance.
\end{problem}

 As we will see this will become even more interesting once we start discussing gradient fields with non-vanishing tilt. Naturally, once these basic convergence issues are settled one can turn to more subtle questions such as, for instance:

\begin{problem}
For a Gaussian field with covariance $(-\cmss L_\kappa)^{-1}$ with zero values on the boundary of a cubic domain $\Lambda_N:=[-N,N]^d\cap\Z^d$, what is the distribution of $\max_{x\in\Lambda_N}\phi_x$? What is the scaling limit of the level sets? And how about the Hausdorff dimension of various exceptional sets (e.g., the so called thick points)?
\end{problem}

These problems have recently been studied for the homogeneous lattice GFF, e.g., by Bolt\-hau\-sen, Deuschel and Zeitouni~\cite{Bolthausen-Zeitouni}, Daviaud~\cite{Daviaud}, Hu, Miller and Peres~\cite{Hu-Miller-Peres}, Schramm and Sheffield~\cite{Schramm-Sheffield}, and also for the uniformly convex interactions (Miller~\cite{Miller}).

\subsection{Outlook}
\noindent
The upshot of the above overview that all of these problems, although quite varied in nature, can be reduced to specific properties of the Random Conductance Model. In particular, many of the solutions boil down to similar technical questions. In the rest of these notes we will attempt to explain the main ideas underlying the existing solutions and point out the obstacles that are known of for the problems that remain unresolved.
\vfill

\section{Limit laws for the RCM}
\label{sec2}\noindent
The goal of this section is to exhibit the main techniques that will allow us to establish the validity of the SLLN (Problem~\ref{P:SLLN}) and the Functional CLT (Problem~\ref{P:CLT}) for rather general Random Conductance Models. We will take a very pedagogical approach that starts off by addressing the simplest non-trivial cases of interest while isolating, as clearly as possible, various technical issues that come up along the way.

\subsection{Point of view of the particle}
\noindent
A  first basic problem that arises in analyzing the Markov chain $(X_n)$ for a fixed realization of the environment~$\omega$ is that the increments of this chain are not stationary. A way to mend this is to invoke the first fundamental idea encountered in the theory of random walks in random environment: the \emph{point of view of the particle}. Namely, instead of making a random walk run through a fixed environment, we will shift the environment around so that the walk remains always at the origin. Technically, this amounts to representing the sequence $(\tau_{X_n}\omega)$ as a trajectory of a Markov chain on the space of all environments.

\begin{lemma}
\label{lemma-stationary}
Suppose $\BbbP$ is translation invariant. Then $(\tau_{X_n}\omega)$ is a sample from a Markov chain on the space of conductances~$\Omega$ with the transition kernel
\begin{equation}
\scrP(\omega,\textd\omega'):=\sum_x\cmss P_\omega(0,x)\delta_{\tau_x\omega}(\textd\omega').
\end{equation}
Moreover, whenever $Z:=\E\pi_\omega(0)<\infty$, this chain has the stationary and reversible measure
\begin{equation}
\Q(\textd\omega):=\frac{\pi_\omega(0)}Z\BbbP(\textd\omega).
\end{equation}
\end{lemma}

\begin{proofsect}{Proof}
The fact that the kernel $\scrP$ generates the Markov chain $(\tau_{X_n}\omega)$ is a trivial calculation. For the second part, we need to invoke a bit of $L^2$-calculus. For any two bounded measurable functions $f=f(\omega)$ and $g=g(\omega)$, define
\begin{equation}
\langle\, f,g\rangle:=E_\Q\bigl(f(\omega)g(\omega)\bigr).
\end{equation}
This is a natural inner product in $L^2(\Q)$.
To show reversibility (and thus stationarity) of~$\Q$, it suffices to show that $\langle\, f,\scrP g\rangle=\langle \scrP f,g\rangle$ for any such bounded non-negative~$f,g$ --- in fact, indicators of measurable events would be enough. For that case we compute
\begin{equation}
\label{E:2.4}
\begin{aligned}
\langle\, f,\scrP g\rangle
&=\frac1Z\sum_x\E\Bigl(\pi_\omega(0)f(\omega)\cmss P_\omega(0,x)g\circ\tau_x(\omega)\Bigr)
\\
&=\frac1Z\sum_x\E\Bigl( f(\omega)\,\omega_{0,x}\,g\circ\tau_x(\omega)\Bigr).
\end{aligned}
\end{equation}
where all sums are meaningful by positivity of all terms and the assumption that $\sum_x\omega_{0,x}$ is integrable. Now apply $\tau_{-x}$ under the expectation to write this~as
\begin{equation}
\langle\, f,\scrP g\rangle
=\frac1Z\sum_x\E\Bigl( f(\tau_{-x}\omega)\,(\tau_{-x}\omega)_{0,x}\,g(\omega)\Bigr).
\end{equation}
A key property of the environment is its symmetry \eqref{E:omega-sym} whereby we get
\begin{equation}
\label{E:shift+sym}
(\tau_{-x}\omega)_{0,x}=\omega_{-x,0}=\omega_{0,-x}.
\end{equation}
Relabeling $-x$ for~$x$, we thus conclude
\begin{equation}
\label{E:2.7}
\langle\, f,\scrP g\rangle
=\frac1Z\sum_x\E\Bigl( f(\tau_{x}\omega)\,\omega_{0,x}\,g(\omega)\Bigr)
\end{equation}
which is, rolling back the first rewrite, exactly $\langle\scrP f,g\rangle$.
\end{proofsect}

It is not hard to check that, for any bounded $f,g$,
\begin{equation}
\label{E:2.8}
\bigl\langle\, f,(\text{\rm id}-\scrP) g\bigr\rangle
=-\frac1{2Z}\sum_x\E\Bigl( \omega_{0,x}\,\bigl(f(\tau_{x}\omega)-f(\omega)\bigr)\,\bigl(g(\tau_{x}\omega)-g(\omega)\bigr)\Bigr).
\end{equation}
This will help us solve:

\begin{hwproblem}
Show that, whenever $Z:=\E\pi_\omega(0)<\infty$, the operator~$\scrL:=\scrP-\text{\rm id}$ with domain
\begin{equation}
\text{\rm Dom}(\scrL):=\biggl\{f\in L^2(\Q)\colon\sum_x\E\bigl( \omega_{0,x}\,\bigl(f(\tau_{x}\omega)-f(\omega)\bigr)^2\bigr)<\infty\biggr\}
\end{equation}
is self-adjoint and negative semi-definite.
\end{hwproblem}

Notice that the stationary measure~$\Q$ and the \emph{a priori} law~$\BbbP$ are mutually absolutely continuous; we in fact even have a very explicit expression for~$\Q$. In the studies of general (non-reversible) random walks in random environments it is  (usually) not too hard  to infer the existence of a stationary measure but a key obstacle is the absolute continuity of~$\BbbP$ with respect to~$\Q$ --- which we often need to conclude that events that occur~$\Q$-a.s.\ also occur~$\BbbP$-a.s. But even in such cases it is unusual to have any sort of explicit handle of~$\Q$.

\smallskip
These considerations move us to the question under what conditions is the Markov chain $(\tau_{X_n}\omega)$ \emph{ergodic}. In order to explain this a bit better, recall that a stationary Markov chain $(Z_n)$ on a general state space with stationary measure~$\pi$ can always be embedded into a \emph{Markov shift} as follows: Sample the initial state~$Z_0$ from $\pi$ and then use the Markov kernel to sample a whole forward trajectory $(Z_n)_{n\ge1}$. If need be, also use the reversed chain to sample the entire backward trajectory $(Z_n)_{n<0}$. This defines --- through the Kolmogorov Extension Theorem --- a law $\mu$ on trajectories of the Markov chain. The canonical shift --- simply use $Z_n$ for the value of $Z_{n-1}$ for all $n$ --- then defines a \emph{measure preserving transformation}.

This construction and the Birkhoff-Khinchine Ergodic Theorem imply that, for $\pi$-almost every $Z_0$ and almost every path of the Markov chain --- in short, for $\mu$-almost every trajectory --- the limit
\begin{equation}
\lim_{n\to\infty}\frac1n\sum_{k=0}^{n-1}f(Z_k)
\end{equation}
exists and is finite for any function $f$ such that $f\in L^1(\pi)$. However, we often wonder whether this limit is in fact almost surely constant --- and this will only be true for a general~$f$ if the chain is ergodic. Explicitly, the above Markov chain is ergodic if any measurable set of trajectories~$A$ satisfies $\mu(A)\in\{0,1\}$.

\smallskip
Ergodicity will in our context be guaranteed by the following condition:

\begin{proposition}
\label{prop-ergodic}
Suppose~$\BbbP$ satisfies the following conditions:
\begin{enumerate}
\item[(1)]
$\BbbP(\pi_\omega(0)>0)=1$ and $\E\pi_\omega(0)<\infty$ --- i.e., $\Q$ exists and is equivalent to~$\BbbP$.
\item[(2)]
$\BbbP$ is irreducible in the sense that, for every~$x\in\Z^d$,
\begin{equation}
\BbbP\biggl(\,\omega\colon\sup_{n\ge0}\cmss P_\omega^n(0,x)>0\biggr)=1.
\end{equation}
\item[(3)]
$\BbbP$ is ergodic with respect to the translations of~$\Z^d$ --- i.e., $\BbbP(A)\in\{0,1\}$ for any event~$A$ such that $\tau_x(A)=A$ for all~$x$.
\end{enumerate}
Then the Markov chain $(\tau_{X_n}\omega)$ with initial law~$\Q$ is ergodic.
\end{proposition}

\begin{proofsect}{Proof}
The proof of this proposition is quite standard --- the result has been used at various levels of explicit detail in the literature --- although the general setting makes the use of ergodicity of~$\BbbP$ a bit subtle. Kozlov~\cite{Kozlov} proves this by way of a functional theoretical argument; we will follow a probabilistic argument from Berger and Biskup~\cite{Berger-Biskup}.

Let $A$ be the event on the space of trajectories $(\omega_n)_{n\in\Z}$ that is shift invariant. Explicitly, if $\theta$ is the Markov shift, $(\theta\omega)_n=\omega_{n+1}$, we have $\theta^{-1}(A)=A$. Let~$\mu$ denote the law of the trajectories induced by the Markov chain with stationary measure~$\Q$. Our goal is to show that $\mu(A)\in\{0,1\}$.

The first part of the proof is the classical approximation argument that drives the proof of more or less every known zero-one law. Define the function
\begin{equation}
f(\omega_0):=E_\mu(\1_A|\omega_0).
\end{equation}
We claim that~$f^2=f$ $\mu$-a.s. To this end approximate~$A$ by a sequence of events $A_n\in\sigma(\omega_{-n},\dots,\omega_n)$ so that
\begin{equation}
\label{E:2.13}
\Vert\1_A-\1_{A_n}\Vert_{L^1(\mu)}\,\underset{n\to\infty}\longrightarrow\,0.
\end{equation}
The shift invariance of~$A$ implies that the same holds for $A_n$ replaced by $\theta^n(A_n)$ and by $\theta^{-n}(A_n)$.

Invoking the general fact $E_\mu|E_\mu(g|\omega_0)|\le E_\mu|g|$ and applying \eqref{E:2.13}, we thus have
\begin{equation}
\bigl\Vert f-E_\mu(\1_{\theta^{\pm n}A_n}|\omega_0)\bigr\Vert_{L^1(\mu)}\,\underset{n\to\infty}\longrightarrow\,0.
\end{equation}
Similarly, replacing $\1_A$ by $\1_A\1_A$ and approximating the first indicator by $\1_{\theta^n(A_n)}$ and the second by $\1_{\theta^{-n}(A_n)}$ we obtain
\begin{equation}
\Vert f-E_\mu(\1_{\theta^{-n}(A_n)}\1_{\theta^{n}(A_n)}|\omega_0)\Vert_{L^1(\mu)}\,\underset{n\to\infty}\longrightarrow\,0.
\end{equation}
But $\theta^n(A_n)\in\sigma(\omega_0,\dots,\omega_{2n})$ and $\theta^{-n}(A_n)\in\sigma(\omega_{-2n},\dots,\omega_0)$ have only one coordinate in common and so, conditional on~$\omega_0$, they are independent. This means
\begin{equation}
E_\mu(\1_{\theta^{-n}(A_n)}\1_{\theta^{n}(A_n)}|\omega_0)=E_\mu(\1_{\theta^{-n}(A_n)}|\omega_0)E_\mu(\1_{\theta^{n}(A_n)}|\omega_0).
\end{equation}
Passing to $n\to\infty$, the right-hand side tends to $f(\omega_0)^2$ in $L^1(\mu)$ thus proving that $f=f^2$ $\mu$-a.s.

The second step is more subtle. Indeed, we claim that $f(\tau_x\omega)=f(\omega)$ for all~$x\in\Z^d$ and~$\Q$-almost every $\omega$. To this end let us note that, by the $\theta$-invariance of~$A$, if~$\omega_0$ is the initial configuration of a path in~$A$, then also $\omega_1$ is the initial step of a path in~$A$ --- namely, the shifted path! A moment's thought shows that this implies $f(\omega_0)=f(\omega_1)$ $\mu$-a.s. and thus
\begin{equation}
f(\tau_{X_n}\omega)=f(\omega)
\end{equation}
for $\Q$-a.e.~$\omega$ and $P_{\omega}^0$-a.e.\ trajectory $(X_n)$ of the Markov chain. The conditions on~$\BbbP$ guarantee that for $\BbbP$-a.e.~$\omega$, with positive probability $(X_n)$ visits any given~$x$ and so we must have $f(\tau_x\omega)=f(\omega)$. The event~$\{f=1\}$ is thus shift invariant and so $\BbbP(f=1)\in\{0,1\}$, by the ergodicity of~$\BbbP$. Then
\begin{equation}
\mu(A)=E_{\Q}(f)\underset{f=0,1}=\Q(f=1)\in\{0,1\},
\end{equation}
where we used that the $\omega_0$-marginal of~$\mu$ is~$\Q$ and that $\Q\sim\BbbP$.
\end{proofsect}

\subsection{Vanishing speed}
\label{sec2.2}\noindent
The conclusion of Lemma~\ref{lemma-stationary} and Proposition~\ref{prop-ergodic} can be formalized in multiple ways. E.g., we thus know that for any $f=f(\omega)$ with $E_{\Q}|f(\omega)|<\infty$,
\begin{equation}
\label{E:ergodic-1arg}
\lim_{n\to\infty}\frac1n\sum_{k=0}^{n-1}f(\tau_{X_k}\omega)=E_{\Q}f(\omega)
\end{equation}
for~$\BbbP$-a.e.~$\omega$ and $P_\omega^0$-a.e.\ path $(X_n)$. But since the convergence comes from the Markov shift, we are not limited to functions of only one argument. Thus, for instance, we also know that for any function $f=f(\omega_0,\omega_1)$ such that $E_{\Q}E_\omega^0 |f(\omega,\tau_{X_1}\omega)|<\infty$,
\begin{equation}
\label{E:ergodic-2arg}
\lim_{n\to\infty}\frac1n\sum_{k=0}^{n-1}f(\tau_{X_k}\omega,\tau_{X_{k+1}}\omega)=E_{\Q}E_\omega^0f(\omega,\tau_{X_1}\omega)
\end{equation}
for~$\BbbP$-a.e.~$\omega$ and $P_\omega^0$-a.e.\ path $(X_n)$. This permits us to prove:

\begin{theorem}[Vanishing speed]
\label{thm-SLLN}
Suppose~$\BbbP$ obeys assumptions (1-3) in Pro\-position~\ref{prop-ergodic} and
\begin{equation}
\label{E:omega-1st-moment}
\E\biggl(\,\sum_x\omega_{0,x}|x|\biggr)<\infty.
\end{equation}
Then for~$\BbbP$-a.e.~$\omega$ and $P_\omega^0$-a.e.~trajectory $(X_n)$,
\begin{equation}
\label{E:SLLN}
\lim_{n\to\infty}\frac{X_n}n=0.
\end{equation}
\end{theorem}

\begin{proofsect}{Proof}
Our key problem is to represent~$X_n$ as an \emph{additive functional} of the Markov chain $(\tau_{X_k}\omega)$. This can be done easily under the assumption that the environment is not periodic:
\begin{equation}
\label{E:non-periodic}
\BbbP\bigl(\{\omega\colon\tau_x\omega=\omega\}\bigr)=0,\qquad x\ne0.
\end{equation}
(Clearly, if the environment is periodic in some direction, there is no way for the walk to ``notice'' its motion through it when it makes a step in that direction.) We will thus prove the theorem only in this case leaving the periodic cases --- which for ergodic~$\BbbP$ are a.s.\ events --- to a (simple) Exercise afterwards.

We claim that, under \eqref{E:non-periodic}, we get
\begin{equation}
X_n=\sum_{k=0}^{n-1}f(\tau_{X_k}\omega,\tau_{X_{k+1}\omega})\quad\text{where}\quad
f(\omega,\omega'):=\sum_zz\1_{\{\omega'=\tau_z\omega\}}.
\end{equation}
Indeed, for almost every environment and any path of the chain, at most one of the indicators in the definition of~$f$ will be non-zero, and it is precisely the one that relates $\omega'$ to the shifted configuration~$\omega$.

We will now apply the conclusion \eqref{E:ergodic-2arg}, but to get the conclusion of the theorem we need to show that
\begin{equation}
E_{\Q}E_\omega^0 |f(\omega,\tau_{X_1}\omega)|<\infty\quad\text{and}\quad
E_{\Q}E_\omega^0f(\omega,\tau_{X_1}\omega)=0.
\end{equation}
This is a matter of a straightforward calculation. First,
\begin{equation}
E_{\Q}E_\omega^0 |f(\omega,\tau_{X_1}\omega)|=E_{\Q}E_\omega^0 |X_1|=\frac1Z\E\biggl(\,\sum_x\omega_{0,x}|x|\biggr)<\infty.
\end{equation}
Second, the absolute summability we just showed permits us to write
\begin{multline}
\qquad
E_{\Q}E_\omega^0f(\omega,\tau_{X_1}\omega)=E_{\Q}E_\omega^0X_1=\frac1Z\E\biggl(\,\sum_x\omega_{0,x}x\biggr)
\\=\frac1{2Z}\sum_x\E\bigl(\,\omega_{0,x}x+\omega_{0,-x}(-x)\bigr).
\qquad
\end{multline}
To see that the last expectation vanishes, recall \eqref{E:shift+sym} to see that $\E\omega_{0,x}=\E(\omega_{0,-x})$.
\end{proofsect}

The minor trouble with periodic configurations disappears if we encode the sequence of environments along with the corresponding (next) step of the walk. This is an approach that was taken in Kozlov~\cite{Kozlov}; however, the above works just as well. Indeed, we pose:

\begin{hwproblem}
Consider a product law on configurations~$(\omega,\sigma)$ where $\omega$ is sampled from~$\BbbP$ and $\sigma=(\sigma_x)_{x\in\Z^d}$ are i.i.d.\ (non-degenerate) Bernoulli. Show that this law is ergodic with respect to the Markov shift
\begin{equation}
(\omega,\sigma)\mapsto(\tau_{X_1}\omega,\tau_{X_1}\sigma) \text{ \rm with~$X_1$ sampled from~$P_\omega^0$}.
\end{equation}
 Find a function of the (joint) environment which encodes~$X_n$ as an additive function of two consecutive environments. Use this to conclude that \eqref{E:SLLN} still holds for almost every path of the Markov chain over~$\omega$, regardless of whether the aperiodicity condition \eqref{E:non-periodic} holds or not.
\end{hwproblem}

Notice that, for a shift-invariant configuration~$\omega$, the condition \eqref{E:omega-1st-moment} reduces exactly to the first moment condition in the SLLN. So \eqref{E:omega-1st-moment} should generally fail once \eqref{E:omega-1st-moment}  is violated, although exact conditions under which this is true do not seem to be available. The same should apply (under a different condition) when only convergence in measure is in question.

The following lemma, which arose in the writing of a proof in Biskup, Louidor, Rozinov and Vandenberg-Rodes~\cite{UCLA-team}, can sometimes be useful in applications:

\begin{lemma}
Let $f\in L\log L(\Q)$ and suppose that $\BbbP$ obeys assumptions (1-3) in Proposition~\ref{prop-ergodic}. Then for $\BbbP$-a.e.~$\omega$,
\begin{equation}
\lim_{n\to\infty}\frac1n E_\omega^0\biggl(\,\sum_{k=0}^{n-1}f(\tau_{X_k}\omega)\biggr) = E_\Q f(\omega).
\end{equation}
In particular, the limit exists $\BbbP$-a.s.
\end{lemma}

\begin{proofsect}{Proof}
Without loss of generality assume that~$f\ge0$ and recall that $L\log L(\Q)$ is the space of functions $f$ such that $f\log|f|\in L^1(\Q)$. By Wiener's Dominated Ergodic Theorem (e.g.,~Petersen~\cite[Theorem~1.16]{Petersen}) these functions are distinguished by the fact that
\begin{equation}
f^\star:=\sup_{n\ge1}\frac1n\,\sum_{k=0}^{n-1}f\circ\tau_{X_k}\in L^1(\Q).
\end{equation}
Since~$\frac1n\,\sum_{k=0}^{n-1}f\circ\tau_{X_k}$ are dominated by $f^\star$ and tend to $E_\Q f(\omega)$ $P_\omega^0$-almost surely for~$\BbbP$-a.e.~$\omega$, the result follows by the Dominated Convergence Theorem.
\end{proofsect}

A subtlety of the above statement is that although the averages $\frac1n\,\sum_{k=0}^{n-1}f\circ\tau_{X_k}$ converge almost surely and in $L^1(\Q\otimes P_\omega^0)$, this is not enough to guarantee convergence in $L^1(P_\omega^0)$, for $\BbbP$-a.e.~$\omega$. A useful step towards understanding this is solving:

\begin{hwproblem}
Construct a sequence of random variables $Z_n$ such that $Z_n\to Z$ almost surely in~$L^1$, but such that, for some $\sigma$-algebra~$\scrA$, the conditional expectations $E(Z_n|\scrA)$ do not converge almost surely.
\end{hwproblem}

The above arguments are useful even for the continuous-time versions of our random walk. Indeed, we can combine Exercise~\ref{ex1.11} with Theorem~\ref{thm-SLLN} to solve:

\begin{hwproblem}
\label{ex2.8}
Let~$\BbbP$ be ergodic with $\BbbP(\pi_\omega(0)>0)=1$ and $\E\pi_\omega(0)<\infty$. Then for~$\BbbP$-a.e.~$\omega$, the VSRW does not escape to infinity --- i.e., no blow-ups occur --- in finite time.
\end{hwproblem}

\subsection{Martingale (Functional) CLT}
\label{sec2.3}\noindent
Once a variant of the Law of Large Numbers has been established the next natural question is that of fluctuations. In order to discuss all aspects of this question in a reasonably pedagogical fashion, for a while we will restrict attention to a class of \emph{toy models} in which the environment has the following properties:

\begin{assumption}[Toy-model assumptions]
\label{toy-ass}
For some $\alpha\in(0,1)$ and~$\BbbP$-almost every~$\omega$,\looseness=1
\begin{enumerate}
\item[(1)] $\omega_{xy}=0$ unless $|x-y|=1$ (nearest-neighbor environment).
\item[(2)] For each coordinate vector~$\hate_i$ and each~$x\in\Z^d$, $\omega_{x-\hate_i,x}=\omega_{x,x+\hate_i}$.
\item[(3)] $\alpha\le\omega_{x,x+\hate_i}\le\frac1\alpha$ for all~$i$ and all~$x$.
\end{enumerate}
In other words, the environments are nearest-neighbor, elliptic and the conductances are constant along the edges on each line of sites in~$\Z^d$.
\end{assumption}

What makes these environments special is:

\begin{lemma}
\label{lemma-martingale}
Let~$\scrF_n:=\sigma(X_0,\dots,X_n)$. For all environments above, $\{X_n,\scrF_n\}$ is a martingale.
\end{lemma}

\begin{proofsect}{Proof}
Any environment satisfying conditions (1-2)  above has the property that the \emph{local drift},
\begin{equation}
\label{E:local-drift}
V(\omega):=E_\omega^0(X_1),
\end{equation}
identically vanishes. To see how this implies the claim we note that, by the Markov property the law of $X_{n+1}-X_n$ conditional on~$X_n$ is that of $X_1$ in distribution $P_{\tau_{X_n}\omega}^0$. Hence,
\begin{equation}
E_\omega^0(X_{n+1}|\scrF_n)=X_n+V(\tau_{X_n}\omega)=X_n
\end{equation}
and so~$X_n$ is a martingale.
\end{proofsect}

We remark that more general (particularly, non-reversible) cases of such \emph{balanced environments} have been treated by Lawler~\cite{Lawler}, Guo and Zeitouni~\cite{Guo-Zeitouni} and, quite recently, Berger and Deuschel~\cite{Berger-Deuschel}. The main issue dealt with in those papers is a construction, and proper control, of an ergodic, invariant law on environments.

\medskip
Returning to the setting of Toy Models, the fact that $X_n$ is a martingale with bounded increments  immediately implies, via Azuma's inequality, Gaussian bounds on its tails. Explicitly, for any unit vector~$\hate\in\R^d$ we will have
\begin{equation}
P_\omega^0\bigl(\hate\cdot X_n>\lambda\sqrt n\bigr)\le\texte^{-\lambda^2/2}.
\end{equation}
However, to get the desired CLT we will have to invoke a more delicate tool which is:

\begin{theorem}[Martingale Functional CLT]
\label{thm-MCLT}
Let~$\{M_n,\scrF_n\}$ be an $\R$-valued, square-integrable martingale such that the following conditions hold:
\settowidth{\leftmargini}{(XXXX)}
\begin{enumerate}
\item[(LF1)]
There is~$\sigma^2\in[0,\infty)$ such that for all~$t>0$,
\begin{equation}
\lim_{n\to\infty}\,\frac1n\,\sum_{k=0}^{\lfloor tn\rfloor}E\bigl(|M_{k+1}-M_k|^2\big|\scrF_k\bigr)
\,\underset{n\to\infty}\longrightarrow\,t\sigma^2
\end{equation}
in probability.
\item[(LF2)]
For each~$\epsilon>0$,
\begin{equation}
\lim_{n\to\infty}\,\frac1n\,\sum_{k=0}^{n}E\bigl(|M_{k+1}-M_k|^2\1_{|M_{k+1}-M_k|>\epsilon\sqrt n}\big|\scrF_k\bigr)
\,\underset{n\to\infty}\longrightarrow\,0
\end{equation}
in probability.
\end{enumerate}
Then for each~$T>0$, the law of
\begin{equation}
t\mapsto \frac1{\sqrt n}\bigl(M_{\lfloor tn\rfloor}+(tn-\lfloor tn\rfloor)(M_{\lfloor tn\rfloor+1}-M_{\lfloor tn\rfloor})\bigr)
\end{equation}
on $C([0,T]$, tends to the Wiener measure with $EB_t=0$ and $E B_t^2=t\sigma^2$.
\end{theorem}

This is what is sometimes referred to as the ``Lindeberg-Feller Functional CLT,'' although this is only thanks to the formulation which is borrowed from the context of sums of independent random variable (the Lindeberg-Feller CLT, see, e.g., Durrett~\cite{Durrett}). The result for martingales is, in this formulation, first due to Brown~\cite{Brown}. Derriennic~\cite{Derriennic} gave a thoughtful survey of these results; unfortunately, the full version of his paper is somewhat hard to get hold of.

A simple way how to understand the scaling of the martingale paths to Brownian motion is via \emph{Skorohod embedding}. Explicitly, we have:

\begin{theorem}[Skorohod~\cite{Skorohod}, Strassen~\cite{Strassen} and Dubins~\cite{Dubins}]
Supp\-ose that $\{M_n,\scrF_n\}$ is a square-inte\-grable (real-valued) martingale with $E(M_0)=0$. Then there is a sequence of integrable stopping times $(T_i)$ with $T_0=0$ and $T_{i+1}\ge T_i$, such that
\begin{equation}
\text{\rm Law of }(M_n)_{n\ge0}\,\,=\,\,\text{\rm Law of }(B_{T_n})_{n\ge0}.
\end{equation}
\end{theorem}

The history of this result is roughly as follows: Skorohod~\cite{Skorohod} noted its validity for sums of independent random variables, Strassen~\cite{Strassen} observed that it holds even for martingales and Dubins~\cite{Dubins} finessed an important technical detail where the construction of the stopping times can be done purely on the path-space of the Brownian motion (i.e., without reliance on additional random~variables).

Returning to the above Martingale CLT, condition (LF1) guarantees that $T_n/n\to\sigma^2$ which means that the time change between the martingale and the Brownian motion is asymptotically linear. The condition (LF2) ensures tightness in the space of continuous paths (i.e., the Brownian motion will not wiggle too far from the piece-wise linear path interpolating the martingale values). The Skorohod representation only applies to~$\R$-valued martingales, hence our restriction to those.

\begin{figure}[t]
\centerline{\includegraphics[width=0.65\textwidth]{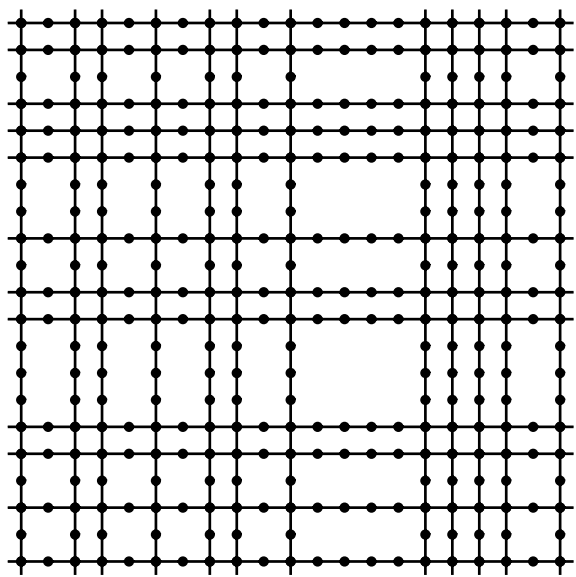}}
\caption{An example of the Random Conductance Model satisfying Assumptions~\ref{toy-ass}(1,2), but not the uniform ellipticity requirement in part~(3). Here, for each horizontal or vertical line of edges in~$\Z^2$, we independently retain, resp., drop all edges with probability~$p$, resp., $1-p$. The resulting random subgraph of~$\Z^2$ is almost-surely connected and the conclusion of Lemma~\ref{lemma-martingale} holds for almost every sample thereof.}
\label{fig3b}
\end{figure}

We can now finish the proof of:

\begin{proposition}
For any shift-ergodic environment law~$\BbbP$ satisfying (Toy Model) Assumptions~\ref{toy-ass} and for $\BbbP$-a.e.~sample from it, the law of $t\mapsto X_{\lfloor tn\rfloor}/\sqrt n$, linearly interpolated into a continuous path, tends to Brownian motion.
\end{proposition}

\begin{proofsect}{Proof}
We already know that~$X_n$ is a martingale for the filtration $\scrF_n:=\sigma(X_0,\dots,X_n)$ so we need to verify the conditions of the above theorem. This will be done again by using the point of view of the particle. By the Cram\'er-Wold device it suffices to prove the convergence for the projections onto all vectors in~$\R^d$. Fix a vector~$\hate\in\R^d$ and consider the function
\begin{equation}
f(\omega):=E_\omega^0\bigl(|\hate\cdot X_1|^2\bigr)
\end{equation}
and define $M_n:=\hate\cdot X_n$. The Markov property  guarantees
\begin{equation}
E_\omega^0\bigl(|M_{k+1}-M_k|^2\big|\scrF_k\bigr)=f(\tau_{X_k}\omega)
\end{equation}
and since $f$ is bounded and the environment is elliptic, (LF1) follows with $\sigma^2:=E_\Q f(\omega)$ by \eqref{E:ergodic-1arg}. The condition (LF2) is trivially satisfied and so we have the result.
\end{proofsect}

Notice the (somewhat counterintuitive) fact that we prove a CLT type of result by invoking a LLN type of result. But this is not so strange when we realize that for convergence to Brownian motion we need three things: asymptotically independent increments, their zero mean/second-moment property and their stationarity. The former two properties can be safely attributed to the use of martingales, but for the last one --- and, in this setting, the most difficult one --- we need to use the Ergodic Theorem and thus the machinery originally developed for the LLN.

\begin{hwproblem}
Consider the example of a random environment in Fig.~\ref{fig3b}. Show that, for almost every realization of this environment, the Martingale CLT applies. Characterize the variance of the limiting Brownian motion.
\end{hwproblem}

\subsection{Martingale approximations and other tricks}
\label{sec2.4}\noindent
The derivations in the preceding sections, however elegant, hinge on the crucial assumption of vanishing drift. Unfortunately, this is not what one can (and wants) to ask from a generic Random Conductance Model. Historically, this puts us somewhere in the first half of 1980s when people made first successful attempts to address the CLT in this level of generality. We will follow Kipnis and Varadhan~\cite{Kipnis-Varadhan} where the following strategy was taken:
\begin{enumerate}
\item[(1)]
Represent~$X_n$ as the sum of a martingale and an additive functional of (a single state of) the Markov chain on environments.
\item[(2)]
Approximate the additive functional by a martingale with an error that can be controlled at the level of the CLT.
\end{enumerate}
The first step can be achieved trivially:
\begin{equation}
X_n=\sum_{k=0}^{n-1}\bigl[X_{k+1}-X_k-E(X_{k+1}-X_k|\scrF_k)\bigr]+\sum_{k=0}^{n-1}E(X_{k+1}-X_k|\scrF_k)
\end{equation}
The first sum on the right is clearly a martingale --- call it $M_n$ --- while $E(X_{k+1}-X_k|\scrF_k)=V(\tau_{X_k}\omega)$ makes the second part an \emph{additive functional} of the Markov chain $(\tau_{X_k}\omega)$. (Note that we already know that~$X_n$ as additive functional of \emph{two} consecutive environments, but for the application of the Martingale Functional CLT the dependence on a single environment is much easier.) Now we need to~write
\begin{equation}
\sum_{k=0}^{n-1}
V(\tau_{X_k}\omega)=M_n'+E_n,
\end{equation}
where $\max_{k\le n}|E_k|/\sqrt n$ tends to zero in probability. This can be done under proper conditions but one then  faces the (rather extreme) difficulty that $M_n$ and $M_n'$ are \emph{not} independent.

\smallskip
To see how an additive functional of a Markov chain can be approximated by a martingale, consider a Markov chain on a state space~$\Omega$ with transition kernel~$\scrP$. Suppose~$g\colon\Omega\to\R$ is a function such that $g\in\text{Ran}(\text{id}-\scrP)$. In other words, we require
\begin{equation}
\label{Poisson-eq}
g=h-\scrP h
\end{equation}
for some function~$h\colon\Omega\to\R$. If $\omega_0,\omega_1,\dots$ denote the successive states of the Markov chain, then a similar trick to the one used above yields
\begin{equation}
\sum_{k=0}^{n-1}g(\omega_k)=h(\omega_0)-h(\omega_n)+\sum_{k=0}^{n-1}\bigl[h(\omega_{k+1})-\scrP h(\omega_k)\bigr].
\end{equation}
Set $E_n:=h(\omega_0)-h(\omega_n)$ and define $M_n'$ to be the sum. By the Markov property,
\begin{equation}
\scrP h(\omega_k) = E\bigl(h(\omega_{k+1})\big|\sigma(\omega_0,\dots,\omega_k)\bigr),
\end{equation}
which implies that $(M_n')$ is a martingale. Of course, in order to have a useful statement, we need that this martingale is properly integrable, which means that the Poisson equation \eqref{Poisson-eq} must be solved with~$h$ in, say, $L^2$. As we will comment in a minute, this may be quite a challenge to prove (and in fact, it is often too much to ask). However, such considerations are entirely unnecessary for finite-state Markov chains:

\begin{hwproblem}
Consider a Markov chain with a finite state space~$\Omega$ and a stationary measure~$\Q$. Let~$g\colon\Omega\to\R$ satisfy $E_\Q g=0$. Show that, for $\Q$-a.e.\ initial state~$\omega_0$, the law of
\begin{equation}
\label{E:2.43}
\frac1{\sqrt n}\sum_{k=0}^{n-1}g(\omega_k)
\end{equation}
tends to a mean-zero normal random variable. Characterize its variance.
\end{hwproblem}

This statement is actually one of the main results of a note due to Gordin and Lif\v sic~\cite{Gordin-Lifshitz}.
It will be easy to see that the result generalizes to arbitrary state spaces under the condition that $g\in\text{Ran}(\text{id}-\scrP)$ --- which we take to mean that \eqref{Poisson-eq} has a solution~$h\in L^2(\Q)$; the error,
\begin{equation}
E_n:=h(\omega_0)-h(\omega_n),
\end{equation}
is then trivially bounded in $L^2$. However, a bounded error is a luxury that we do not need; indeed, for the purpose of the CLT one can tolerate errors up to $o(\sqrt n)$ --- particularly, if that brings the benefit of weaker conditions on~$g$. A milestone achievement in this vain is the result of Kipnis and Varadhan~\cite{Kipnis-Varadhan} who proved the following theorem:

\begin{theorem}
\label{thm-KV}
Suppose that a Markov chain on state space~$\Omega$ with transition kernel~$\scrP$ is reversible with respect to~$\Q$. Consider the law on trajectories $(\omega_n)_{n\ge0}$ where $\omega_0$ is sampled from~$\Q$. Let $g\in L^2(\Q)$ with $E_\Q(g)=0$. Then the law of the \eqref{E:2.43} tends to a (zero-mean, finite-variance) normal random variable $\NN(0,\sigma_g^2)$ if and only if $g\in\text{\rm Ran}([\text{\rm id}-\scrP]^{-1/2})$ or, equivalently,
\begin{equation}
\label{E:2.44}
\sup_{\epsilon>0}\,E_\Q\bigl(\,g(\omega)(1+\epsilon-\scrP)^{-1}g(\omega)\bigr)<\infty.
\end{equation}
Moreover, the supremum equals $\sigma_g^2$ and the convergence extends (with the limit given by Brownian motion) even to paths (linearly) interpolating the values of $t\mapsto n^{-1/2}\sum_{k=0}^{\lfloor tn\rfloor} g(\omega_k)$.
\end{theorem}

Note that the claim concerns the averaged law; no statement about a typical starting point~$\omega_0$ is made. This is one of the deficiencies we will have to address in detail when proving the quenched invariance principle in the next two sections. The original method of proof in~\cite{Kipnis-Varadhan} was to consider the \emph{spectral measure} $\mu_g$ associated with the function~$g$ and the operator~$\scrP$ on $L^2(\Q)$. This measure has the property that, for any $F\in L^1(\mu_g)$,
\begin{equation}
\langle\, g,F(\scrP)g\rangle_{L^2(\Q)}=\int F(\lambda)\mu_g(\textd\lambda).
\end{equation}
The Kipnis-Varadhan condition \eqref{E:2.44} can then be written as
\begin{equation}
\label{E:sigma-g}
\sigma_g^2=\int\frac1{1-\lambda}\,\mu_g(\textd\lambda)<\infty.
\end{equation}
Notice that the spectrum of~$\scrP$, and thus the support of~$\mu_g$, is contained in $[-1,1]$.

\begin{hwproblem}
Show that if $(\omega_n)$ is a stationary Markov chain with~$\omega_0$ distributed according to~$\Q$, and $g\in L^2(\Q)$, then the variance of \eqref{E:2.43} tends to the quantity in \eqref{E:sigma-g}.
\end{hwproblem}

The spectral measure is a very interesting object in its own right due to the connection with the area of random Schr\"odinger operators. What is quite puzzling is that we do not have any substantive information to report on:

\begin{problem}
Describe the connection between the spectral properties of the generator~$\cmss L_\omega$ of the random walk among conductances~$\omega$ --- many of which, as is well known, are same for a.e.~$\omega$ --- and the generator $\scrL:=\scrP-\text{\rm id}$ of the Markov chain on environments.
\end{problem}

Let us make some remarks on how the history of the above ideas seems to have evolved. First, the idea to decompose additive functionals (of general stationary ergodic processes) into a martingale and an error is presumably due to Gordin~\cite{Gordin} who also had the insight to characterize the objects in terms of their functional-analytic (rather than mixing) properties. Gordin and Lif\v sic~\cite{Gordin-Lifshitz} then applied this idea in the specific context of finite-state Markov chains.

The understanding that martingale approximations can be the ultimate passage to limit laws for random walks in random environment seems to have grown out of the work of Papanicolaou and Varadhan~\cite{Papanicolaou-Varadhan}; the predecessors of this work were mostly focused on periodic environments. An alternative approach based on resolvent methods was devised by K\"unnermann~\cite{Kuennermann}.
The above (Kipnis-Varadhan) Theorem~\ref{thm-KV}  more or less closed the matter for the annealed law in reversible cases. Two natural ways to generalize Theorem~\ref{thm-KV} are as follows: One is to go beyond the annealed law and the other is to extend beyond reversible Markov chains. Both of these directions are far from settled and both constitute a subject of intense research.

We will expound on how to go from annealed to quenched laws in the rest of these notes. Concerning departures from reversible situations, two lines of thought are  generally being followed: One approach, drawing on the functional-analytic ideas, goes by imposing (and checking) various \emph{sector conditions} (e.g., Olla~\cite{Olla}, Sethuraman, Varadhan and Yau~\cite{SVY}, Horv\'ath, T\'oth and Vet\H o~\cite{HTV}). The role of these conditions is to control the antisymmetric (``non-reversible'') part of the generator by the symmetric one. Another approach goes by imposing decay-rate conditions on  time-correlations (e.g., Maxwell and Woodroofe~\cite{Maxwell-Woodroofe}, Derriennic and~Lin~\cite{Derriennic-Lin}, Peligrad and Utev~\cite{Peligrad-Utev}, Klicnarov\'a and Voln\'y~\cite{Klicnarova-Volny}, Voln\'y~\cite{Volny}, etc.). However, unlike the reversible situations, it does not seem likely that a single condition will eventually cover all cases of interest.


\section{Harmonic embedding and the corrector}
\label{sec3}
\noindent
Although the subject of martingale approximations is very attractive and useful, in the sequel we will adopt a different approach that emphasizes the geometrical component of the problem over its analytical component. To motivate this approach, consider the explicit example of the simple random walk on the two-dimensional supercritical percolation cluster. When the local drift $V(\omega)$ is non-zero, then this is because there is an odd number of neighbors of the origin and the origin thus no longer lies in the barycenter of its neighbors. The martingale defect can therefore be thought to arise from the use of the \emph{geometric} embedding of the graph, before the edges got removed.

\begin{figure}[t]
\centerline{\includegraphics[width=0.45\textwidth]{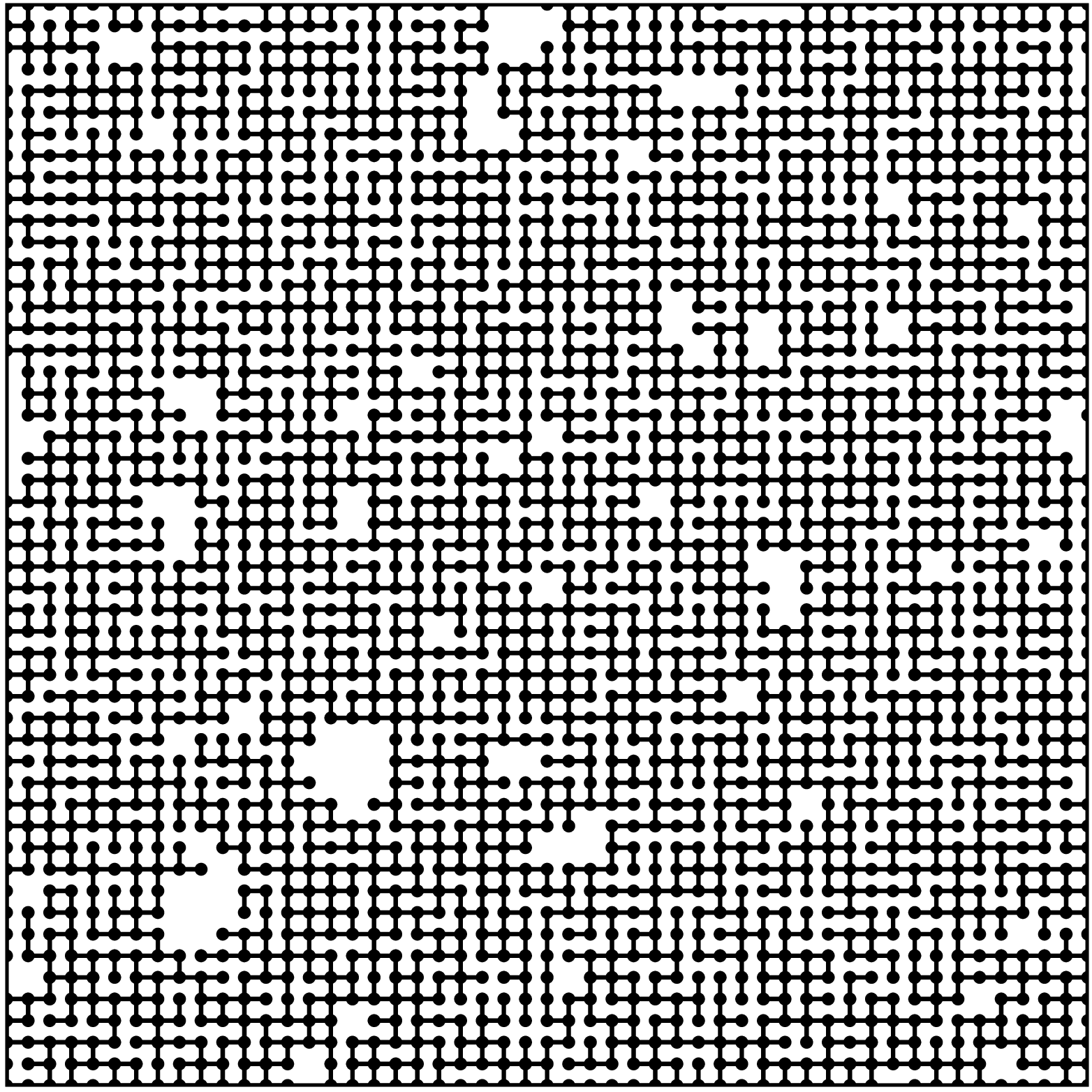}\hfil\includegraphics[width=0.45\textwidth]{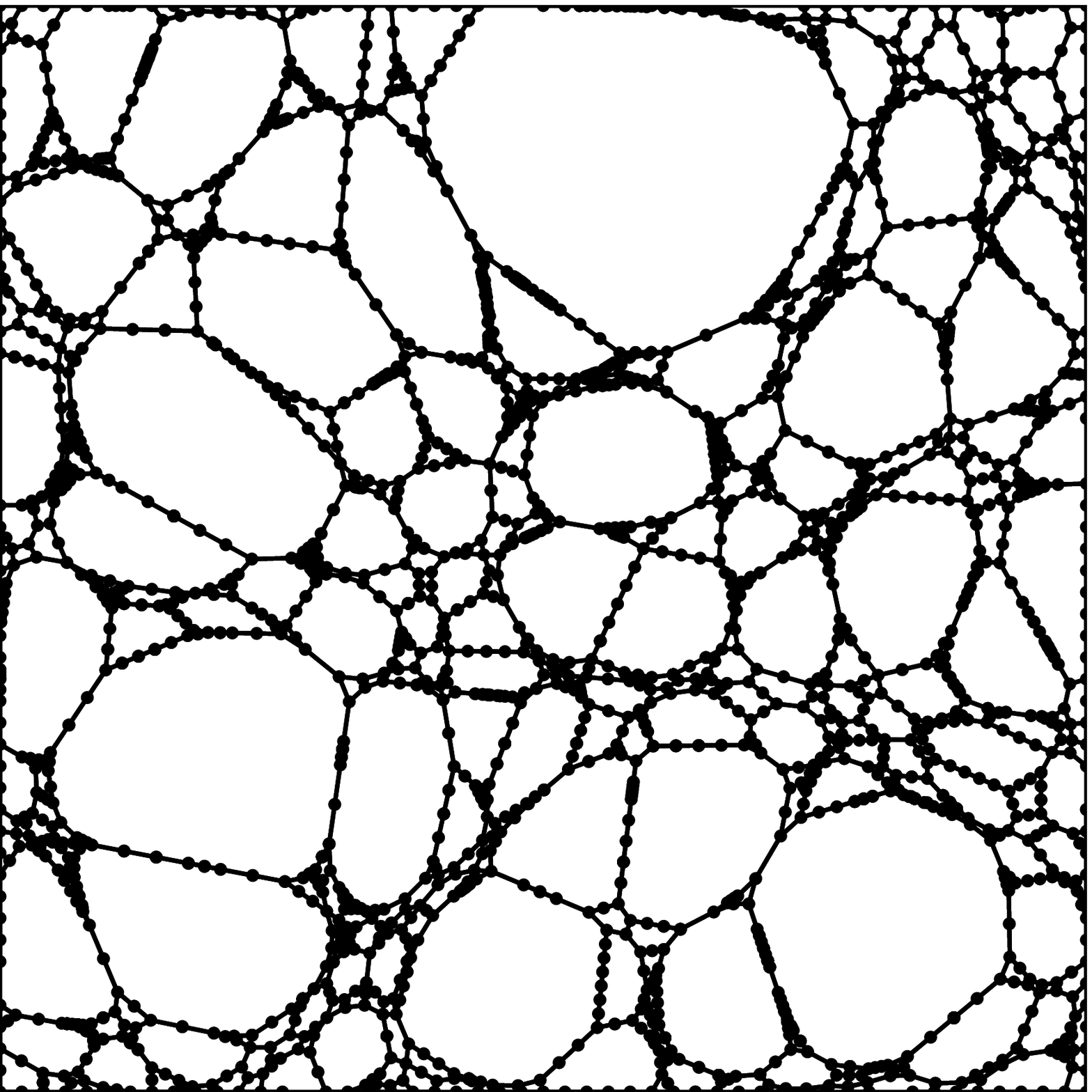}}
\caption{{\sc Left:} A sample of a percolation cluster at $50\times50$ square with bond probability $p=0.65$. Only vertices connected to the boundary are retained. {\sc Right:} The corresponding harmonic deformation obtained by relaxing the positions (except those on the boundary) to make each vertex lie in the ``center of mass'' of its (graph-theoretic) neighbors. Note that all dangling ends --- parts of the cluster attached only by one edge --- collapse to a point while the components attached by exactly two edges line up along a linear segment.}
\label{fig4}
\end{figure}

This suggests an idea that one might instead try to look for a different, \emph{harmonic} embedding for which $V$ would trivially vanish. A moment's thought shows that such an embedding is easy to find in any finite box using a computer --- just freeze the positions on the boundary and then ask the computer to sequentially pass through all vertices and always put them at the center of mass of their (graph-theoretic) neighbors. It turns out that this procedure rapidly converges and  leads to a picture as in Fig.~\ref{fig4}. How such an embedding is generated without recourse to finite volume is a slightly more complicated, although not unsolvable problem. The main new ingredient will be the reliance on homogenization theory.

Here and henceforth we will make repeated use of this notion:

\begin{definition}
We will henceforth say that~$\BbbP$ obeys the ``usual conditions'' if it satisfies the conditions (1-3) in Proposition~\ref{prop-ergodic}.
\end{definition}

These are exactly the conditions that guarantee the existence and ergodicity of the Markov chain on the space of environments.

\subsection{Minimizing Dirichlet energy}
\label{sec3.1}\noindent
We begin with some motivational observations for general reversible Markov chains that will explain in more detail how Fig.~\ref{fig4} was generated.
Suppose a countable set $\scrV$ is given along with the collection of (non-negative) conductances $(\omega_{xy})_{x,y\in\scrV}$ subject to restrictions \twoeqref{E:pi>0}{E:omega-sym}. Suppose in addition the irreducibility condition: for each~$x,y\in\scrV$, there is an $n\ge0$ with $\cmss P_\omega^n(x,y)>0$. For a finite set~$A\subset\scrV$ we then define
\begin{equation}
\label{E:dir-A}
\EE_A(f):=\frac12\sum_{\begin{subarray}{c}
x,y\in\scrV\\\{x,y\}\cap A\ne\emptyset
\end{subarray}}
\omega_{xy}\bigl[f(y)-f(x)\bigr]^2
\end{equation}
to be the Dirichlet energy in~$A$ for the potential~$f$. The following is well known:

\begin{lemma}[Dirichlet principle]
\label{lemma-Dir}
Let~$A\subset\scrV$ be a finite set with $\scrV\setminus A\ne\emptyset$ and let $g\colon\scrV\to\R$ be a bounded function. Then the infimum
\begin{equation}
\inf\bigl\{\EE_A(f)\colon f_{\scrV\smallsetminus A}=g_{\scrV\smallsetminus A}\bigr\}
\end{equation}
is achieved by the unique solution to the Dirichlet problem
\begin{equation}
\label{E:3.2}
\begin{cases}
\cmss L_\omega f(x)=0,\qquad&\text{if }x\in A,
\\
f(x)=g(x),\qquad&\text{if }x\in\scrV\setminus A.
\end{cases}
\end{equation}
\end{lemma}

\begin{proofsect}{Proof}
Pick~$x\in\scrV$ and any function~$f$. Let $f_x$ be defined by
\begin{equation}
f_x(z):=\begin{cases}
f(z),\qquad&\text{if }z\ne x,
\\
\cmss P_\omega f(x),\qquad&\text{if }z=x.
\end{cases}
\end{equation}
We claim that whenever $x\in A$, the ``move'' $f\mapsto f_x$ demonstrably lowers the Dirichlet energy, $\EE_A(f)\ge\EE_A(f_x)$. This is seen from the identity
\begin{multline}
\qquad
\sum_y\omega_{xy}\bigl[f(y)-f(x)\bigr]^2
\\=\sum_y\omega_{xy}\bigl[f(y)-\cmss P_\omega f(x)\bigr]^2+\pi_\omega(x)\bigl[f(x)-\cmss P_\omega f(x)\bigr]^2,
\qquad
\end{multline}
which is proved by optimizing the left-hand side over possible $f(x)$ --- this shows that the minimum is achieved at $\cmss P_\omega f(x)$ --- and  using that, for $h(x):=Ax^2+Bx+C$ with $A>0$, if $x_{\text{min}}$ is the minimizer then $h(x)=h(x_{\text{min}})+A(x-x_{\text{min}})^2$.

The explicit control of  $\EE_A(f)-\EE_A(f_x)$ shows that, applying the averaging $f(x)\mapsto \cmss P_\omega f(x)$ keeps lowering the Dirichlet energy as long as $f(x)\ne\cmss P_\omega(f)$. Furthermore,
\begin{equation}
\label{E:3.7}
\inf_{z\in\scrV}f(z)\le\inf_{z\in\scrV}f_x(z)\le\sup_{z\in\scrV}f_x(z)\le\sup_{z\in\scrV}f(z),
\end{equation}
and so any minimizing sequence of~$\EE_A(f)$ in \eqref{E:3.2} is bounded. Reducing to subsequences if needed, we extract a limit which then obeys $f=\cmss P_\omega f$ on~$A$ and thus solves the Dirichlet problem. To see that the solution is unique, note that $f=\cmss P_\omega f$ on~$A$ implies that~$f$ cannot have (strict) local extrema inside~$A$. In particular, we have the \emph{maximum principle}:
\begin{equation}
\label{E:3.7a}
\inf_{z\in\scrV\smallsetminus A}\,g(z)\le\min_{z\in A}f(z)\le\max_{z\in A}f(z)\le\sup_{z\in\scrV\smallsetminus A}\,g(z).
\end{equation}
Linearity guarantees that the difference between two solutions to \eqref{E:3.2} solves \eqref{E:3.2} with $g:= 0$. The maximum principle ensures that the difference must be zero.
\end{proofsect}

The above proof suggests that we could perhaps use the Dirichlet energy as a kind of measure of distance from a harmonic function. We will explore this very soon in a more general context. However, the argument also highlights a difficulty associated with attempts to ``harmonize'' the linear function $f(x)=x$ in infinite volume. Indeed, the full-lattice Dirichlet energy of such an~$f$ is infinity and so the procedure does not make sense.

This problem is not unknown from other situations and it naturally leads us to a guiding principle of \emph{homogenization theory}: Instead of trying to find the deformation of the linear function $f(x)=x$ that is harmonic with respect to $\cmss L_\omega$ at all locations for one given~$\omega$, we will solve the problem at one specific location --- namely the origin --- but simultaneously for all~$\omega$. Technically, this amounts to replacing the space $\ell^2(\pi_\omega)$ associated with the Markov chain $(X_n)$ by the space $L^2(\Q)$ associated with the chain~$(\tau_{X_k}\omega)$. The advantage of working on~$L^2(\Q)$ is that, unlike~$\pi_\omega$, the measure $\Q$ is finite.

\begin{figure}[t]
\centerline{\includegraphics[width=3.5in]{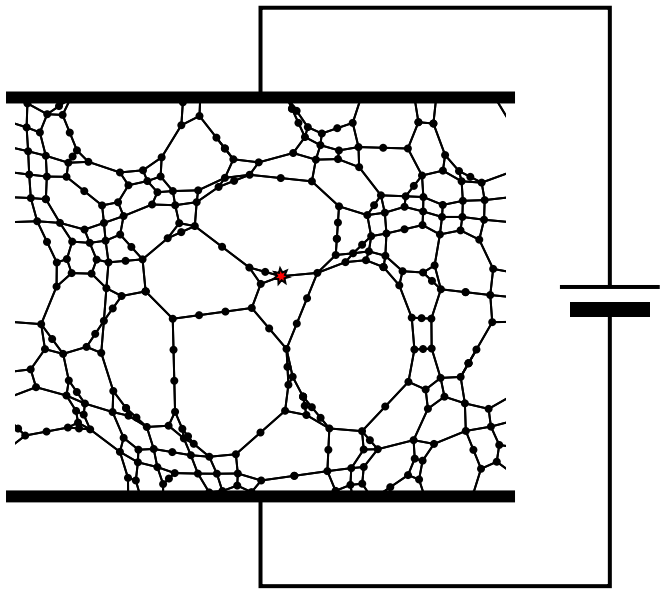}}
\caption{A graphical solution of the electrostatic problem depicted in Fig.~\ref{fig2} based on a harmonic deformation of the underlying graph. The electrostatic potential changes linearly in the height (more precisely, the $y$-coordinate) of the point. In particular, the potential at the vertex marked by the star --- originally, the origin of coordinates --- is proportional to the ratio between its distance to the top and the bottom plates.}
\label{fig7}
\end{figure}

\newcommand\Lvec{L_{\text{\rm cov}}}
\newcommand\GRAD{\text{\rm grad}\,}
\newcommand\DIV{\text{\rm div}}

\subsection{Weyl decomposition and the corrector}
\label{sec3.2}\noindent
To motivate the forthcoming definitions, recall that the process of substituting $\cmss P_\omega f(x)$ for  $f(x)$  applied to the function $f(x):=x$ would replace the value~$x$ by $x+E_{\tau_x\omega}^0(X_1)$. From the point of view of the particle it makes sense to shift this so that the origin of coordinates will not be moved under this action and so we may in fact want to replace $x$ by $x+E_{\tau_x\omega}^0(X_1)-E_\omega^0(X_1)$. The difference $E_{\tau_x\omega}^0(X_1)-E_\omega^0(X_1)$ is in the form of a \emph{gradient},
\begin{equation}
\nabla_x\varphi(\omega):=\varphi(\tau_x\omega)-\varphi(\omega).
\end{equation}
We are thus led to minimizing the functional
\begin{equation}
\label{E:Dir}
\varphi\mapsto \E\Bigl(\,\sum_x\omega_{0,x}\bigl|x+\nabla_x\varphi(\omega)\bigr|^2\Bigr)
\end{equation}
over all, say, local functions $\varphi=\varphi(\omega)$. Here we recall that $\varphi=\varphi(\omega)$ is said to be \emph{local} if it is a bounded, continuous function of a finite number of $\omega_{0,x}$'s. To see how homogenization translates finite-volume quantities to functionals over space of environments, it is instructive to solve:

\begin{hwproblem}
Consider the Dirichlet energy~$\EE_{\Lambda_N}(f)$ from \eqref{E:dir-A} for the set $\Lambda_N:=[-N,N]^d\cap\Z^d$. Fix a local function~$\varphi=\varphi(\omega)$ and set~$f$ to $x+\nabla\varphi$ defined by $(x+\nabla\varphi)(z):=z+\varphi\circ\tau_z(\omega)$ for~$z\in\Z^d$. Assuming~$\omega$ is a sample from an ergodic measure~$\BbbP$, carefully check that
\begin{equation}
\lim_{N\to\infty}\,\frac1{|\Lambda_N|}\EE_{\Lambda_N}(x+\nabla\varphi)=
\E\biggl(\,\sum_x\omega_{0,x}\bigl|x+\nabla_x\varphi(\omega)\bigr|^2\biggr).
\end{equation}
\end{hwproblem}

For technical reasons it will be advantageous to interpret \eqref{E:Dir} as a quadratic form on vector fields. Let~$\NN$ denote the set of admissible jumps of the Markov chain,
\begin{equation}
\NN:=\bigl\{x\in\Z^d:\BbbP(\omega_{0,x}>0)>0\bigr\}\cup\{0\}
\end{equation}
By a \emph{vector field} we will then mean a (measurable) map $u\colon\Omega\times\NN\to\R^d$, i.e., a vector valued function $u=u(\omega,x)$ indexed by environments and points in~$\NN$. We will always set
\begin{equation}
\label{E:convention}
u(\omega,0):=0
\end{equation}
by definition.

An example of a vector field is a \emph{potential field} $\nabla\varphi$ where $\nabla\varphi(\omega,x):=\nabla_x\varphi(\omega)$. Any potential field is \emph{curl-free} in the sense that it obeys the \emph{cycle conditions}. These conditions state that for any sequence $x_0,x_1,\dots,x_n:=x_0$ of vertices in $\Z^d$ such that $x_{i+1}-x_i\in\NN$ for all~$i$ we have
\begin{equation}
\sum_{i=0}^{n-1}u(\tau_{x_i}\omega,x_{i+1}-x_i)=0.
\end{equation}
In light of our convention \eqref{E:convention}, whenever~$\NN$ generates all of~$\Z^d$ (as an additive group), this turns out to be equivalent to
\begin{equation}
\label{E:shift-cov}
u(\omega,x+z)-u(\omega,x)=u(\tau_x\omega,z),\qquad \text{whenever }x,z,x+z\in\NN,
\end{equation}
The vector fields that obey this property (for all~$\omega$) will be called \emph{shift covariant} (sometimes they are called \emph{stationary}). Note that from \eqref{E:shift-cov} we automatically have $u(\omega,x)=-u(\tau_x\omega,-x)$.

As already alluded to, all potential fields are shift covariant. Another example of a shift-covariant field is the \emph{position field}, $x(\omega,z):=z$. As we shall see later, the position field and the potential fields generate the vector space of all shift-covariant fields. The reason for singling out shift-covariant fields is that they correspond to gradients of lattice functions. The following exercise details this connection:

\begin{hwproblem}
Assume the irreducibility condition $\BbbP(\sup_{n\ge1}\cmss P_\omega^n(0,x)>0)=1$, for all $x\in\Z^d$. Show that for any shift covariant~$u$ there is a ($\BbbP$-a.s.) unique function $U=U(\omega,x)$ with~$x\in\Z^d$ such that
\begin{equation}
\label{E:3.15a}
U(\omega,0)=0
\end{equation}
 and
\begin{equation}
\label{E:3.16a}
U(\omega,x+z)-U(\omega,x)=u(\tau_x\omega,z),\qquad x\in\Z^d,\,z\in\NN.
\end{equation}
\end{hwproblem}

To indicate that the vector field~$u(\omega,x)$ and the function~$U(\omega,x)$ are related as in \twoeqref{E:3.15a}{E:3.16a}, we will sometimes write $u=\GRAD U$ or say that~$U$ is an extension of~$u$ to~$\Z^d$.

The expression
\begin{equation}
\langle v,w\rangle:=\E\biggl(\,\sum_{x\in\NN} \omega_{0,x}\,v(\omega,x)\cdot w(\omega,x)\biggr)
\end{equation}
defines a natural inner product on the set of vector fields; the dot in $v(\omega,x)\cdot w(\omega,x)$ stands for the usual (Euclidean) dot product in~$\R^d$.
This inner product defines a natural $L^2$-norm; a minor technical problem --- which has often been overlooked in the literature --- is that $\langle u,u\rangle=0$ does not imply that $u=0$, only that $\omega_{0,x}u(\omega,x)=0$ for all~$x\in\NN$. A standard approach would be to factor the space of vector fields by the equivalence relation $u\sim u'$ whenever $\langle u-u',u-u'\rangle=0$. However, this is unnecessary once we restrict attention to shift-covariant fields (and impose a proper non-degeneracy condition). Indeed, define the set
\begin{equation}
\Lvec^2:=\bigl\{u\colon \text{\rm shift covariant},\,\langle u,u\rangle<\infty\bigr\}
\end{equation}
and set $\Vert u\Vert_{\Lvec^2}:=\langle u,u\rangle^{1/2}$. It is then not too hard to solve:

\begin{hwproblem}
Assume the irreducibility condition $\BbbP(\sup_{n\ge1}\cmss P_\omega^n(0,x)>0)=1$, $x\in\Z^d$. If~$\Vert u\Vert_{\Lvec^2}=0$ then $u(\omega,x)=0$ for all~$x\in\NN$ and $\BbbP$-a.e.~$\omega$.
\end{hwproblem}

Once the $L^2$-structure is in place, we note that potential fields define a natural closed subspace
\begin{equation}
L_\nabla^2:=\overline{\{\nabla\varphi\colon \varphi\text{ local}\}}^{\,\lower1pt\hbox{$\scriptstyle\Lvec^2$}}
\end{equation}
of $\Lvec^2$. With this space comes the orthogonal decomposition
\begin{equation}
\label{E:orth}
\Lvec^2=L_\nabla^2\oplus(L_\nabla^2)^\perp.
\end{equation}
It turns out that the vector fields from $(L_\nabla^2)^\perp$ can be quite well characterized. To see that explicitly, define the \emph{divergence} $\DIV(\omega u)$ by the formula:
\begin{equation}
\DIV(\omega u):=\sum_x\omega_{0,x}\bigl[u(\omega,x)-u(\tau_x\omega,-x)]
\end{equation}
where the bracket simplifies to $2u(\omega,x)$ once~$u$ is shift covariant. Thinking of $u(\omega,x)$ as the flux from~$0$ to~$x$, the first term on the right (including $\omega_{0,x}$) corresponds to the total flux out of the origin and the second one to the flux into the origin.

\begin{lemma}
\label{lemma-sol}
For~$u\in\Lvec^2$, we have $u\in(L_\nabla^2)^\perp$ if and only if $\DIV(\omega u)=0$ for $\BbbP$-a.e.~$\omega$. In particular, if~$U$ is a function such that $u=\GRAD U$, then $\cmss L_\omega\,U(\omega,x)=0$ at all~$x$ and $\BbbP$-a.e.~$\omega$.
\end{lemma}

\begin{proofsect}{Proof}
Pick a local function $\varphi$ and note that
\begin{equation}
\begin{aligned}
\langle u,\nabla\varphi\rangle
&=\sum_x \,\E\Bigl(\omega_{0,x}\,u(\omega,x)\cdot\bigl[\varphi\circ\tau_x(\omega)-\varphi(\omega)\bigl]\Bigr)
\\
&=\sum_x \,\E\Bigl(\varphi(\omega)\bigl[(\tau_{-x}\omega)_{0,x}\,u(\tau_{-x}\omega,x)-\omega_{0,x}\,u(\omega,x)\bigl]\Bigr)
\\
&= -\E\bigl(\varphi(\omega)\text{div}(\omega u)\bigr),
\end{aligned}
\end{equation}
where we used $u\in\Lvec^2$ to split the second expectation into two terms and then relabeled $x$ for~$-x$.
It follows that if $\langle u,\nabla\varphi\rangle=0$ for all local functions, then $\DIV(\omega u)=0$ $\BbbP$-a.s.\ and \emph{vice versa}.

For~$u=\GRAD U$, a simple calculation shows $\DIV(\omega u)=2\cmss L_\omega\,U$. With the help of shift covariance the condition $\DIV(\omega u)=0$ then forces $\cmss L_\omega\,U(\omega,\cdot)=0$.
\end{proofsect}

Lemma~\ref{lemma-sol} shows that the fields in $(L_\nabla^2)^\perp$ are, after multiplication by $\omega$, necessarily \emph{divergence-free} --- and are thus sometimes referred to as \emph{solenoidal} fields. The orthogonal decomposition \eqref{E:orth} is thus an analogue of the \emph{Weyl decomposition} from differential geometry. For readers familiar with basic electrostatics, the function~$U$ --- associated to a shift-covariant field~$u$ --- can be thought of as an electrostatic potential while $\omega u$ plays the role of an electric current. The fact that potential difference and current are related by way of a multiplication by~$\omega$ is a demonstration of \emph{Ohm's law} of electrostatics. See Doyle and Snell~\cite{Doyle-Snell} and/or Sect.~\ref{sec6}.

\smallskip
A natural next question to ask now is whether there are any solenoidal fields at all. For nearest-neighbor, constant conductances, a perfect candidate for a solenoidal field is the \emph{position field} which simply assigns $x(\omega,x):=x$. (Indeed, this function is discrete harmonic with respect to the homogeneous Laplacian on $\Z^d$ and so it obeys the conclusion of previous lemma.) Of course, once the conductances are not constant, $\text{div}(\omega x)$ --- which equals twice the local drift $V(\omega)$ --- is generally non-zero, but one can still hope that~$x$ has a non-trivial projection into the subspace~$(L_\nabla^2)^\perp$. This is all expressed in:

\begin{proposition}
\label{prop-correct}
Suppose $\BbbP$ obeys the ``usual conditions'' and, in addition, assume that
\begin{equation}
\label{E:square}
\E\biggl(\sum_x\omega_{0,x}|x|^2\biggr)<\infty.
\end{equation}
Then there is a function $\Psi=\Psi(\omega,x)$ defined for all~$x\in\Z^d$ with the properties:
\begin{enumerate}
\item[(1)]
Harmonicity: $\cmss L_\omega\Psi(\omega,x)=0$ for all~$x\in\Z^d$ and $\BbbP$-a.e.\ $\omega$.
\item[(2)]
Shift covariance: $\Psi(\omega,0)=0$ and
\begin{equation}
\Psi(\omega,x+z)-\Psi(\omega,x)=\Psi(\tau_x\omega,z),\qquad x,z\in\Z^d.
\end{equation}
\item[(3)]
Square integrability: $E_\Q E_\omega^0|\Psi(\omega,X_1)|^2<\infty$.
\end{enumerate}
In addition, for any minimizing sequence~$\varphi_n$ of the function \eqref{E:Dir}, we have $\nabla\varphi_n\to\chi(\omega,x)$ in $\Lvec^2$ where~$\chi$ is the \emph{corrector} that is given by
\begin{equation}
\label{E:corrector}
\chi(\omega,x):=\Psi(\omega,x)-x.
\end{equation}
Furthermore, $\Psi=\Psi(\omega,x)$ with~$x$ restricted to~$\NN$ is the orthogonal projection
\begin{equation}
\label{E:proj}
\Psi(\omega,\cdot):=\text{\rm proj}_{(L_\nabla^2)^\perp}x(\omega,\cdot).
\end{equation}
The infimum of \eqref{E:Dir} over all $\nabla\varphi\in L_\nabla^2$ is exactly $\Vert\Psi\Vert_{\Lvec^2}^2$.
\end{proposition}

\begin{proofsect}{Proof}
The proof could be simply started by defining $\Psi$ via \eqref{E:proj} and then checking the stated properties based on facts from the theory of abstract Hilbert spaces. However, it will be more instructive to prove some of the those claims directly in the present setting.

First note that the object in \eqref{E:Dir} can be interpreted as
\begin{equation}
\E\biggl(\,\sum_x\omega_{0,x}\bigl|x+\nabla_x\varphi(\omega)\bigr|^2\biggr)=\Vert x+\nabla\varphi\Vert^2_{\Lvec^2}.
\end{equation}
The condition \eqref{E:square} then guarantees that \eqref{E:Dir} takes a finite value for all local functions. Since it is also positive, we can pick a sequence $\varphi_n$ for which it tends to its infimum. The parallelogram law then yields
\begin{multline}
\label{E:3.25}
\qquad
\frac12\Vert\nabla\varphi_n-\nabla\varphi_m\Vert^2_{\Lvec^2}
=\Vert x+\nabla\varphi_n\Vert^2_{\Lvec^2}+\Vert x+\nabla\varphi_m\Vert^2_{\Lvec^2}
\\-2\Bigl\Vert x+\nabla\frac{\varphi_n+\varphi_m}2\Bigr\Vert^2_{\Lvec^2}.
\end{multline}
The first two terms on the right both tend to the infimum while the last term is bounded by twice the infimum. It follows that $\nabla\varphi_n$ is Cauchy in~$\Lvec^2$ and so it converges to a vector field that we denote by~$\chi$. This is the corrector in \eqref{E:corrector}.

Since $\chi$ is a limit of gradients, it is shift-covariant and so it extends to a unique function on~$\Z^d$. Now we define~$\Psi:=x+\chi$ and note that $\Vert \Psi\Vert^2_{\Lvec^2}$ is the infimum of \eqref{E:Dir}. This implies that for all local functions~$\varphi$ and all $\epsilon$,
\begin{equation}
\Vert \Psi+\epsilon\nabla\varphi\Vert^2_{\Lvec^2}\ge\Vert \Psi\Vert^2_{\Lvec^2}
\end{equation}
Expanding the left-hand side and taking~$\epsilon\to0$ yields $\langle\Psi,\nabla\varphi\rangle=0$ for all local functions, i.e., $\Psi\in(L_\nabla^2)^\perp$. By Lemma~\ref{lemma-sol}(1), $\Psi$ is $\cmss L_\omega$-harmonic.
\end{proofsect}

Obviously, the conditions (1-3) in the above proposition can be satisfied by $\Psi:=0$; it is thanks to \eqref{E:proj} that this can generally be excluded. (However, we could still have that~$\Psi$ is identically zero; see Exercise~\ref{E:Psi=0}.)
A question might also arise whether the function $\Psi$ is uniquely determined by the above properties. Biskup and Spohn~\cite{Biskup-Spohn} showed by fairly soft arguments that this is indeed the case. In fact, one even has a stronger statement:
\begin{equation}
\label{E:3.27}
\Lvec^2=L_\nabla^2\oplus\bigl\{A\Psi\colon A\in\text{GL}(\R,d)\bigr\},
\end{equation}
with $A\Psi(\omega,x)$ denoting the vector whose $i$-th Cartesian coordinate is given by $\sum_j a_{ij}\hate_j\cdot\Psi(\omega,x)$ where $A=(a_{ij})$. The position function and the potential fields thus generate all shift-covariant square-integrable ($\R^d$-valued) vector fields. (Notwithstanding, see Problem~\ref{pb4.10} for a very non-trivial generalization of this question.) Quastel~\cite{Quastel} has derived a similar result to \eqref{E:3.27} albeit with the use of Poincar\'e inequality and spectral-gap estimates.

\smallskip
It should be emphasized at this point that the above constructions have been quite standard --- albeit perhaps in different context and using different notations --- in various contributions dealing with homogenization theory. An application of these techniques to random walk in random environment was done somewhat independently in the Western school by Varadhan, Papanicolaou and coauthors and in the Russian school by Kozlov.

In particular, Kozlov's well-known paper~\cite{Kozlov} contains an extended version of the Weyl decomposition of vector fields --- which he calls forms --- into the sum of a gradient field, a harmonic field and a constant field which applies even in non-reversible situations. Apart from strong ellipticity, the main requirements for this decomposition in~\cite{Kozlov} are:
\begin{enumerate}
\item[(1)]
There is an measure~$\Q$ which is invariant for the Markov chain on environments and absolutely continuous with respect to~$\BbbP$.
\item[(2)]
The reciprocal value of the Radon-Nikodym derivative $\frac{\textd\Q}{\textd\BbbP}$ is in $L^1(\BbbP)$.
\end{enumerate}
While the absolute continuity of an invariant measure is usually somewhat challenging, it is the second condition that is invariably nearly impossible to check directly in any realistic (non-reversible) situation. We note that although Kozlov's paper is known to contain inconsistencies, it puts forward a number of good ideas and is thus a very recommended reading for anyone with interest in this subject.

\smallskip
The construction of the harmonic deformation can be performed rather seamlessly even in the case when $\pi_\omega(x)$ is zero at some vertices. What we need to assume is that there is a $\BbbP$-a.s. unique \emph{infinite} component $\CC_\infty$ of vertices with $\pi_\omega(x)>0$ such that the conditional measure
\begin{equation}
\BbbP_0(-)=\BbbP(-|0\in\CC_\infty),
\end{equation}
with expectation denoted by~$\E_0$,
satisfies the following variant of the ``usual conditions'':
\begin{enumerate}
\item[(1')]
$\BbbP_0(\pi_\omega(0)>0)=1$ (which holds trivially) and $\E_0\pi_\omega(0)<\infty$.
\item[(2')]
$\BbbP_0$ is irreducible in the sense that, for every~$x\in\Z^d$ with $\BbbP_0(x\in\CC_\infty)>0$,
\begin{equation}
\BbbP_0\bigl(\,\omega\colon\sup_{n\ge0}\cmss P_\omega^n(0,x)>0\bigl|x\in\CC_\infty\bigr)=1.
\end{equation}
\end{enumerate}
(Condition~(3) for measure~$\BbbP$ is not needed for now, the translation invariance of~$\BbbP$ suffices.)

\begin{hwproblem}
Suppose~that $\CC_\infty$ and $\BbbP_0$ are well defined and assume conditions~(1') and (2') above. Suppose also \eqref{E:square}. If~$\varphi_n$ is any minimizing sequence of the functional
\begin{equation}
\label{E:Dir-perc}
\varphi\mapsto \E_0\biggl(\,\sum_x\omega_{0,x}\bigl|x+\nabla_x\varphi(\omega)\bigr|^2\biggr),
\end{equation}
show that $\nabla\varphi_n(\omega,\cdot)$ still tends to some $\chi(\omega,\cdot)$ in $\Lvec^2$. Use this to define $\Psi=\Psi(\omega,x)$ with~$x\in\CC_\infty$ which is harmonic with respect to~$\cmss L_\omega$.
\end{hwproblem}

The function $\Psi$ constructed in this Exercise is the harmonic embedding of~$\CC_\infty$ that we discussed at the beginning of this section. A construction along the above lines can be found in the paper of Mathieu and Piatnitski~\cite{Mathieu-Piatnitski} for the problem of supercritical percolation cluster and in Biskup and Prescott~\cite{Biskup-Prescott} at the current level of generality. Berger and Biskup~\cite{Berger-Biskup} give a construction which is based on the spectral representation method of Kipnis and Varadhan (see end of Sect.~\ref{sec2.4}). Another way to define the corrector might be a result of:

\begin{hwproblem}
Show that the limit in
\begin{equation}
\lim_{n\to\infty}\bigl[ E_\omega^x(X_n)-E_\omega^0(X_n)\bigr]
\end{equation}
exists and equals $\Psi(\omega,x)$ for~$\BbbP$-a.e.~$\omega$.
\end{hwproblem}

It would be of much interest to find a solution to this problem without a recourse to the functional-analytic methods discussed above.

\subsection{Quenched Invariance Principle on deformed graph}
\noindent
Let us now turn attention back to the problem of a random walk among random conductances. A simple consequence of the above constructions is:

\begin{corollary}
\label{cor-def}
Suppose~$\BbbP$ satisfies the ``usual conditions'' and, in addition, \eqref{E:square} holds. Define $M_n:=\Psi(\omega,X_n)$. Then for $\BbbP$-a.e.~$\omega$ and each~$T>0$, the law of
\begin{equation}
\label{E:3.35a}
t\mapsto \frac1{\sqrt n}\bigl(M_{\lfloor tn\rfloor}+(tn-\lfloor tn\rfloor)(M_{\lfloor tn\rfloor+1}-M_{\lfloor tn\rfloor})\bigr)
\end{equation}
induced by~$P_\omega^0$ on the space $C([0,T]$, tends to the Brownian motion~$B_t$ with $EB_t=0$ and the covariance structure determined by
\begin{equation}
\label{E:Cov}
\frac1t E\bigl[(\lambda\cdot B_t)^2\bigr]=E_\Q E_\omega^0\bigl[[\lambda\cdot\Psi(\omega,X_1))^2\bigr],\qquad\lambda\in\R^d.
\end{equation}
\end{corollary}

\begin{proofsect}{Proof}
By the Cram\'er-Wold device it suffices to prove the convergence in law for the projection of the process onto any vector. We will denote this projection (with some abuse of notation) also by $M_n:=\lambda\cdot\Psi(\omega,X_n)$. The filtration is as before: $\scrF_n=\sigma(X_0,\dots,X_n)$.

First, the $\cmss L_\omega$-harmonicity of $\Psi$ guarantees that $M_n$ is a martingale so we just need to verify the conditions (LF1-LF2) of the Martingale Functional CLT. We will take care of both of these by considering the function
\begin{equation}
f_K(\omega):=E_\omega^0\bigl(|M_1|^2\1_{\{|M_1|\ge K\}}\bigr).
\end{equation}
Indeed, by property~(3) in Proposition~\ref{prop-correct}, $f_K\in L^2(\Q)$ for all $K\ge0$. Next, the shift-covariance of $\Psi$ implies $M_{k+1}-M_k=\lambda\cdot\Psi(\tau_{X_k}\omega,X_{k+1}-X_k)$ and so, by the Markov property,
\begin{equation}
E_\omega^0\bigl(|M_{k+1}-M_k|^2\1_{|M_{k+1}-M_k|\ge K}\big|\scrF_k\bigr)=f_K(\tau_{X_k}\omega).
\end{equation}
It follows that the left-hand side of (LF1) equals
\begin{equation}
\label{E:erg-sum}
\frac1n\sum_{k=0}^{n-1}f_K(\tau_{X_k}\omega)
\end{equation}
for~$K:=0$, while the left-hand side of the expression in (LF2) is bounded by this term from above as soon as~$n$ is so large that $\epsilon\sqrt n>K$.

Ergodicity of~$\BbbP$ with respect to translations ensures via \eqref{E:ergodic-1arg} that the expression \eqref{E:erg-sum} tends to $E_\Q f_K(\omega)$ as $n\to\infty$. This verifies (LF1) with $\sigma^2$ given by the right-hand side of \eqref{E:Cov}, and it also proves (LF2) because, thanks to the Dominated Convergence Theorem, we have
\begin{equation}
\lim_{K\to\infty}E_\Q f_K(\omega)=0.
\end{equation}
The result now follows by applying Theorem~\ref{thm-MCLT}.
\end{proofsect}

The above argument can be pushed through even in the case when the walk is restricted to an infinite connected component~$\CC_\infty$, as described above. One just needs to carefully check that the current proof of Proposition~\ref{prop-ergodic} still applies (details are spelled out in Berger and Biskup~\cite{Berger-Biskup}). However, later arguments might be seriously hampered by the fact that~$\BbbP_0$ is no longer shift invariant. This can be circumvented by the introduction of an \emph{induced shift}. Namely, for each~$i=1,\dots,d$, let
\begin{equation}
\theta_i\omega:=\tau_{n_i(\omega)\hate_i}\omega
\end{equation}
where
\begin{equation}
n_i(\omega):=\inf\bigl\{n\ge1\colon n\hate_i\in\CC_\infty(\omega)\bigr\}.
\end{equation}
The collection of maps $(\theta_1,\dots,\theta_d)$ defines shifts which preserve~$\BbbP_0$ and, in fact, make~$\BbbP_0$ ergodic. To see why these are well defined and the last property is true, consider the following exercise from abstract ergodic theory:

\begin{hwproblem}
Let $(\scrX,\scrF,\mu)$ be a probability space and let~$A\in\scrF$ be such that $\mu(A)>0$. Let $\tau\colon\scrX\to\scrX$ be a $\mu$-preserving bijection and suppose that $\mu$ is ergodic with respect to~$\tau$. Let $n_A(x):=\inf\{n\ge1\colon\tau^n(x)\in A\}$ for each~$x\in\scrX$. Do the following:
\begin{enumerate}
\item[(1)]
Show that $n_A<\infty$ $\mu$-a.s.
\end{enumerate}
This permits us to define $\theta(x):=\tau^{n_A(x)}(x)$. Next:
\begin{enumerate}
\item[(2)]
Show  that $\theta(A)=A$ $\mu$-a.s.\ and that $\theta$ preserves~$\mu_A(-):=\mu(-|A)$.
\item[(3)]
Prove that $\mu_A$ is ergodic with respect to~$\theta$.
\end{enumerate}
\end{hwproblem}

\smallskip
We will close this section with an exercise that illustrates the above abstract setting in one situation where explicit calculations are possible.

\begin{hwproblem}
\label{HW:d=1}
Suppose $d=1$ and only nearest-neighbor conductances. Assume that $\BbbP$ is ergodic with respect to the canonical shift on~$\Z$ and suppose that
\begin{equation}
\E(\omega_{0,1})<\infty\quad\text{\rm and}\quad C^{-1}:=\E\Bigl(\frac1{\omega_{0,1}}\Bigr)<\infty.
\end{equation}
 Verify that
\begin{equation}
\label{E:psi-d=1}
\Psi(\omega,x)=\begin{cases}
\displaystyle C\sum_{i=0}^{x-1}\frac1{\omega_{i,i+1}},\qquad&\text{if }x>0,
\\
\displaystyle-C\sum_{i=x}^{-1}\frac1{\omega_{i,i+1}},\qquad&\text{if }x<0,
\end{cases}
\end{equation}
defines a function satisfying properties~(1-3) in Proposition~\ref{prop-correct}. Conclude that the random walk $(\Psi(\omega,X_n))_{n\ge0}$ satisfies the (quenched) invariance principle.
\end{hwproblem}

We remark that the one-dimensional Random Conductance Model have quite intensely been studied, e.g., by Comets and Popov~\cite{Comets-Popov}, Gallesco and Popov~\cite{Gallesco-Popov}, Gallesco, Gantert, Popov and Vachovskaia~\cite{GGPV}, etc. A related problem is that of the random walk on random trees (with or without random conductances); see e.g., Lyons, Pemantle and Peres~\cite{Lyons-Pemantle-Peres}, Peres and Zeitouni~\cite{Peres-Zeitouni}, Gantert, M\"uller, Popov and Vachovskaia~\cite{GMPV}.


\section{Taming the deformation}
\label{sec4}\noindent
In this section our main goal is to finish the discussion of the essential steps of the proof of the quenched invariance principle. We will do this while leaving the most technically involved part, heat-kernel estimates, to the next section. Most of the material discussed here is quite standard; a possible exception is Theorem~\ref{thm-2D-sublinear} which has not appeared in this generality before.

\subsection{Remaining issues}
\noindent
Let us quickly review what we have accomplished so far. First, we used the examples of the balanced environments to isolate the martingale property as the key vehicle that  will get us to the CLT (Section~\ref{sec2.3}). Then, in the situations which are not balanced, we introduced a new embedding of~$\Z^d$ --- described by the function~$\Psi$ above --- that again makes the random walk into a martingale (Proposition~\ref{prop-correct}). On this embedding we succeeded in proving the convergence to Brownian motion (Corollary~\ref{cor-def}). However, two issues remained unresolved:
\begin{enumerate}
\item[(1)] The limiting Brownian motion may be degenerate to a point.
\item[(2)] The harmonic embedding may be quite distorted from the original lattice.
\end{enumerate}
Although the answer to (1) is ultimately related to the answer to (2), we will first focus on (1) as it is easier. We will start by solving Exercise~\ref{HW:d=1}.

It is easy to check that the function $\Psi$ from \eqref{E:psi-d=1} is harmonic with respect to $\cmss L_\omega$. This follows from the calculation
\begin{equation}
\begin{aligned}
\cmss L_\omega\Psi(\omega,x)
&=\omega_{x,x+1}\bigl[\Psi(\omega,x+1)-\Psi(\omega,x)\bigr]
\\
&\qquad\qquad\qquad\qquad\qquad+\omega_{x-1,x}\bigl[\Psi(\omega,x-1)-\Psi(\omega,x)\bigr]
\\
&=\omega_{x,x+1}\frac{C}{\omega_{x,x+1}}+\omega_{x-1,x}\Bigl(-\frac C{\omega_{x-1,x}}\Bigr)=C-C=0.
\end{aligned}
\end{equation}
The shift-covariance is a consequence of the additive form of the expressions in \eqref{E:psi-d=1} while integrability follows from
\begin{equation}
\begin{aligned}
E_{\Q}E_\omega^0\Psi(\omega,X_1)^2
&=\frac1Z\,\E\biggl[\omega_{0,1}\Bigl(\frac C{\omega_{0,1}}\Bigr)^2+\omega_{-1,0}\Bigl(\frac C{\omega_{-1,0}}\Bigr)^2\biggr]
\\
&=\frac1{\E(\omega_{0,1})\E(1/\omega_{0,1})},\end{aligned}
\end{equation}
which is finite and \emph{positive} by our assumptions. Applying the arguments in the proof of Corollary~\ref{cor-def}, $\Psi(\omega,X_n)$ satisfies an invariance principle with a non-degenerate limiting Brownian motion. The remainder of the Exercise is now embedded into:

\begin{proposition}
\label{prop-positive}
Suppose $\BbbP$ obeys the ``usual assumptions'' and \eqref{E:square}. In addition, assume that
\begin{equation}
\label{E:4.3}
\E\bigl(1/\omega_{0,\hate_i}\bigr)<\infty,\quad i=1,\dots,d.
\end{equation}
Then the limiting Brownian motion in Corollary~\ref{cor-def} is non-degenerate.
\end{proposition}

\begin{proofsect}{Proof}
We need to show that the right-hand side of \eqref{E:Cov} is bounded below by $c|\lambda|^2$ for some $c>0$ and all~$\lambda\in\R^d$. To this end we write
\begin{equation}
\label{E:4.4}
\begin{aligned}
E_\Q E_\omega^0\bigl[[\lambda\cdot\Psi(\omega,X_1)]^2\bigr]
&=\inf_\varphi\frac1Z\E\biggl(\sum_x\omega_{0,x}\bigl[\lambda\cdot(x+\nabla_x\varphi(\omega))]^2\biggr)
\\
&\ge\frac1Z\,\sum_{i=1}^d\inf_\varphi\,\E\biggl(\sum_{x=\pm\hate_i}\omega_{0,x}\bigl[\lambda\cdot(x+\nabla_x\varphi(\omega))]^2\biggr).
\end{aligned}
\end{equation}
The expectation on the extreme right now involves only edges in the $i$-th coordinate direction and thus effectively becomes a one-dimensional problem. To overcome a possible lack of separate ergodicity, let $\scrA_i$ be the $\sigma$-algebra of $\tau_{\hate_i}$-invariant events and let $\lambda_i:=\lambda\cdot\hate_i$. The second infimum in \eqref{E:4.4} is then bounded below by $2\lambda_i^2\E( 1/\E(\omega_{0,\hate_i}^{-1}|\scrA_i))$ which by Jensen's inequality is at most $2\lambda_i^21/E(1/\omega_{0,\hate_i})$. It follows that
\begin{equation}
E_\Q E_\omega^0\bigl[[\hate_i\cdot\Psi(\omega,X_1)]^2\bigr]\ge\frac2{\E(\pi_\omega(0))}\sum_{i=1}^d\frac{\lambda_i^2}{\E(1/\omega_{0,1})}.
\end{equation}
By \eqref{E:square} and \eqref{E:4.3} we conclude that this exceeds $c|\lambda|^2$ for some $c>0$.
\end{proofsect}

Note that the same argument would apply whenever the set
\begin{equation}
\label{E:gen-set}
\bigl\{x\in\Z^d\colon\E(1/\omega_{0,x})<\infty\bigr\}
\end{equation}
generates all of~$\Z^d$ (as an additive group). This still does not cover the case of supercritical percolation (which can nonetheless be covered by an alternate argument) so we pose:

\begin{problem}
Give general (and natural) conditions under which the limiting Brownian motion is non-degenerate. In particular, if $R(x,y)$ denotes the effective resistivity \eqref{E:1.20} from~$x$ to~$y$ in environment $\omega$, does it suffice to replace $\E(1/\omega_{0,x})<\infty$ in \eqref{E:gen-set} by $\E R(0,x)$?
\end{problem}

Note that we have a pointwise bound $R(0,x)\le1/\omega_{0,x}$ with $R(0,1)=1/\omega_{0,1}$ in $d=1$ with nearest-neighbor conductances. This suggests also:

\begin{hwproblem}
\label{E:Psi=0}
Suppose $d=1$ and let $\BbbP$ be a measure on i.i.d.~positive and nearest-neighbor conductances such that $\E(\omega_{0,1})<\infty$ and~$\E(1/\omega_{0,1})=\infty$. Show that the infimum of \eqref{E:Dir} over local functions is zero. Conclude that we must have $\chi(\omega,x)=-x$.
\end{hwproblem}

A proof of (an analogue of) Proposition~\ref{prop-positive} appeared in Kozlov~\cite{Kozlov} and in de Masi, Ferrari, Goldstein and Wick~\cite{demas1,demas2}. With a bit more effort one can develop a variational characterization of the \emph{inverse} of the limiting covariance matrix by minimizing a (version of) Dirichlet energy over nearly linear flows (Biskup~\cite{Biskup-grad}). This in principle allows one to numerically approximate the covariance matrix with arbitrary precision from above and below.

Approximation arguments for the diffusion constants are at the core of the Kipnis-Varadhan approach sketched in Sect.~\ref{sec2.4}. Caputo and Ioffe~\cite{Caputo-Ioffe} studied periodized versions of the Random Conductance Model and the convergence of the effective diffusion coefficient to the infinite volume object; related work in a continuum context can be found in Owhadi~\cite{Owhadi}.

\subsection{Sublinearity of the corrector}
\noindent
Having addressed non-degeneracy of the limiting Brownian motion, we are ready to move to the second --- and considerably more involved --- issue. The important thing is to realize that for our purposes it would suffice to show that
\begin{equation}
\label{E:4.6}
\chi(\omega,X_n)=\Psi(\omega,X_n)-X_n=o(X_n)
\end{equation}
asymptotically along a typical path of the random walk. Indeed, once we know that $\Psi(\omega,X_k)-X_k=o(X_k)$ we can use the martingale CLT to get $\Psi(\omega,X_k)=o(\sqrt n)$ for all~$k\le n$ which then implies that also $X_k=O(\sqrt n)$. But then we will have $\Psi(\omega,X_k)-X_k=o(\sqrt n)$ for all $k\le n$, which means that the change of embedding of the graph has a vanishing effect at the diffusive scale.

A more general version of \eqref{E:4.6} would be to require this for all positions in the lattice, not just those visited by the path. In $d=1$, this is not hard to get:

\begin{hwproblem}
\label{ex4.2}
Suppose that $\BbbP$ is an ergodic law on nearest-neighbor conductances in $d=1$. Assume $\E(\omega_{0,1})<\infty$ and~$\E(1/\omega_{0,1})<\infty$. Show that
\begin{equation}
\Psi(\omega,x)-x=o(|x|),\qquad |x|\to\infty,
\end{equation}
and prove that the corresponding random walk satisfies a quenched invariance principle. (Compare also with Exercise~\ref{E:Psi=0}.)
\end{hwproblem}

However, the situation in higher dimensions is quite more subtle. While the technical details of derivations in the paper of Kipnis and Varadhan follow a different route, their methods can be used to show:

\begin{theorem}
\label{thm-4.5}
Under the ``usual assumptions,'' \eqref{E:square} and \eqref{E:4.3}, for each~$\epsilon>0$,
\begin{equation}
\label{E:4.9a}
E_\Q\,P_\omega^0\bigl(\,\max_{k\le n}|\chi(\omega,X_k)|>\epsilon\sqrt n\bigr)\,\underset{n\to\infty}\longrightarrow\,0.
\end{equation}
\end{theorem}

This statement will imply the so called \textit{Annealed Invariance Principle}, sometimes also called a \textit{functional CLT in probability}. We will choose to formulate this in a form of a coupling. Here we recall that, given two probability measures $P$ and~$P'$, their \emph{coupling} is a probability measure~$Q$ on the product space whose first, resp., second marginal is given by $P$, resp.,~$P'$.

\begin{corollary}[Annealed Invariance Principle]
\label{cor-AIP}
Given a path $X=(X_n)$, a time $t\ge0$ and $n\in\N$ let
\begin{equation}
W_t^{(n)}:=\frac1{\sqrt n}\bigl(X_{\lfloor tn\rfloor}+(tn-\lfloor tn\rfloor)(X_{\lfloor tn\rfloor+1}-X_{\lfloor tn\rfloor})\bigr).
\end{equation}
Under the assumptions of Theorem~\ref{thm-4.5}, for~$\BbbP$-a.e.~$\omega$, there exists a coupling $Q_\omega^0$ of the law of $t\mapsto W_t^{(n)}$ induced by~$P_\omega^0$ and a Brownian motion $t\mapsto B_t$ with mean zero and covariance \eqref{E:Cov} so that, for each $T>0$ and each $\epsilon>0$,
\begin{equation}
\label{E:4.11a}
\E\,Q_\omega^0\biggl(\,\sup_{0\le t\le T}|W_t^{(n)}-B_t|>\epsilon\biggr)\,\underset{n\to\infty}\longrightarrow\,0.
\end{equation}
\end{corollary}

\begin{proofsect}{Proof} (Sketch)
First let us note that both \eqref{E:4.9a} and \eqref{E:4.11a} hold equivalently with expectation $E_\Q$ or expectation~$\E$. (This is because $\Q$ and $\BbbP$ are equivalent and the quantity under expectation is bounded.) To prove \eqref{E:4.11a}, we will use the fact, implied by the Skorohod embedding, that such a coupling exists between the Brownian motion and the analogue of $t\mapsto W_t$ defined using the martingale $M_n:=\Psi(\omega,X_n)$. Let $Z_t^{(n)}$ denote the expression on the right of \eqref{E:3.35a}. Then we have
\begin{equation}
\label{E:4.11b}
\E\,Q_\omega^0\bigl(\,\sup_{0\le t\le T}|Z_t^{(n)}-B_t|>\epsilon\bigr)\,\underset{n\to\infty}\longrightarrow\,0
\end{equation}
where $Q_\omega^0$ is induced by the Skorohod embedding. As to \eqref{E:4.11a}, we note that
\begin{equation}
\sup_{0\le t\le T}|W_t^{(n)}-Z_t^{(n)}|\le \frac1{\sqrt n}\,\max_{k\le\lfloor Tn\rfloor+1}\bigl|\chi(\omega,X_k)\bigr|
\end{equation}
Since the event on the right does not depend on the second marginal of~$Q_\omega^0$, we thus have
\begin{equation}
\label{E:4.14a}
\E\,Q_\omega^0\biggl(\,\sup_{0\le t\le T}|W_t^{(n)}-Z_t^{(n)}|>\epsilon\biggr)
\le \E\,P_\omega^0\biggl(\,\max_{k\le\lfloor Tn\rfloor+1}\bigl|\chi(\omega,X_k)\bigr|>\epsilon\sqrt n\biggr)
\end{equation}
which tends to zero as~$n\to\infty$ by Theorem~\ref{thm-4.5}. Combining \twoeqref{E:4.11b}{E:4.14a} the result follows.
\end{proofsect}

We remark that when the supremum is dropped from \eqref{E:4.9a}, we talk about an \emph{annealed CLT}. The averaging over the invariant measure~$\Q$ in Theorem~\ref{thm-4.5} is not a mere technical convenience as the statement is  not strong enough to infer \eqref{E:4.11a} without the expectation over environment.
It actually took nearly 20 years after Kipnis-Varadhan's result before this issue was first successfully addressed and a proper \emph{quenched} invariance principle proved. This was done in the work of Sidoravicius and Sznitman~\cite{Sidoravicius-Sznitman} who realized that one can get further with the help of the \emph{heat kernel estimates}. However, Berger and Biskup~\cite{Berger-Biskup} were later able to avoid the use of these in their argument for the two-dimensional supercritical percolation cluster. We will present a sketch of their argument in a slightly more general, albeit non-percolative, setting:

\begin{theorem}
\label{thm-2D-sublinear}
Let~$d=2$ and suppose~$\BbbP$ is an ergodic law on nearest-neighbor conductances subject to the conditions $\E(\omega_{0,\hate_i})<\infty$ and $\E(1/\omega_{0,\hate_i})<\infty$ for $i=1,2$. Then
\begin{equation}
\lim_{n\to\infty}\,\frac1n\,{\max_{|x|\le n}|\chi(\omega,x)|}=0,\qquad\BbbP\text{\rm-a.s.}
\end{equation}
\end{theorem}

Our proof of Theorem~\ref{thm-2D-sublinear} begins by a lemma that generalizes Exercise~\ref{ex4.2} --- and that even in $d=1$, when the conductances are no longer just nearest neighbor --- to all dimensions:

\begin{lemma}[Directional sublinearity]
\label{lemma-direct-sub}
Suppose $d\ge1$ and assume $\BbbP$ is an ergodic law subject to the restriction \eqref{E:square}. Assume $\E(1/\omega_{0,\hate})<\infty$ for some~$\hate\in\{\pm\hate_i\colon i=1,\dots,d\}$. Then
\begin{equation}
\label{E:dir-sub}
\lim_{n\to\infty}\frac{\chi(\omega,n\hate)}n=0,\qquad\BbbP\text{\rm-a.s.}
\end{equation}
\end{lemma}

Before we set out to prove this, we note that there is a small technical subtlety that arises from the distinction between ergodicity and directional ergodicity. To make this distinction clearer, we invite the reader to first solve:

\begin{hwproblem}
Construct a law~$\BbbP$ on nearest-neighbor conductances that is (jointly) ergodic with respect to translations --- i.e., $\BbbP(A)=0$ for all $A$ with $\tau_x(A)=A$ for all~$x\in\Z^d$ --- but not separately ergodic in the sense that there is a set $B$ of environments which is invariant under translations in the first coordinate direction and for which $0<\BbbP(B)<1$.
\end{hwproblem}

\begin{proofsect}{Proof of Lemma~\ref{lemma-direct-sub}}
Using shift covariance we get
\begin{equation}
\chi(\omega,n\hate)=\sum_{k=0}^{n-1}\chi(\tau_{k\hate}\omega,\hate)
\end{equation}
We would like to use the (pointwise) Ergodic Theorem and $\tau_{\hate}$-invariance of~$\BbbP$ to extract the limit
\begin{equation}
f_{\hate}(\omega):=\lim_{n\to\infty}\frac{\chi(\omega,n\hate)}n
\end{equation}
and prove that it vanishes $\BbbP$-a.s. For that we will need to establish three things:
\begin{enumerate}
\item[(1)] $\E|\chi(\omega,\hate)|<\infty$.
\item[(2)] $\E\chi(\omega,\hate)=0$.
\item[(3)] $f_{\hate}$ is translation invariant.
\end{enumerate}
The first two items will follow from the construction of the corrector. Recall that we are guaranteed that $\nabla\varphi_n\to\chi$ in $\Lvec^2$ --- which is a kind of weighted $L^2$-space. Since for $\varphi\in L^\infty(\BbbP)$,
\begin{equation}
\E\nabla_{\hate}\varphi(\omega)=\E\varphi\circ\tau_{\hate}-\E\varphi=0
\end{equation}
it suffices to show that $\chi(\cdot,\hate)\in L^1$ and $\nabla_{\hate}\varphi_n\to\chi(\cdot,\hate)$ in $L^1$. (The former actually follows from the latter, but we find this order more instructive.) And, indeed, by the Cauchy-Schwarz inequality we get
\begin{equation}
\bigl[\E|\chi(\omega,\hate)|\bigr]^2\le\E\Bigl(\frac1{\omega_{0,\hate}}\Bigr)\E\Bigl(\omega_{0,\hate}|\chi(\omega,\hate)|^2\Bigr)\le\E\Bigl(\frac1{\omega_{0,\hate}}\Bigr)\Vert\chi\Vert_{\Lvec^2}^ 2<\infty
\end{equation}
and similarly we derive
\begin{equation}
\bigl[\E|\chi(\omega,\hate)-\nabla_{\hate}\varphi_n(\omega)|\bigr]^2\le\E\Bigl(\frac1{\omega_{0,\hate}}\Bigr)\Vert\chi-\nabla\varphi_n\Vert_{\Lvec^2}^2
\end{equation}
which tends to zero as $n\to\infty$ because $\nabla\varphi\to\chi$ in $\Lvec^2$.

Finally, in order to link the limit to the expectation, we also need to show that $f_{\hate}$ is translation invariant. To that end pick another lattice direction $\hate'$ and note that, by translation covariance,
\begin{equation}
\chi(\tau_{\hate'}\omega,n\hate)=\chi(\omega,n\hate)+\chi(\tau_{n\hate}\omega,\hate')-\chi(\omega,\hate').
\end{equation}
Dividing by~$n$, the $L^1$-limit of the last two terms is zero and so from the above $L^1$-inclusions we conclude that $f_{\hate}(\tau_{\hate'}\omega)=f_{\hate}(\omega)$ for $\BbbP$-a.e.~$\omega$. Putting all pieces together the claim follows.
\end{proofsect}

We remark that the fact that the conditions in Lemma~\ref{lemma-direct-sub} are the same as in Proposition~\ref{prop-positive} is not a coincidence. Indeed we have:

\begin{hwproblem}[Sublinearity implies nondegeneracy]
Show that if \eqref{E:dir-sub} holds for vector~$\hate$, then the component of the limiting Brownian motion  --- constructed, at this point, by the Martingale Convergence Theorem --- in direction of~$\hate$ is non-degenerate.
\end{hwproblem}

Our next goal is to boost the directional subadditivity --- which we may assume for both lattice directions under the conditions of Theorem~\ref{thm-2D-sublinear} --- into a corresponding statement over a box of side~$n$. To this end, let us say that the origin is {$(K,\epsilon)$-good in~$\omega$} if for all $\hate\in\{\pm\hate_i\colon i=1,2\}$ and all $n\ge1$,
\begin{equation}
\bigl|\chi(\omega,n\hate)\bigr|\le K+\epsilon n.
\end{equation}
A point~$x$ is then called $(K,\epsilon)$-good in~$\omega$ if $0$ is $(K,\epsilon)$-good in $\tau_x\omega$. By Lemma~\ref{lemma-direct-sub} we know that
\begin{equation}
\BbbP\bigl(0\text{ is } (K,\epsilon)\text{-good}\bigr)\,\underset{K\to\infty}\longrightarrow\,1.
\end{equation}
It is now an exercise to show that:

\begin{hwproblem}
Fix~$\epsilon>0$. Show that for each~$\rho\in(0,1)$ and for $\BbbP$-a.e.\ $\omega$ there is a number $K=K(\rho,\omega)<\infty$ such that
\begin{enumerate}
\item[(1)]
$0$ is $(K,\epsilon)$-good in $\omega$.
\item[(2)]
The density of $(K,\epsilon)$-good vertices on the lines $\{n\hate_i\colon n\in\Z\}$, $i=1,2$, is at least~$\rho$.
\end{enumerate}
\end{hwproblem}

\begin{figure}[t]
\centerline{\includegraphics[width=0.8\textwidth]{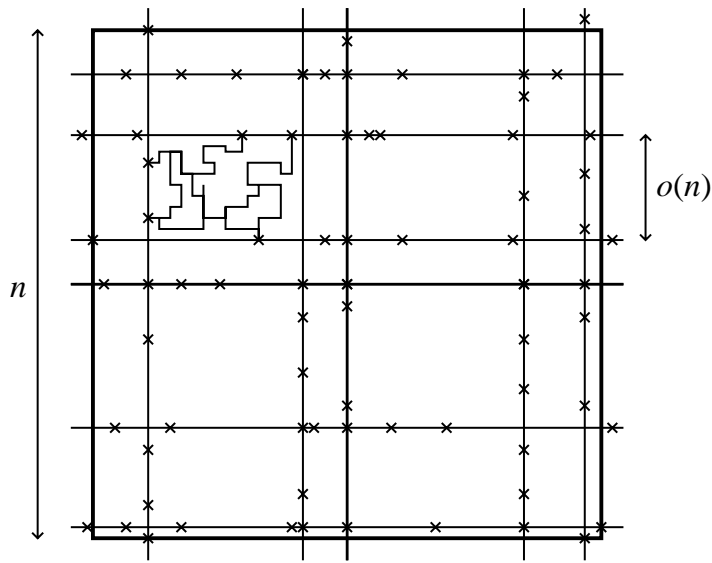}}
\caption{A figure illustrating the last part of the proof of Theorem~\ref{thm-2D-sublinear} where one needs to control the corrector on a component away from the good grid by the values on its boundary. The existence and  positivity of the densities of the good vertices along coordinate directions ensures that the largest such component intersecting the box $[-n,n]^2$ is $o(n)$ in diameter. The argument relies on the maximum principle for $x\mapsto\Psi(\omega,x)$.}
\label{fig5}
\end{figure}

These observations permit us to define a \emph{good grid} as follows. Take the two lines $\{n\hate_i\colon n\in\Z\}$, $i=1,2$, and add to them all vertices of the form $n_1\hate_1+n_2\hate_2$ with $n_1,n_2\in\Z$ such that either $n_1\hate_1$ or $n_2\hate_2$ is $(K,\epsilon)$-good. Call the resulting (random) set of vertices $\G_{K,\epsilon}(\omega)$. Then we note:

\begin{lemma}
\label{lemma-good-grid}
We have:
\begin{equation}
\max_{\begin{subarray}{c}
|x|_\infty\le n\\x\in\G_{K,\epsilon}(\omega)
\end{subarray}}
\bigl|\chi(\omega,x)\bigr|\le 2K+2\epsilon n.
\end{equation}
\end{lemma}

\begin{proofsect}{Proof}
Let $x:=n_1\hate_1+n_2\hate_2$ be a vertex in~$\G_{K,\epsilon}(\omega)$. This means that, e.g., $n_1\hate_1$ is $(K,\epsilon)$-good in~$\omega$. Since the origin is $(K,\epsilon)$-good as well, we can write
\begin{equation}
\begin{aligned}
\bigl|\chi(\omega,x)\bigr|
&\le\bigl|\chi(\omega,x)-\chi(\omega,n_1\hate_1)\bigr|+\bigl|\chi(n_1\hate_1,x)\bigr|
\\
&\le K+\epsilon |n_2|+K+\epsilon|n_2|
\end{aligned}
\end{equation}
But $|x|_\infty\le n$ implies $|n_1|,|n_2|\le n$ and so the claim follows.
\end{proofsect}

We now know how to control the corrector at the vertices of the good grid --- which can be made arbitrary dense --- but we still have to worry about those in the complement thereof. An important fact is that the connected components of~$\Z^2\setminus\G_{K,\epsilon}(\omega)$ are finite and, in fact, that any such component intersecting the box $[-n,n]^2$ has diameter $o(n)$. This can be justified by solving:

\begin{hwproblem}
Consider any shift invariant, ergodic, zero-one valued process on~$\Z$ with the densities of $0$'s and $1$'s both positive. Show that in almost-every sample, the size of the largest consecutive block of $1$'s intersecting the interval $[-n,n]$ is only $o(n)$ as~$n\to\infty$.
\end{hwproblem}

We can now finish the proof of sublinearity of the corrector:

\smallskip
\begin{proofsect}{Proof of Theorem~\ref{thm-2D-sublinear}}
Pick~$x\in\Z^2\setminus\G_{K,\epsilon}(\omega)$ with $|x|\le n$. Let~$C(x)$ denote the component containing~$x$. We claim that
\begin{equation}
\bigl|\chi(\omega,x)\bigr|\le \text{diam}\,C(x)+\max_{z\in\partial C(x)}
\bigl|\chi(\omega,z)\bigr|
\end{equation}
This is a consequence of $\cmss L_\omega$-harmonicity of $\Psi$ and the maximum principle. Indeed, define the first hitting time
\begin{equation}
T:=\inf\{n\ge0\colon X_n\not\in C(x)\}
\end{equation}
 of the complement of~$C(x)$. Then
\begin{equation}
\label{E:4.22}
\Psi(\omega,x)=E_\omega^x\bigl(\Psi(\omega,X_T)\bigr),
\end{equation}
 which we can rewrite as
\begin{equation}
\label{E:4.22b}
\chi(\omega,x)=E_\omega^x\bigl(X_T-x\bigr)+E_\omega^x\bigl(\chi(\omega,X_T)\bigr)
\end{equation}
But $|X_T-x|\le \text{diam}\,C(x)$ and $|\chi(\omega,X_T)|\le \max_{z\in\partial C(x)}
\bigl|\chi(\omega,z)\bigr|$ so the bound follows.

To finish the argument, we recall that $\text{diam}\,C(x)=o(n)$ and so we may assume that~$n$ is so large that $C(x)\subset[-2n,2n]^2$. In that case $\max_{z\in\partial C(x)}
\bigl|\chi(\omega,z)\bigr|
$ is bounded by the maximum from Lemma~\ref{lemma-good-grid} with~$n$ replaced by~$2n$. We get
\begin{equation}
\bigl|\chi(\omega,x)\bigr|\le 2K+4\epsilon n+o(n)
\end{equation}
thus proving the claim.
\end{proofsect}

\subsection{Above two dimensions}
\noindent
The above reasoning can be boosted to cover all ergodic two-dimensional environments with a finite range of jumps that satisfy the condition $\E(1/\omega_{0,\hate_i})<\infty$. However, there is an inherent problem with this approach in higher dimension; indeed, one can still define a good grid but this grid will no longer partition $\Z^d$ into finite components. In an attempt to adapt the argument based on \twoeqref{E:4.22}{E:4.22b}, one thus has to worry about two things: How long does it take to hit the good grid and how far will $X_T$ be from~$x$. This can be done but (insofar) only with the help of heat-kernel technology. We paraphrase a theorem from Biskup and Prescott~\cite{Biskup-Prescott}:

\begin{theorem}
\label{thm-BP}
Fix~$\omega$ such that $\pi_\omega(x)\in(0,\infty)$ for all~$x$ and suppose~$\chi=\chi(x)$ is a function and $\theta>0$ is a number such that the following holds:
\settowidth{\leftmargini}{(111a)}
\begin{enumerate}
\item[(1)] (Harmonicity)
The function~$\Psi(x):=x+\chi(x)$ obeys $\cmss L_\omega\Psi(x)=0$ for all~$x$.
\item[(2)] (Sublinearity on average)
For every~$\epsilon>0$,
\begin{equation}
\label{on-average}
\lim_{n\to\infty}\frac1{n^d}\sum_{x\colon |x|\le n}
\1_{\{|\chi(x)|\ge\epsilon n\}}=0.
\end{equation}
\item[(3)] (Polynomial growth)
\begin{equation}
\label{polygrowth}
\lim_{n\to\infty}\,\max_{|x|\le n}
\frac{|\chi(x)|}{n^\theta}=0.
\end{equation}
\end{enumerate}
Let~$Y=(Y_t)$ be the variable-speed continuous-time random walk with generator~$\cmss L_\omega$ and suppose that the following estimates hold:
\begin{enumerate}
\item[(4)] (Diffusive upper bounds) For a sequence~$b_n=o(n^2)$,
\begin{equation}
\label{diffusive}
\sup_{n\ge1}\,\,\max_{|x|\le n}\,\,
\sup_{t\ge b_n}\,\frac{E_{\omega}^x|Y_t-x|}{\sqrt t}<\infty
\end{equation}
and
\begin{equation}
\label{on-diag}
\sup_{n\ge1}\,\,\max_{|x|\le n}\,\,
\sup_{t\ge b_n}\, t^{d/2}P_{\omega}^x(Y_t=x)<\infty.
\end{equation}
\end{enumerate}
Then
\begin{equation}
\label{sublinear}
\lim_{n\to\infty}\,\max_{|x|\le n}
\frac{|\chi(x)|}n=0.
\end{equation}
\end{theorem}

We remark that most of the proof of this theorem goes through even when the variable-speed random walk is replaced by the constant-speed walk (for which the bounds \twoeqref{diffusive}{on-diag} may be easier to prove). This is because $\Psi(X_t)$ is a martingale for both walks. The sole point where the variable speed walk seems to be used is formula (5.13) on page 1338 of~\cite{Biskup-Prescott}.

\smallskip
In an earlier work (e.g.,~Berger and Biskup~\cite[Appendix~A2]{Berger-Biskup}) the same conclusion as given by Theorem~\ref{thm-BP} could be achieved --- although perhaps in a less transparent way --- by using the full heat-kernel upper bounds of the form
\begin{equation}
\cmss P_\omega^n(x,y)\le\frac{c_1}{n^{d/2}}\texte^{-c_2|x-y|^2/n}.
\end{equation}
The point of reducing the heat-kernel input to the statements \twoeqref{diffusive}{on-diag} is that these are easier to verify than the actual heat-kernel upper bounds. We also note that Sidoravicius and Sznitman~\cite{Sidoravicius-Sznitman} have used the heat-kernel bounds mainly to control the tightness of the limiting process, while here we are using it to control the deformations of the harmonic embedding. (Tightness follows in our case from the Martingale Functional CLT.)

\smallskip
A key input in Theorem~\ref{thm-BP} is the sublinearity-on-average claim which we formalize as:

\begin{proposition}[Sublinearity on average]
\label{lemma-average}
Suppose $d\ge1$ and assume $\BbbP$ is an ergodic law subject to the restriction \eqref{E:square}. Assume $\E(1/\omega_{0,\hate})<\infty$ for all~$\hate\in\{\pm\hate_i\colon i=1,\dots,d\}$. Then for each~$\delta>0$,
\begin{equation}
\label{limit-on-average}
\lim_{n\to\infty}\,\frac1{n^d}\sum_{\begin{subarray}{c}
|x|\le n
\end{subarray}}
\1_{\{|\chi(x,\omega)|\ge\delta n\}}=0, \qquad \BbbP\text{\rm-a.s.}
\end{equation}
\end{proposition}

The proof is based on the commutative structure of~$\Z^d$ and a bootstrapping of the one-dimensional sublinearity established in Lemma~\ref{lemma-direct-sub} by induction along dimension. Recall the notion of a good grid~$\G_{K,\epsilon}$ introduced (in $d=2$) earlier. The induction argument is contained in the following deterministic ``pigeon-hole-principle'' lemma:

\begin{lemma}
\label{lemma-determine}
Let $\Lambda_n:=[-n,n]^d\cap\Z^d$, fix~$\epsilon>0$ and~$K<\infty$. For each~$\omega$ there exists a set~$A=A(\omega)\subset\G_{K,\epsilon}(\omega)\cap\Lambda_n$ with the properties
\begin{equation}
|A|\ge(2^d\eta-2^d+1)|\Lambda_n|\quad\text{where}\quad \eta:=\frac{|\G_{K,\epsilon}(\omega)\cap\Lambda_n|}{|\Lambda_n|}
\end{equation}
and
\begin{equation}
\label{2.13}
x,y\in A\quad\Rightarrow\quad \bigl|\chi(y,\omega)-\chi(x,\omega)\bigr|\le 2d\bigl[K+\epsilon(2n+1)\bigr].
\end{equation}
\end{lemma}

\begin{proofsect}{Proof}
We will prove this by induction on dimension. Fix~$\omega$ and for~$\nu\in\{1,\dots,d\}$ define sets $\Lambda_n^{(\nu)}$ of the form
\begin{equation}
\Lambda_n^{(\nu)}:=\Lambda_n\cap\bigcap_{j=\nu+1}^d\bigl\{x=(x_1,\dots,x_d)\colon x_j=m_j\bigr\}
\end{equation}
for some $m_2,\dots,m_d\in[-n,n]\cap\Z$ as follows: We define~$\Lambda_n^{(d)}:=\Lambda_n$ and if~$\Lambda^{(\nu+1)}_n$ has been defined, we use~$\Lambda_n^{(\nu)}$ to denote a $\nu$-dimensional set of the above form which contains the maximum number of good sites. Note that if~$\eta$ is as in the statement, we have
\begin{equation}
\label{eta-bd}
\frac{|\G_{K,\epsilon}(\omega)\cap\Lambda_n^{(\nu)}|}{|\Lambda_n^{(\nu)}|}\ge\eta
\end{equation}
because the ratio on the left decreases in~$\nu$.

Next we set~$\eta_\nu:=2^\nu\eta-2^\nu+1$ and note that
\begin{equation}
\label{eta-ind}
2\eta_\nu-1=\eta_{\nu+1}
\quad\text{and}\quad
\eta\ge\eta_\nu.
\end{equation}
Assuming without loss of generality that~$\eta>1-2^{-d}$ --- otherwise we can take~$A:=0$ in the statement of the lemma --- we have~$\eta_\nu>0$ for all~$1\le\nu\le d$.
We will prove by induction the following claim: For each~$\nu=1,\dots,d$, there exists a set $A^{(\nu)}\subset\G_{K,\epsilon}(\omega)\cap\Lambda_n^{(\nu)}$ such that
\begin{equation}
\label{A-bd}
|A^{(\nu)}|\ge\eta_\nu|\Lambda_n^{(\nu)}|
\end{equation}
and
\begin{equation}
\label{cor-bd1}
x,y\in A^{(\nu)}\quad\Rightarrow\quad \bigl|\chi(y,\omega)-\chi(x,\omega)\bigr|\le 2\nu\bigl[K+\epsilon(2n+1)\bigr].
\end{equation}
For $\nu=d$ this clearly implies the desired claim.

For~$\nu=1$, we define $A^{(1)}:=\G_{K,\epsilon}(\omega)\cap\Lambda_n^{(1)}$. As $\eta_1\le\eta$, this obeys \eqref{A-bd}. The bound \eqref{cor-bd1} is then a direct consequence of the definition of a good line. Suppose now that the claim holds for~$\nu$ and let us prove it for~$\nu+1$. To this extent, let $\Pi$ denote the natural projection of $\Lambda_n^{(\nu+1)}$ onto $\Lambda_n^{(\nu)}$ and, given the set~$A^{(\nu)}$ with the above properties, let
\begin{equation}
A^{(\nu+1)}:=\bigl\{x\in\G_{K,\epsilon}(\omega)\cap\Lambda_n^{(\nu+1)}\colon \Pi(x)\in A^{(\nu)}\bigr\}.
\end{equation}
We now verify that this $A^{(\nu+1)}$ obeys \twoeqref{A-bd}{cor-bd1}. As to \eqref{A-bd}, the same bound for $A^{(\nu)}$ tells us that at most $(1-\eta_\nu)|\Lambda_n^{(\nu+1)}|$ sites in $\Lambda_n^{(\nu+1)}$ do not project into~$A^{(\nu)}$. Hence
\begin{equation}
\begin{aligned}
|A^{(\nu+1)}|&\ge\bigl|\G_{K,\epsilon}(\omega)\cap\Lambda_n^{(\nu+1)}\bigr|-(1-\eta_\nu)|\Lambda_n^{(\nu+1)}|
\\
&\ge(\eta+\eta_\nu-1)|\Lambda_n^{(\nu+1)}|,
\end{aligned}
\end{equation}
where we used \eqref{eta-bd} to get the second inequality. In light of \eqref{eta-ind} this implies \eqref{A-bd}.

To prove also \eqref{cor-bd1}, we pick two sites $x,y\in A^{(\nu+1)}$ and let $\bar x=\Pi(x)$ and $\bar y=\Pi(y)$. The claim for~$\nu$ then implies
\begin{equation}
\bigl|\chi(\bar y,\omega)-\chi(\bar x,\omega)\bigr|\le 2\nu\bigl[K+\epsilon(2n+1)\bigr]
\end{equation}
while the fact that~$x$ is a good site yields
\begin{equation}
\bigl|\chi(x,\omega)-\chi(\bar x,\omega)\bigr|\le K+\epsilon(2n+1)
\end{equation}
and similarly for the pair~$y$ and~$\bar y$. Combining these bounds and using the triangle inequality then implies \eqref{cor-bd1} for~$x$ and~$y$ --- with, of course, $2\nu$ replaced by $2(\nu+1)$.
\end{proofsect}

Lemma~\ref{lemma-determine} now implies that the corrector is sublinear on average:

\smallskip
\begin{proofsect}{Proof of Proposition~\ref{lemma-average}}
Suppose without loss of generality that~$\delta<8^{-d}$, fix $\epsilon<\frac1{32d}\delta$ and note that we can choose~$K$ so large that $\BbbP(0\in\G_{K,\epsilon})\ge1-\ffrac\delta2$. By the Spatial Ergodic Theorem and ergodicity of~$\BbbP$ we thus have
\begin{equation}
\frac{|\G_{K,\epsilon}(\omega)\cap\Lambda_n|}{|\Lambda_n|}\ge1-\delta
\end{equation}
once~$n\ge n_0$ for some a.s.\ finite~$n_0=n_0(\omega)$. We will assume that $n_0$ is so large that also
\begin{equation}
\label{delta}
\delta n> 16d\bigl[K+\epsilon(2n+1)\bigr]
\end{equation}
holds for all~$n\ge n_0$.

By Lemma~\ref{lemma-determine}, for each~$n\ge n_0$ there exists~$A_n=A_n(\omega)\subset\Lambda_n$ with
\begin{equation}
|A_n|\ge(1-2^d\delta)|\Lambda_n|
\end{equation}
and \eqref{2.13} valid for all~$x,y\in A_n$. As $\delta<8^{-d}$, we have $|A_n|\ge(1-4^{-d})|\Lambda_n|$ while $|\Lambda_{2n}\setminus A_{2n}|\le 4^{-d}|\Lambda_{2n}|=2^{-d}|\Lambda_n|$. In particular, $A_n\cap A_{2n}\ne\emptyset$ for each~$n\ge n_0$. Let~$k_0$ be the smallest integer such that $2^{k_0}\ge n_0$ and let us pick a site $x_k\in A_{2^k}\cap A_{2^{k+1}}$ for each~$k\ge k_0$. The bounds \eqref{2.13} and \eqref{delta} then give us
\begin{equation}
\bigl|\chi(x_k,\omega)-\chi(x_{k_0},\omega)\bigr|
\le\sum_{\ell=k_0}^{k-1}\bigl|\chi(x_{\ell+1},\omega)-\chi(x_\ell,\omega)\bigr|
<\frac\delta8\sum_{\ell=0}^{k-1}2^{\ell+1}\le \delta 2^{k-2}
\end{equation}
Choosing $k_1=k_1(\omega)\ge k_0$ so that $|\chi(x_{k_0},\omega)|<\delta 2^{k_1-2}$, this and \eqref{2.13} imply
\begin{equation}
x\in A_{2^k}\quad\Rightarrow\quad \bigl|\chi(x,\omega)\bigr|<\delta 2^k,\qquad k\ge k_1.
\end{equation}
But this means that for~$n\in\{2^k\colon k\ge k_1\}$,
\begin{equation}
\sum_{|x|\le n}\1_{\{|\chi(x,\omega)|\ge\delta n\}}\le|\Lambda_n\setminus A_n|\le \delta 2^d|\Lambda_n|.
\end{equation}
As~$\delta$ was arbitrary, this proves \eqref{limit-on-average} for~$n$ increasing along powers of two. A moment's thought now reveals that the same then holds for the unrestricted limit as well.
\end{proofsect}

As for Theorem~\ref{thm-BP}, we refer the reader to  Biskup and Prescott~\cite{Biskup-Prescott}. It should be emphasized that, although the assumptions to all the above are those of the annealed invariance principle, we in addition require the validity of the diffusive bounds \twoeqref{diffusive}{on-diag}. These are by no means guaranteed for a general ergodic~$\BbbP$, so the problem whether the annealed and quenched invariance principle hold simultaneously remains open.

We close this subsection with a simple exercise concerning the invariance principle for the variable-speed continuous-time version of our random walk.

\begin{hwproblem}
\label{ex-VSRW}
Suppose the ``usual assumptions'' and assume that $(X_n)$ obeys the Quenched Invariance Principle with the limiting Brownian motion having covariance \eqref{E:Cov}. Show that the variable-speed continuous time walk~$X_t$ obeys a Quenched Invariance Principle with the limiting Brownian motion having covariance
\begin{equation}
E\bigl[(\lambda\cdot B_t)^2\bigr]=\E\Bigl(\,\sum_x\omega_{0,x}\bigl[\lambda\cdot\Psi(\omega,x)\bigr]^2\Bigr).
\end{equation}
\end{hwproblem}

Note that the quantity on the right-hand side is closely related to the infimum of \eqref{E:Dir}, which was used to define the corrector. The appearance of expectation~$\E$ instead of~$E_\Q$ is due to the fact that~$\BbbP$ is invariant for the point of view of the particle induced by the VSRW. As to the constant-speed walk, here the quenched invariance principle follows from the discrete-time case by a strong asymptotic concentration of a sum of i.i.d.~exponential times.

\subsection{Known results and open problems}
\label{sec4.4}\noindent
The following sums up the principal steps in the progress towards proving quenched invariance principle in the class of Random Conductance Models:
\settowidth{\leftmargini}{(111)}
\begin{itemize}
\item
\emph{Strongly elliptic, ergodic~$\BbbP$}: proved by Sidoravicus and Sznitman~\cite{Sidoravicius-Sznitman}.
\item
\emph{Nearest-neighbor, i.i.d.\ conductances in~$d\ge2$ subject to the conditions}:
\begin{equation}
\E(\omega_{0,\hate})<\infty\quad\text{and}\quad\BbbP(\omega_e>0)>p_\cc(d)
\end{equation}
where $p_\cc(d)$ is the bond-percolation threshold. Here the quenched CLT has been proved in a sequence of papers by Sidoravicus and Sznitman~\cite{Sidoravicius-Sznitman}, Berger and Biskup~\cite{Berger-Biskup}, Mathieu and Piatnitski~\cite{Mathieu-Piatnitski}, Mathieu~\cite{Mathieu-CLT}, Biskup and Prescott~\cite{Biskup-Prescott}, Barlow and Deuschel~\cite{Barlow-Deuschel} with all approaches synthesized together by Andres, Barlow, Deuschel and Hambly~\cite{ABDH}.
\item
\emph{Nearest-neighbor, ergodic~$\BbbP$ with in $d=1,2$ with}
\begin{equation}
\E(\omega_{0,\hate})<\infty\quad\text{and}\quad\E(1/\omega_{0,\hate})<\infty
\end{equation}
The $d=2$ case is proved in these notes; the $d=1$ is the content of Exercises~\ref{HW:d=1} and~\ref{ex4.2} and goes back at least to Kawazu and Kesten~\cite{Kawazu-Kesten}.
\end{itemize}
We remark that that the condition $\E(\omega_{0,e})<\infty$ is essentially necessary; indeed Barlow and \v Cern\'y~\cite{Barlow-Cerny} ($d\ge3$) and \v Cern\'y~\cite{Cerny} ($d=2$) proved that for i.i.d.\ nearest neighbor conductances with $\alpha$-stable upper tail, $\alpha<1$, the law of $X_{nt}$ is under proper scaling described by $B_{W_t}$, where $B_t$ is a Brownian motion and $W_t$ is the inverse of an independent stable subordinator with index~$\alpha$. In other words, the paths are still Brownian but the heavy edges introduce a non-trivial trapping effects thus rendering the time parametrization non-linear and, in fact, stochastic. We remark that in physics, the limiting process is referred by as the \emph{fractional kinetics process}.

An important open problem concerns the rate of convergence and quantification of errors in martingale approximations. Although optimal results are probably far from reach, interesting ideas have been developed and quantitative results derived by Mourrat~\cite{Mourrat} and Gloria and Mourrat~\cite{Gloria-Mourrat}. The aforementioned work of Gloria and Otto~\cite{Gloria-Otto} gives integrability estimates on the corrector in~$d\ge3$ under strong ellipticity.

The Random Conductance Model has been also studied over other base-graphs than just~$\Z^d$. For instance, Caputo, Faggionato and Prescott~\cite{Caputo-Prescott} have investigated the random walks over various point processes in~$\R^d$. Independent studies for random walks on Voronoi/Delaunay triangulations have been announced  by Buckley~\cite{Buckley}. Ferrari, Grisi and Groisman~\cite{Ferrari} have constructed the harmonic coordinates on such triangulations by means of an interacting particle system; namely, a \emph{harness process}, which is basically a full-space stochastic version of the algorithm described for the finite boxes in Sect.~\ref{sec3.1}. The methods of Kipnis and Varadhan can be applied even to some deterministic quasiperiodic structures; see, e.g., Telcs~\cite{Telcs} who recently established an annealed invariance principle for the simple random walk on Penrose tilings.

\smallskip
Although we are able to control the corrector to the level required for the quenched invariance principle, the object itself remains rather mysterious and many open questions remain. For instance, regardless of what has been said at the end of Section~\ref{sec3.2}, the following problem remains of great interest both from the perspective of probability and analysis:

\begin{problem}
\label{pb4.10}
Is it true that a.e.~realization of random conductances satisfying the ``usual conditions'' admits no non-constant, sublinear harmonic functions?
\end{problem}

Recently, Benjamini, Duminil-Copin, Kozma and Yadin~\cite{Banff-theorem} have shown that that on the supercritical percolation cluster in~$\Z^d$, the space of linear harmonic functions is exactly $d+1$-dimensional. In particular, a typical supercritical percolation cluster supports no non-constant sublinear harmonic functions. We expect this to hold for all i.i.d.\ nearest-neighbor Random Conductance Models; for general environments the problem remains open.

Another open question concerns the scaling limit of the corrector:

\begin{problem}
Assume i.i.d.\ nearest-neighbor, uniformly elliptic Random Conductance Model. Show that the corrector scales to a Gaussian Free Field. More precisely, for any function $f\colon\R^d\to\R$ which is smooth and of compact support and satisfies $\int f(x)\textd x=0$, the law of
\begin{equation}
\chi_\epsilon(f):=\epsilon^{1+d/2}\int\chi\bigl(\omega,\lfloor x\rfloor\bigr)f(\epsilon x)\textd x
\end{equation}
scales, as $\epsilon\downarrow0$, to a Gaussian with mean zero and variance proportional to $(f,-\Delta^{-1}f)_{L^2(\R^d)}$.
\end{problem}

Progress in the uniformly elliptic case has been achieved in recent work of Gloria and Otto~\cite{Gloria-Otto} who have been able to prove that the corrector is in $L^q(\BbbP)$ for all~$q<\infty$, and thus a tight random variable, in all dimensions $d\ge3$. This settled an open problem from~\cite{Berger-Biskup}.

Another, perhaps somewhat related, question is that of the very definition of the corrector. Indeed, the corrector is defined almost surely for every ergodic law on environments~$\BbbP$. However, as different ergodic laws are  singular with respect to one another, it is not clear how to mesh the various correctors together. And yet it seems this should be possible:

\begin{problem}[Universal corrector]
Consider the set of nearest-neighbor environments $\Omega:=[a,b]^{\B(\Z^d)}$ where $0<a<b<\infty$. Define a function $\chi\colon\Omega\times\Z^d\to\R^d$ such that, for every ergodic law~$\BbbP$ on~$\Omega$, it agrees with the corrector corresponding to measure~$\BbbP$.
\end{problem}

We remark that this would be solved if one could find a sequence of local functions $\varphi_n$  such that $\nabla\varphi_n\to\chi$ almost surely for every~$\BbbP$. Note that, although may find functions~$\varphi_n$ for which the convergence takes place in~$\Lvec^2$ for any given~$\BbbP$, almost sure convergence requires reduction to subsequences which may be strongly $\BbbP$-dependent.

\smallskip
The understanding of the Markov chain permits one to consider more complicated questions. One such question concerns the typical number of points visited by the random walk in a given time. This was recently addressed by~Rau\cite{Rau}. Another question is the Law of Iterated Logarithm; this was established by Duminil-Copin~\cite{Duminil-Copin}. Next is the question of the behavior of the random walk on very thin percolation clusters. This can be studied directly in the case when $p=p_\cc$ where, technically speaking, the percolation cluster does not exist but one can still enforce it by conditioning. For the resulting \emph{incipient infinite cluster} (IIC) in sufficiently high dimensions, Nachmias and Kozma~\cite{Nachmias-Kozma} proved the Alexander-Orbach conjecture in all dimensions $d\ge7$ --- modulo caveats regarding the existing level of lace-expansion technology. This conjecture, due to Alexander and Orbach~\cite{Alexander-Orbach}, states that, on IIC,
\begin{equation}
\label{AO-conj}
\cmss P_\omega^{2n}(0,0)\asymp n^{-4/3},\qquad n\to\infty.
\end{equation}
Notably, this is expected to be false in low spatial dimensions. Related to this would be the decay of the diffusive constant for the simple random walk on the supercritical cluster for parameter~$p$, as $p\downarrow p_\cc$. Here we pose:

\begin{problem}
Suppose $d\ge7$ and let $D(p)$ denote the limiting variance of the simple random walk on the supercritical percolation cluster on~$\Z^d$ for parameter~$p>p_\cc(d)$. Show that
\begin{equation}
\label{D-conj}
D(p)\asymp (p-p_\cc)^2,\qquad p\downarrow p_\cc.
\end{equation}
\end{problem}

This problem is closely related to the existence of effective conductivity which was studied in, e.g., Grimmett and Kesten~\cite{Grimmett-Kesten}, Chayes and Chayes~\cite{Chayes-Chayes} and Kesten's monograph~\cite{Kesten} on percolation. See also Sect.~\ref{sec6}. A rather convincing argument can be obtained for this by analyzing the formula \eqref{E:Cov} and making plausible assumptions on the structural properties of the percolation cluster. Resorting to the electrostatic interpretation, the electric current should be carried only by the backbone of the cluster --- which, in the limit $p\downarrow p_\cc$, becomes a ``net'' of fractal curves. The exponent in \eqref{D-conj} then comes from realizing that in $d\ge7$, these fractals have Hausdorff dimension~$2$ (although the relation is not so straightforward as a simple equality of these numbers). This intuition seems be confirmed by observations made in the physics literature, see, e.g., Schr\o der and Dyre~\cite{Schroder-Dyre}. A main puzzle that remains is whether, and how exactly should the exponent $2$ in \eqref{D-conj} be related to the exponent~$4/3$ in \eqref{AO-conj}.

We remark that the amount of physics literature written on this and related subject is absolutely overwhelming; just see the articles citing the review by Dyre and Schr\o der~\cite{Dyre-Schroder-2000}.

\smallskip
Another very interesting class of applications of the above techniques is the random walk in \emph{dynamical} (albeit still reversible) random environments. We will not go into details here, but let us just say that much of Kipnis-Varadhan theory carries to this case and so annealed limit theorems are available. However, the understanding of \emph{quenched} invariance principles is far less evolved. Much can be said when the dynamics of the environment is Markovian and there is enough mixing; one can then get enough control via regeneration arguments. However, even here it is far from clear how to formulate convenient, and very general, conditions under which  invariance principles can be obtained.

From the perspective of this text, one specific class of dynamical random environments is of special interest. Consider a function $V\colon\R\to\R$ which is twice continuously differentiable and define a collection of coupled diffusions $(\phi_x(t))_{x\in\Z^d}$ via
\begin{equation}
\textd\phi_x(t)=\sum_{y\colon|y-x|=1}V'\bigl(\phi_y(t)-\phi_x(t)\bigr)\textd t+\sqrt2\,\textd B_t(x),
\end{equation}
where $B_t(x)$ are independent standard Brownian motions. As it turns out, any gradient Gibbs measure for the potential~$V$ is stationary under this dynamics. Assuming that~$V$ is convex, and thus~$V''\ge0$, we can now define a random walk~$X=(X_t)$ which at time~$t$ at position~$X_t=x$ takes a jump to a neighbor~$y$ at rate $V''(\phi_y(t)-\phi_x(t))$.

An attractive feature of this setting is that it permits us to analyze gradient models with convex interactions. For instance, we have the following formula
\begin{equation}
\label{E:HS}
\text{Cov}_\mu\bigl(\phi_0,\phi_x)=E_\mu E^{0,\phi}\Bigl(\,\int_0^\infty \1_{\{X_t=x\}}\textd t\Bigr)
\end{equation}
for the covariance of the (static) field in two locations with respect to a gradient Gibbs measure~$\mu$ by means of the expected number of visit to~$x$ by the above random walk started at~$0$ --- we expect this to be finite only in $d\ge3$ but other formulas exists in~$d=1,2$. Obviously, this generalizes the well-known formula from the Gaussian case which is distinguished by the fact that the random walk is not coupled to the evolution of the fields.

The formula \eqref{E:HS} is one instance of the \emph{Helffer-Sjostrand random walk representation} of correlation functions for the gradient model. These have been indispensable in the study of gradient models with convex interactions (e.g.,~Naddaf and Spencer~\cite{Naddaf-Spencer}, Giacomin, Olla and Spohn~\cite{Giacomin-Olla-Spohn}, Funaki~\cite{Funaki}, etc).


\section{Heat-kernel decay and failures thereof}
\label{sec5}\noindent
As discussed at length in the previous section, our current strategy of the proof of the quenched invariance principle seems to generally require the use of rather precise estimates on the probability that the Markov chain moves from~$x$ to~$y$ in~$n$ steps. We emphasize that this is conceptually flawed because we seem to need a \emph{local-CLT} type of result to finish a plain~CLT. Notwithstanding, the study of the heat kernel is interesting in its own right. We will only review the techniques that are ultimately relevant for the applications at hand and refer to, e.g., the upcoming textbook by Kumagai~\cite{Kumagai} for a more in-depth treatment of that well-developed area.

\subsection{Some general observations}
\noindent
To set the vocabulary straight, let us first remark that by the \emph{heat kernel} one usually means the quantity
\begin{equation}
\cmss q_n(x,y):=\frac{\cmss P_\omega^n(x,y)}{\pi_\omega(y)}.
\end{equation}
As one can expect, $\cmss P_\omega(x,\cdot)$ will for large~$n$ approach (a multiple of) the stationary measure~$\pi_\omega$. So $\cmss q_n$, being in fact the Radon-Nikodym derivative of $\cmss P_\omega^n(x,\cdot)$ with respect to $\pi_\omega$, is a very natural object to consider. Note that reversibility implies $\cmss q_n(x,y)=\cmss q_n(y,x)$.

Theorem~\ref{thm-BP} required in \eqref{on-diag} that the return probability generally decays as~$n^{-d/2}$. It turns out that, should the CLT hold, we cannot hope for a faster decay than this:

\begin{lemma}
\label{lemma-reg-LB}
Suppose $(X_n)$ is satisfies a CLT with non-degenerate diffusion constant~$\sigma^2$. Assume that $\pi^\star:=\sup_x\pi_\omega(x)<\infty$. Then there is $c=c(d,\sigma^2,\pi^\star)>0$ such that for~$n$ sufficiently large,
\begin{equation}
\cmss P_\omega^{2n}(0,0)\ge\frac{c}{n^{d/2}}\pi_\omega(0).
\end{equation}
\end{lemma}

\begin{proofsect}{Proof}
We use reversibility and simple estimates to get
\begin{equation}
\begin{aligned}
\cmss P_\omega^{2n}(0,0)
&=\sum_x\cmss P_\omega^n(0,x)\cmss P_\omega^n(x,0)
\\
&=\sum_x\cmss P_\omega^n(0,x)^2\frac{\pi_\omega(0)}{\pi_\omega(x)}
\\
&\ge \frac{\pi_\omega(0)}{\pi^\star}\sum_{|x|\le \sqrt n}\cmss P_\omega^n(0,x)^2
\end{aligned}
\end{equation}
The sum on the right-hand can be further bounded using the Cauchy-Schwarz inequality:
\begin{equation}
\cmss P_\omega^{2n}(0,0)\ge\frac{\pi_\omega(0)}{\pi^\star}\frac{P_\omega^0(|X_n|\le\sqrt n)^2}{|\{x\colon|x|\le \sqrt n\}|}
\end{equation}
But the CLT ensures that $P_\omega^0(|X_n|\le\sqrt n)\ge\frac12 P(|B_t|\le 1/\sigma)$ for $n$ large, where $B_t$ is the standard $d$-dimensional Brownian motion, and $|\{x\colon|x|\le \sqrt n\}|\le c'n^{d/2}$ for some~$c'=c'(d)<\infty$.
\end{proofsect}

We remark that a general method of getting such (including ``near-diagonal'') lower bounds in elliptic random environments has been put forward by Nash~\cite{Nash} and Fabes and Stroock~\cite{Fabes-Stroock}.

For reasons discussed earlier, the main technical problem is to find natural conditions on the Markov chain so that an $n^{-d/2}$ \emph{upper} bound can be guaranteed. This problem has been studied for over half a century, starting from proofs of regularity of elliptic PDEs with irregular coefficients (De Giorgi~\cite{DeGiorgi}, Nash~\cite{Nash}, Aronson~\cite{Aronson}) and validity and consequences of Faber-Krahn, Sobolev and Nash inequalities for diffusions on manifolds and Markov chains (e.g., Varopoulos~\cite{Varopoulos}, Carlen, Kusoka and Stroock~\cite{Carlen-Kusoka-Stroock}). A method to get \emph{off-diagonal} bounds --- i.e., for $\cmss q_n(x,y)$ with $x\ne y$ --- has been put forward by Davies~\cite{Davies-off-diag} based on the \emph{Carne-Varopoulos bound} (Carne~\cite{Carne}, Varopoulos~\cite{Varopoulos-CV}).

In the course of time it has been realized that there is a close connection between the desired upper bound and the \emph{geometric} properties of the underlying state-space. The key property to check is the validity of the \emph{isoperimetric inequality} (Cheeger~\cite{Cheeger}) or, more generally, the character of the \emph{isopertimetric profile} (Grigoryan~\cite{Grigoryan}). This connection was later transferred to the context of (discrete-space) Markov chains by Lawer and Sokal~\cite{Lawler-Sokal} and Jerrum and Sinclair~\cite{Jerrum-Sinclair} (invoking isoperimetric inquality) and, later, by Lov\'asz and Kannan~\cite{Lovasz-Kannan} and Morris and Peres~\cite{Morris-Peres} (based on isoperimetric profile).

We will not try to delve deeper into the details of historical developments of the subject; instead, the reader should consult the many texts that have been written on this (e.g., by Coulhon and Grigor'yan~\cite{Coulhon-Grigoryan}, Davies~\cite{Davies}, Kumagai~\cite{Kumagai}, Montenegro and Tetali~\cite{Montenegro-Tetali}, Varopoulos~\cite{Varopoulos-book}, Varopoulos, Saloff-Coste and Coulhon~\cite{Varopoulos-book}, Woess~\cite{Woess}, etc). For us the key fact is that with many Markov chains we may associate a natural graph structure --- simply put an edge between any two states in the state spaces that have a positive transition probability of a jump from one to the other. This permits us to connect the mixing properties of the chain with facts about geometry of this graph.

To illustrate this on an example, consider a graph that consists of two bulky components connected only by a few edges. Clearly, it will take quite a long time to exit one component and discover the other. Naturally, one is thus lead to comparing the size of a set with the size of its boundary which is expressed very well in terms of aforementioned {isoperimetric inequalities}.

\smallskip
In what follows we will rely on a result from a recent work by Morris and Peres~\cite{Morris-Peres} which we find particularly attractive for its probabilistic flavor. Consider a countable state Markov chain with state space~$V$, transition kernel $\cmss P$ and a stationary reversible measure~$\pi$. For a finite set~$A\subset V$, we will measure the boundary via
\begin{equation}
\cmss Q(A,A^\cc):=\sum_{\begin{subarray}{c}
x\in A\\y\in A^\cc
\end{subarray}}
\pi(x)\cmss P(x,y)
\end{equation}
and the volume via
\begin{equation}
\pi(A):=\sum_{x\in A}\pi(x).
\end{equation}
Define the function
\begin{equation}
\label{E:5.7eq}
\phi(r):=\inf\biggl\{\frac{\cmss Q(A,A^\cc)}{\pi(A)}\colon \pi(A)\le r\biggr\}
\end{equation}
that expresses the size of the least possible surface-to-volume ratio for all sets with volume less than~$r$. We can call this function the \emph{isoperimetric profile}. Its computation is often facilitated by the following fact:

\begin{hwproblem}
Show that in \eqref{E:5.7eq} we can restrict to~$A$ that are connected --- in the sense that for every~$x,y\in A$ there is a time~$n$ with $\cmss P^n(x,y)>0$.
\end{hwproblem}

We now quote verbatim Theorem~2 of~\cite{Morris-Peres}:

\begin{theorem}
\label{thm-MP}
Suppose that $\cmss P(x,x)\ge\gamma$ for some $\gamma\in(0,\ffrac12)$. For all~$\epsilon>0$, all $x,y\in V$ and all~$n$ satisfying
\begin{equation}
\label{E:n>integral}
n\ge1+\Bigl(\frac{1-\gamma}{\gamma}\Bigr)^2\int_{4[\pi(x)\wedge\pi(y)]}^{4/\epsilon}\frac{\textd r}{r\phi(r)^2}
\end{equation}
we have
\begin{equation}
\cmss P^n(x,y)\le\epsilon\pi(y).
\end{equation}
\end{theorem}

The restriction to uniformly positive holding probability, $\cmss P(x,x)\ge\gamma$, is a technical nuisance in applications that often requires analyzing a modified chain that has this property.

\subsection{Heat kernel on supercritical percolation cluster}
\noindent
It is quite instructive to check how Theorem~\ref{thm-MP} implies the usual bound for the simple random walk and/or elliptic nearest-neighbor environments. However, we will instead do something far less trivial; namely, we will show how this theorem applies in the case of the random walk on the supercritical percolation cluster.

\begin{theorem}
\label{thm-HK-perc}
Suppose $d\ge2$ and $p>p_\cc(d)$. There is a constant $c=c(d,p)<\infty$ and a random variable $n_0=n_0(\omega)$ such that for almost every sample of the bond-percolation cluster $\CC_\infty$ containing the origin, we have
\begin{equation}
\label{E:HK-perc}
\cmss P_\omega^{2n}(0,0)\le\frac c{n^{d/2}},\qquad n\ge n_0.
\end{equation}
\end{theorem}

For a finite set~$A\subset\CC_\infty(\omega)$, let $\partial^\omega A$ denote the set of open edges in~$\omega$ with exactly one endpoint in~$A$. A simple observation yields
\begin{equation}
\frac{\cmss Q(A,A^\cc)}{\pi(A)}\ge\frac1{2d}\frac{|\partial^\omega A|}{|A|}
\end{equation}
If $\omega:=1$ for all edges, then $\partial^\omega A=\partial A$. In such circumstances, one has the \emph{isoperimetric inequality} of the form: There is a constant $c=c(d)>0$, such that
\begin{equation}
\label{E:Euclid-iso}
|\partial A|\ge c|A|^{\frac{d-1}d},\qquad A\subset\Z^d\quad\text{finite}.
\end{equation}
This inequality cannot hold on $\CC_\infty$ because the infinite component contains arbitrarily long one dimensional (and other) pieces. However, we can have this for connected sets that are not too small compared to their distance to the~origin:

\begin{lemma}
\label{thm-isoperimetry}
For all $d\ge2$ and $p>p_\cc(d)$, there are positive and finite constants $c_1=c_1(d,p)$ and $c_2=c_2(d,p)$ and an a.s.\ finite random variable $R_0=R_0(\omega)$ such that for each $R\ge R_0$ and each $\omega$-connected~$A$ satisfying
\begin{equation}
A\subset\CC_\infty\cap[-R,R]^d\quad\text{and}\quad |A|\ge (c_1\log R)^{\frac d{d-1}}
\end{equation}
we have
\begin{equation}
\label{isoperimetry}
|\partial^\omega A|\ge c_2|A|^{\frac{d-1}d}.
\end{equation}
\end{lemma}

There have been a number of proofs of this and/or related results, see e.g. Benjamini and Mossel~\cite{Benjamini-Mossel}, Heicklen and Hoffman~\cite{Heicklen-Hoffman}, Mathieu and Remy~\cite{Mathieu-Remy}, Barlow~\cite{Barlow}, Berger-Biskup-Hoffman-Kozma~\cite{BBHK}, Pete~\cite{Pete}. We will not prove this claim here for all $p>p_\cc(d)$ as the proof uses non-trivial facts from percolation theory. However, for $p$ very close to 1 there is a much simpler argument due to Benjamini and Mossel:

\begin{hwproblem}
Show that once $p$ is sufficiently close to one, there is a constant~$c_1\in(0,\infty)$ and a random variable~$R_0=R_0(\omega)<\infty$ such that for all $R\ge R_0$,
\begin{equation}
A\subset[-R,R]^d\cap\Z^d\,\,\text{ and }\,\,|A|\ge(c_1\log R)^{\frac d{d-1}}\,\,\text{ imply }\,\,|\partial^\omega A|\ge\frac12|\partial A|.
\end{equation}
\end{hwproblem}

Note that from here we will immediately have \eqref{isoperimetry} via \eqref{E:Euclid-iso}. In order to see how \eqref{isoperimetry} feeds into Theorem~\ref{thm-MP}, note that the Markov chain by time~$2n$ will not leave the box $[-2n,2n]$. Thus set $R:=2n+1$, pick $\theta\in(0,\ffrac12)$ and for $A\in[-R,R]^d\cap\Z^d$ connected let us estimate the ratio in the definition of $\phi(r)$ by $c|A|^{-1/d}$ when $|A|\ge R^\theta$ and by $cR^{-\theta}$ when $|A|\le R^\theta$. (In the second step we used that $|\partial^\omega A|\ge1$.) It follows that
\begin{equation}
\phi(r)\ge c\bigl(r^{-1/d}\wedge R^{-\theta}\bigr)
\end{equation}
Plugging this into \eqref{E:n>integral}, the integral is at most $cR^{2\theta}\log R+c\epsilon^{-2/d}$. This will be less than~$n$ for $\epsilon:=cn^{-d/2}$. The inequality \eqref{E:HK-perc} then follows by applying Theorem~\ref{thm-MP}.

A natural consequence of Theorem~\ref{thm-HK-perc} is the result that was first proved by Grimmett, Kesten and Zhang~\cite{Grimmett-Kesten-Zhang} by rather different methods (see also Problem~\ref{pb1.15}):

\begin{corollary}
The simple random walk on (a.e.~realization of) the supercritical percolation cluster is recurrent in dimension $d=2$ and transient in dimensions $d\ge3$.
\end{corollary}

\begin{proofsect}{Proof}
As explained in Sect.~\ref{sec1.2}, it suffices to resolve the $d=3$ case, but we can cover all~$d$'s just as well. From Lemma~\ref{lemma-reg-LB} and Theorem~\ref{thm-HK-perc} we know that $\cmss P_\omega^0(0,0)\asymp n^{-d/2}$. This is summable in dimensions $d\ge3$ and non-summable in $d=1,2$. The summability is then equivalent to the finiteness of the full-lattice Green's function which via~\eqref{E:1.38} is then equivalent to transience.
\end{proofsect}

\subsection{Anomalous decay}
\noindent
From the perspective of nearest-neighbor random walks on~$\Z^d$, the case of the supercritical percolation is a prototype of a non-elliptic situation. However, when we think of this walk as the simple random walk on the \emph{graph}~$\CC_\infty$, it is as elliptic as the one can ever hope for. Indeed, any edge in the graph~$\CC_\infty$ has conductance one and the \emph{ellipticity contrast} --- the difference between a maximal and minimal possible value of the conductance over each edge --- is zero.  The difficulties in the understanding of this walk on~$\CC_\infty$ is thus not the lack of ellipticity but the intricacies of its random geometry.

From this point of view it is natural to ask what happens when ellipticity gets violated in a robust way. This naturally leads to consideration of i.i.d.\ nearest-neighbor environments where the law of the individual conductances is  unbounded either from zero or from infinity (or both). The point is that both situations can lead to trapping effects although each of them for a slightly different reason. We will henceforth focus on the former case and refer to Barlow and Deuschel~\cite{Barlow-Deuschel} for the latter.

Suppose, from now on, that the $\omega$'s are nearest-neighbor, i.i.d.\ with $\BbbP(0<\omega_b\le1)=1$ and
\begin{equation}
\text{essinf}(\omega_b)=0.
\end{equation}
Our assumption implies that~$\BbbP(\omega_b\in\cdot)$ has no atom at zero. Thus all nearest-neighbor jumps on~$\Z^d$ are allowed for the random walk, but some of them may be very unlikely.

It is easy to check that for i.i.d.\ distribution with these properties, the isoperimetry methods sketched above yield a vacuous conclusion. The situation becomes even more suspicious after an inspection of the work of Fontes and Mathieu~\cite{Fontes-Mathieu} in which they design a family of models --- not with i.i.d.\ conductances but close enough --- in which the \emph{expected} diagonal heat kernels, $\E\cmss P_\omega^{2n}(0,0)$, decay arbitrarily slowly with~$n$. Of course, this could be just a result of taking an average over the environment (remember that we are talking about events whose probabilities decay to zero) so one is naturally intrigued by what the typical (quenched) decay of~$\cmss P_\omega^{2n}(0,0)$ might be.

It will not be too surprising that in $d=1$ the trapping can be quite severe even for typical~$\omega$. Indeed, the following is an interesting exercise:

\begin{hwproblem}
Suppose $d=1$ and nearest-neighbor, i.i.d.\ conductances with values in $(0,1]$. For each sequence $\lambda_n\to\infty$ construct a law~$\BbbP$ such that
\begin{equation}
\cmss P_\omega^{2n}(0,0)\ge\frac1{\lambda_n}
\end{equation}
for~$n$ large, along a deterministic subsequence $n_k\to\infty$.
\end{hwproblem}

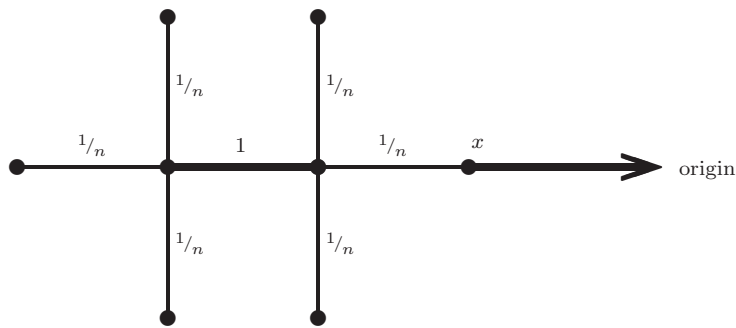
\begin{figure}[t]
\refstepcounter{obrazek}
\vglue0.2cm
\begin{center}
\hglue-1cm
\setlength{\unitlength}{20mm}
\begin{picture}(2.5,1)
\linethickness{1.2pt}
\put(0,0){\line(0,1){1.0}}
\put(0,0){\line(0,-1){1.0}}
\put(0,0){\line(-1,0){1.0}}
\put(1,0){\line(1,0){1.0}}
\put(1,0){\line(0,1){1.0}}
\put(1,0){\line(0,-1){1.0}}
\put(0,0){\circle*{0.1}}
\put(0,1){\circle*{0.1}}
\put(0,-1){\circle*{0.1}}
\put(-1,0){\circle*{0.1}}
\put(1,0){\circle*{0.1}}
\put(1,0){\circle*{0.1}}
\put(2,0){\circle*{0.1}}
\put(1,1){\circle*{0.1}}
\put(1,-1){\circle*{0.1}}
\linethickness{2.6pt}
\put(0,0){\line(1,0){1.0}}
\linethickness{2.6pt}
\put(2,0){\line(1,0){1.2}}
\thicklines
\put(3.22,0){\line(-3,-1){0.2}}
\put(3.24,0){\line(-3,-1){0.22}}
\put(3.26,0){\line(-3,-1){0.24}}
\put(3.28,0){\line(-3,-1){0.26}}
\put(3.22,0){\line(-3,1){0.2}}
\put(3.24,0){\line(-3,1){0.22}}
\put(3.26,0){\line(-3,1){0.24}}
\put(3.28,0){\line(-3,1){0.26}}
\put(0.45,0.1){$1$}
\put(1.4,0.1){$\ffrac1n$}
\put(-0.6,0.1){$\ffrac1n$}
\put(0.05,0.5){$\ffrac1n$}
\put(0.05,-0.55){$\ffrac1n$}
\put(1.05,0.5){$\ffrac1n$}
\put(1.05,-0.55){$\ffrac1n$}
\put(3.4,-0.05){origin}
\put(2.02,0.13){$x$}
\end{picture}
\end{center}
\vskip2cm
\caption{A trap capable of capturing the walk for times of order~$n$. Here an edge of conductance 1 is
separated by edges of conductance $\ffrac1n$ from a path of edges of conductance~1 to the origin.
Once the walk enters the trap, incurring a cost $\ffrac1n$ of probability, it will stay bouncing back and
forth for time of order~$n$ with a uniformly positive probability. The exit at a particular time costs
again $\ffrac1n$.}
\label{fig6}
\end{figure}

A moment's thought --- and a right idea --- then shows that interesting new behavior may actually occur even in high-enough dimensions. Consider the following example from the paper of Berger, Biskup, Hoffman and Kozma~\cite{BBHK}: Fix a sequence $\lambda_n\to\infty$ and define a \emph{trap of order~$n$} to be the configuration in Fig.~\ref{fig6} and suppose that the distance of this trap to the origin is $\ell_n$. If such a trap occurs,  we can estimate $\cmss P_\omega^{2n}(0,0)$ as follows. The cost of getting to vertex~$x$ from the origin is exponential in the distance, i.e., $\texte^{O(\ell_n)}$. Entering the trap at the next step costs order $\ffrac1n$ of probability. The walk can then be made to stay there for the time $2n$ minus twice the distance to the trap; this costs only $\texte^{O(1)}$ of probability. Exiting the trap at the required time costs one more $\ffrac1n$ and the trip back to the origin $\texte^{O(\ell_n)}$. In total, we thus have
\begin{equation}
\cmss P_\omega^{2n}(0,0)\ge\texte^{O(\ell_n)}\,\frac1n\,\texte^{O(1)}\,\frac1n\,\texte^{O(\ell_n)}=\frac{\texte^{O(\ell_n)}}{n^2}
\end{equation}
Now one just beefs up the lower tail of~$\BbbP$ so that, along a deterministic subsequence $n_k\to\infty$, we have $\ell_n=o(\log\lambda_n)$. We have a~proof~of:

\begin{theorem}
Suppose~$d\ge5$. For each~$\lambda_n\to\infty$ there exists an i.i.d.\ conductance law~$\BbbP$ satisfying $\BbbP(0<\omega_b\le1)=1$, a deterministic sequence $n_k\to\infty$ and a $\BbbP$-a.s.\ positive random variable $C=C(\omega)>0$ such that for each $n\in\{n_k\}$,
\begin{equation}
\label{E:d>4LB}
\cmss P_\omega^{2n}(0,0)\ge\frac{C(\omega)}{n^2\lambda_n}.
\end{equation}
\end{theorem}

Notice that the above argument yields a similar bound in all dimensions $d\ge2$, but this bound has no significant value in dimensions $d=2,3,4$ as (by the CLT proved by Mathieu~\cite{Mathieu-CLT} and, independently, Biskup and Prescott~\cite{Biskup-Prescott}) $\cmss P_\omega^{2n}(0,0)$ decays at least as $n^{-d/2}$; cf Lemma~\ref{lemma-reg-LB}. But in $d\ge5$ this shows that the heat kernel may decay \emph{more slowly} than $n^{-d/2}$ and, in particular, there is no way that a diffusive heat kernel upper bounds would generally hold.

An interesting question is whether \eqref{E:d>4LB} is the worst one can do. The answer turns out to be, more or less, in the affirmative:

\begin{theorem}
For any nearest-neighbor, i.i.d.\ conductance law~$\BbbP$  with $\BbbP(0<\omega_{0,\hate}\le1)=1$ there is a random variable $C=C(\omega)<\infty$ such that
\begin{equation}
\cmss P_\omega^{2n}(0,0)\le C(\omega)\,\begin{cases}
n^{-d/2},\qquad&d=2,3,
\\
n^{-2}\log n,\qquad&d=4,
\\
n^{-2},\qquad&d\ge5.
\end{cases}
\end{equation}
In addition, we have
\begin{equation}
\label{E:5.22}
\lim_{n\to\infty}n^2\cmss P_\omega^{2n}(0,0)=0,\qquad\BbbP\text{\rm-a.s.\ in }d\ge5.
\end{equation}
and
\begin{equation}
\label{E:5.23}
\lim_{n\to\infty}\frac{n^2}{\log n}\cmss P_\omega^{2n}(0,0)=0,\qquad\BbbP\text{\rm-a.s.\ in }d=4.
\end{equation}
\end{theorem}
All except \eqref{E:5.23} in this result is due to Berger, Biskup, Hoffman and Kozma \cite{BBHK}; the property \eqref{E:5.23} was derived only recently in Biskup, Louidor, Rozinov and Vandenberg-Rodes~\cite{UCLA-team}. The latter group has also shown that, in many cases where the heat kernel decays subdiffusively, the trapping phenomenon described in the example above actually occurs: the path spends $n-o(n)$ of time in a very small spatial region.

Notice that \eqref{E:d>4LB} and \eqref{E:5.22} nicely complement each other: anything up to, but no worse than, $o(n^{-2})$ decay can occur in $d\ge5$. A question remains whether the $\log n$ factor in~$d=4$ is an artifact of the proof or a real phenomenon. This was solved recently by Biskup and Boukhadra~\cite{Biskup-Boukhadra} who constructed an environment, for each sequence $\lambda_n\to\infty$, such that
\begin{equation}
\cmss P_\omega^{2n}(0,0)\ge\frac{\log n}{n^2\lambda_n},
\end{equation}
eventually, along a deterministic subsequence $n_k\to\infty$. The construction is quite involved because in $d=4$ the trapping occurs more or less equally likely over a whole range of exponentially-growing spatial scales (hence the $\log n$ factor).

\subsection{Conclusions}
\noindent
The upshot of the above results and derivations is that with the random conductance models we are finding ourselves in a somewhat unusual situation when the path distribution satisfies a non-degenerate functional CLT and yet the heat kernel decays anomalously; i.e., we have a \emph{CLT without local CLT}. Although this may contradict intuition, there is nothing wrong about this: a CLT is a statement about the bulk of the distribution and a local-CLT is a statement about the tails. There is no particular reason why these should match one another.



\section{Applications}
\label{sec6}\noindent
In this section we will try to address some aspects of the applications that were introduced in the first section of these notes. Specifically, we will discuss homogenization of discrete parabolic (random) problems, scaling limit of associated Green's functions, convergence of random Gaussian gradient models to Gaussian Free Field and, finally, applications to electrostatics.

\subsection{Some homogenization theory}
\noindent
The phrase "homogenization theory" usually refers to a diverse set of methods and ideas that address one of the fundamental problems of material science: the computation of macroscopic material constants and characteristics (e.g., heat or electric conductivity, resistivity, etc) from the microscopic properties. One of the typical mathematical issues resolved by homogenization theory concerns differential equations: Although the microscopic quantities evolve according to an differential equation with rapidly varying coefficients, properly rescaled macroscopic versions thereof are governed by equations with smooth coefficients.

We will not go into the subject and history of homogenization theory in any further detail; these can be found in the literature, e.g., the monograph by Jikov, Kozlov and Oleinik~\cite{JKO}. Instead, we will attempt to demonstrate the conclusions on an example of heat conduction.

\smallskip
Suppose that some material of a rapidly varying internal microscopic internal structure --- described at the lattice level of spacing $\epsilon$ by a configuration of conductances $\omega$ --- is put in a macroscopic temperature profile at time~$0$. At the lattice level, the evolution of the temperature profile with time is described by the Cauchy problem
\begin{equation}
\label{E:6.1}
\left\{
\begin{alignedat}{3}
&\frac\partial{\partial t} u(t,x)=\cmss L_\omega u(t,x),&\quad\qquad &t\ge0,\,x\in\Z^d,\
\\*[1mm]
&u(0,x)=f(x),&\qquad &x\in\Z^d,
\end{alignedat}
\right.
\end{equation}
where $\cmss L_\omega$ is the operator \eqref{L-omega} (acting only on the $x$ coordinate) that represents the microscopic diffusive properties of the material and~$f$ is the initial temperature profile. Our first question concerns the existence and uniqueness of the solution. We note the classical fact:

\begin{lemma}
Suppose~$\omega$ is a sample from an ergodic measure~$\BbbP$ with $\E\pi_\omega(0)<\infty$. Let~$X_t$ denote the variable-speed continuous-time Markov chain on~$\Z^d$ with generator~$\cmss L_\omega$. Pick~$f\colon\Z^d\to\R$ bounded. Then
\begin{equation}
\label{prob-sol}
u(t,x):=E^x_\omega\bigl(f(X_t)\bigr)
\end{equation}
is the unique solution to \eqref{E:6.1} which is bounded in both~$t$ and~$x$.
\end{lemma}

\begin{proofsect}{Proof}
By Exercise~\ref{ex2.8} and the general theory expounded in, e.g., Liggett~\cite{Liggett}, the conditions on~$\omega$ guarantee that a stochastic solution to the backward Komogorov equations \twoeqref{E:backward-Kolmogorov}{E:initial-cond} exits and the semigroup for the VSRW is well defined. The fact that \eqref{prob-sol} is a solution is then a consequence of a direct calculation. Indeed, we have
\begin{equation}
u(t,x)=\sum_{z\in\Z^d}\cmss R_\omega^t(x,z)f(z).
\end{equation}
and the boundedness of~$f$ and finiteness of~$\pi_\omega$ permit us to exchange the sum over~$z$ with the time-derivative and $\cmss L_\omega$. Hence, $u$ satisfies \eqref{E:6.1}.

The remaining issue is thus a proof of uniqueness among bounded solutions. Let $\tilde u(t,x)$ be such a solution and, for $0\le s\le t$, consider the random variable
\begin{equation}
M_s:=\tilde u(t-s,X_s)
\end{equation}
and let $\scrF_s:=\sigma(X_r\colon0\le r\le t)$. Then $\{M_s,\scrF_s\}_{0\le s\le t}$ is a martingale. Indeed, by the Markov property, on the event $\{X_s=z\}$ we have
\begin{equation}
\begin{aligned}
E^x_\omega(&M_{s+\delta}-M_s|\scrF_s)
=E^z_\omega\tilde u(t-s-\delta,X_\delta) - \tilde u(t-s,z)
\\
&\quad=E^z_\omega\bigl[\tilde u(t-s-\delta,X_\delta)-\tilde u(t-s-\delta,z)\bigr]
\\
&\qquad\qquad\qquad\qquad\qquad+\bigl[\tilde u(t-s-\delta,z)-\tilde u(t-s,z)\bigr].
\end{aligned}
\end{equation}
This yields
\begin{equation}
\lim_{\delta\downarrow0}\frac1\delta\bigl[E^x_\omega(M_{s+\delta}-M_s|\scrF_s)\bigr]
=\cmss L_\omega \tilde u(t-s,z)-\frac\partial{\partial t}\tilde u(t-s,z) = 0.
\end{equation}
almost surely for every~$s$. Integrating over final intervals and applying the Bounded Convergence Theorem proves that $\{M_s,\scrF_s\}_{0\le s\le t}$ is a martingale. (At $s=t$ we apply continuity from the left.)

The Optional Stopping Theorem then yields
\begin{equation}
E^x_\omega M_0=E^x_\omega M_t
\end{equation}
which reads
\begin{equation}
\tilde u(t,x) = E^x_\omega \tilde u(t,0)=E^x_\omega f(X_t)=u(t,x).
\end{equation}
The uniqueness is proved as well.
\end{proofsect}

\begin{hwproblem}
Construct a configuration of nearest-neighbor conductances on~$\Z$ for which there is a non-zero solution to \eqref{E:6.1} with $u(0,\cdot):=0$. 
\end{hwproblem}

Our next goal is to describe the asymptotic of the solution for the situation when~$f$ is a macroscopic profile over a lattice of spacing~$\epsilon$. Fix a function $f\colon \R^d\to\R$ in $L^{1,\text{loc}}(\textd x)$ and let $u^{(\epsilon)}(t,x)$ denote the unique bounded solution to \eqref{E:6.1}  with initial data
\begin{equation}
u(0,x):=\int_{[0,1]^d}\textd z\, f(\epsilon x+\epsilon z),\qquad x\in\Z^d.
\end{equation}
Under diffusive scaling of space and time, we get the quantity
\begin{equation}
u_\epsilon(t,x):=u^{(\epsilon)}\bigl(t\epsilon^{-2},\lfloor x\epsilon^{-1}\rfloor\bigr),\qquad t\ge0,\,x\in\R^d.
\end{equation}

\begin{theorem}
\label{thm6.3}
Suppose $f\colon\R^d\to\R$ obeys $\Vert f\Vert_{L^2(\R^d)}^2+\Vert \nabla f\Vert_{L^2(\R^d)}^2<\infty$ and let~$\BbbP$ be a law on the conductances satisfying the ``usual conditions,'' \eqref{E:square} and \eqref{E:4.3}. Let~$Q$ denote the generator of the (annealed) limiting Brownian motion and let $\bar u$ be the solution to the Cauchy problem
\begin{equation}
\label{E:6.1-cont}
\left\{
\begin{alignedat}{3}
&\frac\partial{\partial t} \bar u(t,x)=Q \bar u(t,x),&\quad\qquad &t\ge0,\,x\in\R^d,\
\\*[1mm]
&\bar u(0,x)=f(x),&\qquad &x\in\R^d.
\end{alignedat}
\right.
\end{equation}
Then for each~$t\ge0$,
\begin{equation}
u_\epsilon(t,\cdot)\,\,\underset{\epsilon\downarrow0}\longrightarrow\,\,\bar u(t,\cdot)
\qquad\text{ in }L^2(\textd x)\otimes L^2(\BbbP).
\end{equation}
\end{theorem}

\begin{proofsect}{Proof}
Let~$B_t$ be the Brownian motion with generator~$Q$. Then
\begin{equation}
\bar u(t,x)=E^0\bigl(f(x+B_t)\bigr).
\end{equation}
Similarly, resolving the above scaling relations yields
\begin{equation}
u_\epsilon(t,x):=E^{\lfloor x\epsilon^{-1}\rfloor}_\omega\Bigl(\,\int_{[0,1]^d}\textd z\, f(\epsilon X_{t\epsilon^{-2}}+\epsilon z)\Bigr).
\end{equation}
By translation-invariance of~$\BbbP$ and the Cauchy-Schwarz inequality,
\begin{multline}
\qquad
\E\int_\R \bigl|
u_\epsilon(t,x)-\bar u(t,x)\bigr|^2
\textd x
\\
\le\int_{[0,1]^d}\textd z\,\int\E\Bigl(\bigl|E^0_\omega f(x_\epsilon(z)+\epsilon X_{t\epsilon^{-2}})-E^0f(x+B_t)\bigr|^2\Bigr)\textd x,
\end{multline}where $x_\epsilon(z):=\epsilon\lfloor x\epsilon^{-1}\rfloor+\epsilon z$.

Our first step is to replace~$x_\epsilon(z)$ by~$x$ in the argument of the first $f$ on the right-hand side. The difference tends to zero when $\epsilon\downarrow0$ because we have
\begin{equation}
\int\bigl|E^0_\omega f(x_\epsilon(z)+\epsilon X_{t\epsilon^{-2}})-E^0_\omega f(x+\epsilon X_{t\epsilon^{-2}})\bigr|^2\textd x\le4\epsilon^2\Vert\nabla f\Vert_{L^2(\R^d)}^2.
\end{equation}
To control the remaining difference, we note that, by the Annealed CLT (in analogy with Corollary~\ref{cor-AIP}) there exists a coupling $Q_\omega^0$ of the random walk $X_{t\epsilon^{-2}}$ and the Brownian motion $B_t$ such that, for any~$\delta>0$ and any~$t>0$,
\begin{equation}
a_\epsilon(\delta):=\E Q_\omega^0\bigl(\,|\epsilon X_{t\epsilon^{-2}}-B_t|>\delta\bigr)\,\,\underset{\epsilon\downarrow0}\longrightarrow\,\,0.
\end{equation}
Picking an arbitrary~$\delta>0$, the bound
\begin{multline}
\label{E:6.18a}
\qquad
\E\, E_{Q_\omega^0}\int\bigl|f(x+\epsilon X_{t\epsilon^{-2}})-f(x+B_t)\bigr|^2\textd x
\\
\le4a_\epsilon(\delta)\Vert f\Vert_{L^2(\R^d)}^2+\delta^2\Vert\nabla f\Vert_{L^2(\R^d)}^2
\qquad
\end{multline}
then shows that the expectation on the left tend to zero as~$\epsilon\downarrow0$ (followed by $\delta\downarrow0$). The proof is then finished by noting that
\begin{equation}
\E\int\Bigl(\bigl|E^0_\omega f(x+\epsilon X_{t\epsilon^{-2}})-E^0f(x+B_t)\bigr|^2\Bigr)\textd x
\le\text{LHS of \eqref{E:6.18a}},
\end{equation}
as implied by using Cauchy-Schwarz one last time.
\end{proofsect}

Theorem~\ref{thm6.3} exemplifies a statement in homogenization theory. Indeed, a solution to the parabolic problem with rapidly varying coefficients does behave, at a large scale, as a solution to a parabolic problem with constant coefficients. As is seen from Exercise~\ref{ex-VSRW}, the coefficients in the equation, namely, the entries in the symmetric, positive semi-definite matrix $(q_{ij})$ in
\begin{equation}
\label{Q-op}
Qf(x)=\sum_{i,j=1}^d q_{ij}\frac{\partial^2 f}{\partial x_i\partial x_i}
\end{equation}
are given by
\begin{equation}
q_{ij}:=\E\Bigl(\,\sum_x\omega_{0,x}\bigl(\hate_i\cdot\Psi(\omega,x)\bigr)\bigl(\hate_j\cdot\Psi(\omega,x)\bigr)\Bigr),
\end{equation}
where $\Psi(\omega,x)$ is the harmonic coordinate discussed at length in Section~\ref{sec3}. Notice that these are characterized by a variational problem
\begin{equation}
\sum_{i,j=1}^d\lambda_i\lambda_j q_{ij}=\inf_\varphi\,
\E\Bigl(\,\sum_x\omega_{0,x}\bigl(\lambda\cdot x+\nabla_x\varphi(\omega)\bigr)^2\Bigr)
\end{equation}
where~$\lambda=(\lambda_1,\dots,\lambda_d)\in\R^d$ and where $\varphi\colon\Omega\to\R$ runs over all local functions. This is the same variational problem that defines the corrector. This is the desired formula that at least in principle allows us to compute material coefficients from its microscopic properties.

\subsection{Green's functions and gradient fields}
\noindent
The arguments in the previous section can be cast in a more symmetric form provided we are willing to invoke some functional analysis. Given an operator~$\scrO$ on $\ell^2(\Z^d)$ with coefficients $\scrO(x,y):=\langle\delta_x,\scrO\delta_y\rangle$, we can interpret it as an operator on~$L^2(\R^d)$ by way of
\begin{equation}
\langle\, f,\scrO g\rangle_{L^2(\R^d)}:=\int_{\R^d\times\R^d}\scrO\bigl(\lfloor x\rfloor,\lfloor y\rfloor\bigr) f(x)g(y)\,\textd x\,\textd y.
\end{equation}
For any $f\in L^2(\R^d)$ define
\begin{equation}
\label{f-epsilon}
f_\epsilon(x):=\epsilon^{d/2+1}f(x\epsilon).
\end{equation}
In this notation, the statement of Theorem~\ref{thm6.3} implies:

\begin{corollary}
\label{cor6.5}
For any smooth functions $f,g\colon\R^d\to\R$ of compact support,
\begin{equation}
\epsilon^{-2}\bigl\langle\, g_\epsilon,\texte^{t\epsilon^{-2}\cmss L_\omega} f_\epsilon\bigr\rangle_{L^2(\R^d)}
\,\,\underset{\epsilon\downarrow0}\longrightarrow\,\,
\langle\, g,\texte^{tQ}f\rangle_{L^2(\R^d)},\qquad \text{\rm in }L^2(\BbbP).
\end{equation}
\end{corollary}

\begin{proofsect}{Proof}
Just note that, in the notation of Theorem~\ref{thm6.3}, $\epsilon^{-2}\langle\, g_\epsilon,\texte^{t\epsilon^{-2}\cmss L_\omega} f_\epsilon\rangle=\langle g,u_\epsilon(t,\cdot)\rangle$ while $\langle\, g,\texte^{tQ}f\rangle=\langle g,\bar u(t,\cdot)\rangle$. These tend to each other as $u_\epsilon(t,\cdot)\to\bar u(t,\cdot)$ in $L^2(\textd x)\otimes L^2(\BbbP)$.
\end{proofsect}

Corollary~\ref{cor6.5} supplies the core idea underlying the proof of our next result:

\begin{theorem}
\label{thm6.6}
Consider any ergodic law $\BbbP$ on nearest-neighbor elliptic conductances and pick any $f,g\colon\R^d\to\R$ that are smooth and of compact support. In $d=1,2$ assume in addition that the integral of~$f$ and~$g$ over~$\R^d$ equals zero. Then
\begin{equation}
\label{E:6.23}
\bigl\langle\, g_\epsilon,(-\cmss L_\omega)^{-1}f_\epsilon\bigr\rangle_{L^2(\R^d)}
\,\,\underset{\epsilon\downarrow0}\longrightarrow\,\,
\bigl\langle\, g,(-Q)^{-1}f\bigr\rangle_{L^2(\R^d)},\qquad \text{\rm in }L^2(\BbbP).
\end{equation}
\end{theorem}

\begin{proofsect}{Proof} (Sketch)
We only sketch the main ideas; details for this setting can be found in the work of Biskup and Spohn~\cite{Biskup-Spohn}. All inner products will be those in~${L^2(\R^d)}$ so we will not make this notationally explicit.

First let us note that both inner products are well defined. Indeed, $-\cmss L_\omega$ is self-adjoint and positive semi-definite with empty kernel (in $\ell^2(\Z^d)$). Moreover, it is invertible on all functions of finite support in~$\Z^d$ subject to --- in dimensions $d=1,2$ --- the condition of a vanishing total sum. Uniform ellipticity gives us the  following inequality between norms:
\begin{equation}
\bigl\langle\, f,(-\cmss L_\omega)^{-1}f\bigr\rangle
\le\Vert f\Vert^2_2+c\bigl\langle\, f,(-\Delta)^{-1}f\bigr\rangle,
\end{equation}
where $\Delta$ is a continuum Laplacian and where the passage from discrete to continuum Laplacian is due to~\cite[Lemma~2.2]{Biskup-Spohn}. As is not hard to check, replacing~$f$ by~$f_\epsilon$ on the left and using that $\Vert f_\epsilon\Vert_2=\epsilon^2\Vert f\Vert_2$ while $\langle\, f,(-\Delta)^{-1}f_\epsilon\rangle=\langle\, f,(-\Delta)^{-1}f\rangle$, the bound still holds all~$\epsilon\le1$. The family in \eqref{E:6.23} is thus uniformly bounded.

By the polarization identity, it suffices to prove the claim for $g:=f$. To this end we notice the following representation
\begin{equation}
\label{E:6.25}
\begin{aligned}
\bigl\langle\, f_\epsilon,(-\cmss L_\omega)^{-1}f_\epsilon\bigr\rangle
&=\int_0^\infty
\bigl\langle\, f_\epsilon,\texte^{t\cmss L_\omega} f_\epsilon\rangle\,\textd t,
\\
&=\int_0^\infty
\epsilon^{-2}\bigl\langle\, f_\epsilon,\texte^{t\epsilon^{-2}\cmss L_\omega} f_\epsilon\rangle\,\textd t,
\end{aligned}
\end{equation}
where we  scaled~$t$ by~$\epsilon^{-2}$ in the second line. By Corollary~\ref{cor6.5}, the integrand on the right-hand side tends to $\langle\, f,\texte^{tQ} f\rangle$ so, ignoring the important issue whether we are able to interchange the limit and the integral, we should have
\begin{equation}
\bigl\langle\, f_\epsilon,(-\cmss L_\omega)^{-1}f_\epsilon\bigr\rangle
\,\,\underset{\epsilon\downarrow0}\longrightarrow\,\,
\int_0^\infty\langle\, f,\texte^{tQ}g\rangle\,\textd t.
\end{equation}
The right hand side is again bounded by the fact that~$Q$ is uniformly elliptic, and it equals the term $\langle\, f,(-Q)^{-1}f\rangle$.

The key technical point of the proof is thus the control of the tails of the integral in \eqref{E:6.25}. This is a non-trivial problem where we will have to invoke, once again, heat-kernel estimates. This is easier in dimensions $d\ge3$ where it suffices to invoke the result of Delmotte~\cite{Delmotte}:
\begin{equation}
\cmss P^t_\omega(x,y)\le \frac{c}{t^{d/2}},\qquad x,y\in\Z^d,\, t\ge0,
\end{equation}
with some constant $c\in(0,\infty)$, uniformly in $\omega$ --- subject to the strong-ellipticity condition. For~$f\in L^1(\textd x)$ this yields
\begin{equation}
\epsilon^{-2}\bigl\langle\, f_\epsilon,\texte^{t\epsilon^{-2}\cmss L_\omega} f_\epsilon\rangle
\le \Vert f\Vert_1^2\,\frac c{t^{d/2}}.
\end{equation}
This is uniformly integrable in all dimensions $d\ge3$.

In dimension~$d=1,2$ one needs a corresponding bound on the \emph{gradient} of the heat kernel. Such bounds were proved in the annealed setting by Delmotte and Deuschel~\cite{Delmotte-Deuschel}. See Corollary~4.3 in~\cite{Biskup-Spohn} for details.
\end{proofsect}

The above conclusions permit a statement on the random Gaussian field introduced in Problem~\ref{pb19}. Indeed, let $(\phi_x)$ be a sample from the Gaussian measure with zero mean and covariance $(-\cmss L_\kappa)^{-1}$, for a collection of nearest-neighbor elliptic conductances~$\kappa$. Recall the notation $\phi_\epsilon(f)$ from \eqref{E:phi_epsilon}. Then we have:

\begin{corollary}
Suppose~$f$ is smooth with compact support and (in $d=1,2$) of zero total integral. As $\epsilon\downarrow0$, the law of~$\phi_\epsilon(f)$ tends to that of a Gaussian with mean zero and limiting variance
\begin{equation}
\label{E:6.29}
\text{\rm Var}\bigl(\phi_\epsilon(f)\bigr)
\,\,\underset{\epsilon\downarrow0}\longrightarrow\,\,
\bigl\langle\, f,(-Q)^{-1}f\bigr\rangle_{L^2(\R^d)}
\end{equation}
in~$\BbbP$-probability, where~$Q$ is the generator of the limiting Brownian motion.
\end{corollary}

\begin{proofsect}{Proof}
As $\phi_\epsilon$ is Gaussian, it suffices to prove the convergence of the variances, i.e., \eqref{E:6.29}. This is  \eqref{E:6.23} in disguise.
\end{proofsect}

\newcommand{\Hspace}{\HH}

The key point is that since the limit is non-random, the same will be true even if the law of the~$\phi$'s is further averaged over~$\kappa$. This permits the main conclusion of the paper of Biskup and Spohn~\cite{Biskup-Spohn} which repharse as follows:

\begin{theorem}[Scaling to GFF]
\label{thm-main}
Suppose~$V$ is as in \eqref{V-def} with~$\varrho$ compactly supported in~$(0,\infty)$. Let~$\mu$ be a gradient Gibbs measure for the potential~$V$ which we assume to be ergodic with respect to the translations of~$\Z^d$ and to have zero tilt. Then for every~$f\in \text{\rm Dom}((-\Delta)^{-1/2})$, the law of $\phi_\epsilon(f)$ tends to a Gaussian with mean zero and covariance
\begin{equation}
\label{limit}
\sigma_f^2:=\bigl\langle\, f,(-Q)^{-1}f\bigr\rangle_{L^2(\R^d)},
\end{equation}
where $Q^{-1}$ is the inverse of the operator \eqref{Q-op}.
\end{theorem}

\begin{proofsect}{Proof} (Sketch)
Conditional on the $\kappa$'s, the law of the $\phi$'s is Gaussian. Hence
\begin{equation}
\label{E:6.31}
\begin{aligned}
E_\mu(\texte^{\texti\phi_\epsilon(f)}\bigr)
&=E_\mu(E_\mu(\texte^{\texti\phi_\epsilon(f)}|\kappa)\bigr)
\\
&=E_\mu\bigl(\texte^{\texti E_\mu(\phi_\epsilon(f)|\kappa)-\frac12\Var(\phi_\epsilon(f)|\kappa)}\bigr).
\end{aligned}
\end{equation}
The above tells us that $\Var(\phi_\epsilon(f)|\kappa)\to\sigma_f^2$ in probability, but we cannot use it unless we can simultaneously deal with the conditional mean $E_\mu(\phi_\epsilon(f)|\kappa)$. The most substantive part of the result is the following representation: Suppose~$\mu$ has tilt~$t$. Then
\begin{equation}
E_\mu(\phi_x-\phi_y|\kappa)=t\cdot\bigl[\Psi(\kappa,x)-\Psi(\kappa,y)\bigr],
\end{equation}
where~$\Psi$ is the harmonic coordinate. This proves that the conditional mean is identically zero when~$t=0$. Then
\eqref{E:6.31} gives
\begin{equation}
E_\mu(\texte^{\texti\phi_\epsilon(f)}\bigr)
\,\,\underset{\epsilon\downarrow0}\longrightarrow\,\,
\texte^{-\frac12\sigma_f^2}.
\end{equation}
Since this holds for all multiples of~$f$ as well, Levy's characterization of convergence in law implies the desired claim.
\end{proofsect}

We note that the above Gaussian field with random (ergodic) covariance structure has been (probably first introduced and) studied by Caputo and Ioffe~\cite[Section~4.5]{Caputo-Ioffe}. Their motivation was to provide a link between the derivative of the exponential rate function for changing the tilt of the field --- the so called \emph{surface tension} --- and the diffusivity of the corresponding random walk among random conductances. For the above Gaussian case, this link is verified by a direct calculation, but for general uniformly-convex interactions --- for which one still has a random-walk representation (Naddaf and Spencer~\cite{Naddaf-Spencer}, Giacomin, Olla and Spohn~\cite{Giacomin-Olla-Spohn}) --- it remains conjectural despite serious effort.

\subsection{Random electric networks}
\noindent
Theorem~\ref{thm6.6} can be understood as an application of homogenization theory to electrostatic equilibrium. Indeed, given an assignment of charge $\rho(x)$ at vertex~$x$, we wish to find an electrostatic potential $\varphi_\omega\colon\Z^d\to\R$ satisfying the \emph{Poisson equation}
\begin{equation}
\label{E:Poisson}
\cmss L_\omega\varphi_\omega(x)=\rho(x),\qquad x\in\Z^d,
\end{equation}
with the normalization~$\varphi_\omega(0)=0$.
As to the existence of  solutions, we have:

\begin{lemma}
Suppose $\rho\in\text{\rm Dom}((-\cmss L_\omega)^{-\ffrac12})$ which is equivalent to
\begin{equation}
\sup_{\epsilon>0}\,\bigl\langle \rho,(\epsilon-\cmss L_\omega)^{-1}\rho\bigr\rangle_{\ell^2(\Z^d)}<\infty.
\end{equation}
Then $\{\varphi\colon\EE(\varphi)<\infty\}$ contains exactly one function $\varphi_\omega$ satisfying \eqref{E:Poisson} and $\varphi_\omega(0)=0$. Moreover, we have
\begin{equation}
\label{E:3.36}
\inf\bigl\{\EE(\varphi)+\langle\varphi,\rho\rangle\colon \EE(\varphi)<\infty\bigr\}=-\frac12\sup_{\epsilon>0}\,\bigl\langle \rho,(\epsilon-\cmss L_\omega)^{-1}\rho\bigr\rangle_{\ell^2(\Z^d)}
\end{equation}
and every minimizing sequence~$\varphi^{(n)}$ on the left has the property that $\EE(\varphi^{(n)}-\varphi_\omega)\to0$. In particular, the infimum is achieved by $\varphi_\omega$.
\end{lemma}

\begin{proofsect}{Proof} (Sketch)
The containment $\rho\in\text{\rm Dom}((-\cmss L_\omega)^{-\ffrac12})$ guarantees that the supremum in \eqref{E:3.36} equals $\langle \rho,(-\cmss L_\omega)^{-1}\rho\rangle_{\ell^2(\Z^d)}$. A completion of square yields
\begin{equation}
\label{E:complete-square}
\EE(\varphi)+\langle\varphi,\rho\rangle=\EE\bigl(\varphi-(-\cmss L_\omega)^{-1/2}\rho\bigr)
-\frac12\,\bigl\langle \rho,(-\cmss L_\omega)^{-1}\rho\bigr\rangle_{\ell^2(\Z^d)}
\end{equation}
which by the fact that the Dirichlet energy is non-negative implies that $\ge$ holds in \eqref{E:3.36}. Let thus~$\varphi^{(n)}$ be a minimizing sequence. The parallelogram law (see \eqref{E:3.25}) then immediately gives that $\EE(\varphi^{(n)}-\varphi^{(m)})\to0$ as $m,n\to\infty$. Moreover, $\EE(\varphi^{(n)})$ must remain bounded because if we had $\EE(\varphi^{(n)})\to\infty$, then $\langle\varphi^{(n)},\rho\rangle$ would tend to~$-\infty$ at the same rate as $\EE(\varphi^{(n)})$ tends to $+\infty$. But this is not possible as, by $\rho\in\text{\rm Dom}((-\cmss L_\omega)^{-\ffrac12})$ and Cauchy-Schwarz, $|\langle\varphi^{(n)},\rho\rangle|$ grows at most as the square root of $\EE(\varphi^{(n)})$.

Passing to $m\to\infty$ we thus construct a minimizer $\varphi_\omega$ on~$\Z^d$ with $\EE(\varphi_\omega)<\infty$. By adding small perturbations, we find that $\varphi_\omega$ solves \eqref{E:Poisson}. The identity \eqref{E:complete-square} and the fact that $\EE(\varphi)=0$ only for constants then shows that~$\varphi_\omega$ with a prescribed value  at one lattice site is unique.
\end{proofsect}

Having dismissed the questions of existence and uniqueness, let us now investigate what happens when we scale the lattice to have spacing~$\epsilon$ and scale the charge density to maintain a fixed macroscopic profile. As we will see, the following is  just a rewrite of results proved earlier:

\begin{theorem}
Suppose~$\omega$ is a sample from an ergodic measure on elliptic nearest-neighbor conductances and let~$Q$ denote generator of the (annealed) limiting Brownian motion for this environment.
Suppose~$f\colon\R^d\to\R$ is smooth and of compact support with $\int f(x)\textd x=0$ and let \begin{equation}
\rho_\epsilon(x):=\int_{[0,1]^d}f(\epsilon x+\epsilon z)\,\textd z.
\end{equation}
Define
\begin{equation}
\varphi_\omega^{(\epsilon)}(x):=\epsilon^{2} (\cmss L_\omega^{-1}\rho_\epsilon)(\lfloor x/\epsilon\rfloor),\qquad x\in\R^d.
\end{equation}
Then $\varphi_\omega^{(\epsilon)}\to\bar\varphi$ weakly in~$L^2(\textd x)$ in $\BbbP$-probability, where $\bar\varphi$ is the solution to the Poisson equation
\begin{equation}
Q\bar\varphi(x)=f(x),\qquad x\in\R^d.
\end{equation}
(This is well defined as $f\in\text{\rm Dom}((-Q)^{-1})$.
\end{theorem}

\begin{proofsect}{Proof}
Let $g\colon\R^d\to\R$. Then
\begin{equation}
\bigl\langle g,\varphi_\omega^{(\epsilon)}\bigr\rangle_{L^2(\R^d)}=
\sum_{x\in\Z^d}\int_{[0,1]^d}\textd z\,\epsilon^d \,g(\epsilon x+\epsilon z)\varphi^{(\epsilon)}(\epsilon x)
\end{equation}
where we used that $\varphi^{(\epsilon)}(\epsilon x+\epsilon z)=\varphi^{(\epsilon)}(\epsilon x)$ for all~$x\in\Z^d$ and all $z\in[0,1)^d$. Invoking the definition of $\rho_\epsilon$ and some elementary rewrites, we then get
\begin{equation}
\bigl\langle g,\varphi_\omega^{(\epsilon)}\bigr\rangle_{L^2(\R^d)}=\bigl\langle g_\epsilon,(\cmss L_\omega)^{-1}f_\epsilon\bigr\rangle_{L^2(\R^d)}.
\end{equation}
By Theorem~\ref{thm6.6}, the right-hand side tends to $\langle g,Q^{-1}f\rangle_{L^2(\R^d)}$ in $L^2(\BbbP)$ --- and thus in probability. It follows that $\varphi^{(\epsilon)}_\omega\to Q^{-1}f$ weakly in~$L^2(\textd x)$.
\end{proofsect}

\end{document}